\documentclass[12pt]{amsart}
\usepackage[margin=1in]{geometry}
\usepackage{enumerate}		
\usepackage{amssymb,amsfonts}
\usepackage{hyperref}
\usepackage{amsthm}
\usepackage{amsmath}
\usepackage{xcolor}

\newtheorem{theorem}{Theorem}[section]
\newtheorem*{theorem*}{Theorem}
\newtheorem{cor}[theorem]{Corollary}
\newtheorem*{cor*}{Corollary}

\newtheorem*{prop*}{Proposition}
\newtheorem{lemma}[theorem]{Lemma}
\newtheorem*{lemma*}{Lemma}

\newtheorem*{quest*}{Question}

\newtheorem*{fact*}{Fact}
\newtheorem{defn}[theorem]{Definition}
\newtheorem*{defn*}{Definition}

\theoremstyle{definition}
\newtheorem{rmk}[theorem]{Remark}
\newtheorem*{rmk*}{Remark}

\newtheorem*{exs*}{Examples}
\newtheorem{ex}[theorem]{Example}
\newtheorem*{ex*}{Example}

\newtheorem*{notn*}{Notation}

\newtheorem*{assumption*}{Assumption}

\DeclareMathOperator{\an}{an}
\DeclareMathOperator{\transexp}{transexp}

\DeclareMathOperator{\Supp}{Supp}
\DeclareMathOperator{\Lm}{Lm}
\DeclareMathOperator{\Lt}{Lt}

\DeclareMathOperator{\Lv}{Lv}

\DeclareMathOperator{\RR}{\mathbb{R}}
\DeclareMathOperator{\NN}{\mathbb{N}}

\DeclareMathOperator{\kk}{\boldsymbol{k}}

\DeclareMathOperator{\Init}{Init}

\title{Sublogarithmic-transexponential series}
\author[]{Adele Padgett}

\begin{document}

\begin{abstract}
    We adapt the construction of the field of logarithmic-exponential transseries of van den Dries, Macintyre, and Marker to build an ordered differential field of sublogarithmic-transexponential series. 
    We use this structure to 
    build a transexponential Hardy field closed under composition.
    Specifically, we prove that the germs at $+\infty$ of $\mathcal{L}_{\mathrm{transexp}}$-terms in a single variable are ordered, where $\mathcal{L}_{\transexp}$ is a language containing $\mathcal{L}_{\an}(\exp,\log)$ with new symbols for a transexponential function, its derivatives, and their compositional inverses.
\end{abstract}

\maketitle

\section{Introduction}

After Wilkie proved that $\RR_{\exp}$ is o-minimal \cite{WilkieRexp} and Miller determined that any o-minimal expansion of the real field is polynomially bounded or defines the exponential function \cite{MillerDichotomy}, van den Dries and Miller raised the question of whether an o-minimal structure can define a function that grows faster than any iterate of exponentiation \cite{FirstTransexpO-minimalQuestionReference}. 
We call such a function---and any structure that defines such a function---\textit{transexponential}.
Marker and Miller initially suggested that an expansion of a structure known to be o-minimal by a new function symbol that represents a real analytic function $f : (a, +\infty) \to \mathbb{R}$ satisfying
    $$f(x+1) = e^{f(x)}$$
might also be o-minimal, noting that any function satisfying this difference equation must be transexponential \cite{LeveledO-minimalSTructures}.
Kneser constructed a real analytic solution to this difference equation \cite{Kneser}, and Boshernitzan showed that there are solutions, including Kneser's, whose germs at $+\infty$ belong to a Hardy field \cite{Boshernitzan}.

The following is a two-part strategy for proving an expansion $\mathcal{R}$ of the real field is o-minimal \cite[Lemma 5.2]{vdDMM}:
\begin{enumerate}
    \item Find a language $\mathcal{L}$ in which $\mathrm{Th}(\mathcal{R})$ has quantifier elimination.
    
    \item Prove that the germs at $+\infty$ of $\mathcal{L}$-terms in a single variable with parameters in $\mathbb{R}$ are ordered.
\end{enumerate}
Let $E$ be Kneser's real-analytic solution to $f(x+1) = e^{f(x)}$.
To apply the above strategy to $\mathbb{R}_{\mathrm{an},\exp}(E)$, first consider Part (1). 
Just as $\mathbb{R}_{\mathrm{an},\exp}$ needs a symbol for log to have quantifier elimination, we will include a symbol $L = E^{-1}$ for the functional inverse of $E$.
Derivatives are definable with quantifiers using an order symbol, so we add symbols $E', E'', E''',\dots$, for the derivatives of $E$.
These are new function symbols, so we will add symbols for their functional inverses $(E')^{-1}, (E'')^{-1}, (E''')^{-1},\dots$ too.
We do not need new symbols for the derivatives of $L, (E')^{-1}, (E'')^{-1}, (E''')^{-1},\dots$ because we can express them in terms of the symbols we already have without quantifiers.
For example,
    $$L'(x) = \frac{1}{E'(L(x))}.$$
Thus, we suggest the following as a possible language in which to try to show $\mathrm{Th}\left(\mathbb{R}_{\mathrm{an},\exp}(E)\right)$ has quantifier elimination: 

\begin{defn}
    Let $\mathcal{L}_{\mathrm{transexp}} = \mathcal{L}_{\mathrm{an}}(\exp,\log) \cup \left\{E^{(d)} : d \in \mathbb{N}\right\} \cup \left\{\big(E^{(d)}\big)^{-1} : d \in \mathbb{N}\right\}$. 
\end{defn}

We intend to interpret $E^{(d)}$ as the $d$th derivative of $E$, for $d \in \mathbb{N}$, though we will also use the notation $E = E^{(0)}$, $E' = E^{(1)}$, $E'' = E^{(2)}$, $L = \left(E^{(0)}\right)^{-1}$ etc. when convenient.

The main result of the paper is to carry out Part (2) of the strategy outlined above for $\mathcal{L}_{\mathrm{transexp}}$.

\begin{theorem}\label{main theorem}
    The germs at $+\infty$ of $\mathcal{L}_{\mathrm{transexp}}$-terms in a single variable with parameters in $\mathbb{R}$ form a Hardy field closed under composition.
\end{theorem}

To prove Theorem \ref{main theorem}, we will build an ordered differential field of sublogarithmic-transex\-ponential series based on the construction of the field of logarithmic-exponential transseries in \cite{LEseries}.
The sublogarithmic-transexponential series field is built to embed the germs at $+\infty$ of $\mathcal{L}_{\mathrm{transexp}}$-terms as a substructure, so its order will induce an order on the germs.

The sublogarithmic-transexponential series satisfy a theory $T_{\mathrm{transexp}}$, defined as follows.
Let $T_{\mathrm{an}}(\exp,\log)$ be the universal axiomatization for $R_{\mathrm{an}}(\exp)$ given in \cite{vdDMM} in the language $\mathcal{L}_{\mathrm{an}}(\exp,\log)$.
Let $T_{\mathrm{transexp}}$ be the following set of universal axioms:
\begin{enumerate}
    \item Universal axioms for $T_{\mathrm{an}}(\exp,\log)$;
    
    \item Axioms identifying the restrictions of each of the new function symbols to $[0,1]$ with the corresponding restricted analytic function symbols, e.g.,
    \begin{enumerate}
        \item $0 \le x \le 1 \rightarrow E^{(d)}(x) = \widetilde{E^{(d)}}(x)$, where $\widetilde{E^{(d)}}$ is the restricted analytic function symbol corresponding to the function $E^{(d)}$;
    \end{enumerate}
    
    \item For all $x \ge 0$, $E(x+1) = \exp E(x)$; 
    
    \item For all $x<0$, $E^{(d)}(x) = 0$, for $d \in \mathbb{N}$;
    
    \item Axioms stating that the symbols for the inverse functions are indeed inverses, i.e.,
    \begin{enumerate}
        \item $\left(x \ge 1 \rightarrow E(L(x)) = x \right)\wedge \left(x < 1 \rightarrow L(x) = 0\right)$ and
        \item $\left(x \ge 1 \rightarrow E^{(d)}\left((E^{(d)})^{-1}(x)\right) = x \right) \wedge \left(x < 1 \rightarrow (E^{(d)})^{-1}(x) = 0\right)$;
    \end{enumerate}
    
    \item An axiom for each $d \in \mathbb{N}$ stating that $E^{(d+1)}$ is indeed equal to the derivative of $E^{(d)}$, i.e.,
    \begin{align*}
        \forall x > 0, \forall \epsilon > 0, \forall y > 0 \bigg(&0 < |x-y| < \min\left(1,\frac{\epsilon}{E(x+1)}\right) \rightarrow \\
            &\left|\frac{E^{(d)}(x) - E^{(d)}(y)}{x-y} - E^{(d+1)}(x)\right| < \epsilon \bigg)
    \end{align*}
    Each axiom in this schema is just the usual $\epsilon-\delta$ statement of the derivative, with $\delta := \min\left(1,\frac{\epsilon}{E(x+1)}\right)$.
\end{enumerate} 

Given $F \vDash T_{\mathrm{transexp}}$, we build a field $M_F \vDash T_{\mathrm{transexp}}$ of sublogarithmic-transexponential series with coefficients in $F$.
We also define a derivation on $M_F$ that works like differentiation with respect to the variable of the series field.
The difficulty is to build $M_F$ so that a unique ordering is maintained at each stage of the construction.

\subsection{Outline}

In Section 2, we present some basic computations involving the difference equation $E(x+1) = \exp E(x)$.
We also introduce the main challenge in defining an ordering on $M_F$, which arises from trying to compare asymptotic expressions involving $E$ and its derivatives.
We show how to rewrite such expressions in terms of an alternate ``basis" $E_0, E_1, E_2, \dots$ of functions so that a dominant term arises.
This idea first appears in \cite{Boshernitzan}.
We conclude Section 2 with some useful computations involving $E_0, E_1, E_2, \dots$.

The construction of $M_F$ is divided into Sections 3-5.
Note that all of the symbols used in Sections 3-5 are purely formal---we will use functional notation such as $E(x)^a$ and $\log E'(x)$ suggestively, but we do not assume they represent germs of functions. 
Instead, we will impose relations among these formal symbols through the course of the construction.

In Section 3, we introduce multiplicative groups $G_{X}$ and $\Lambda_{X}$ of monomials of the form $E^{(d_1)}(x_1)^{a_1} \cdots E^{(d_n)}(x_n)^{a_n}$ and $\log E'(y)^b$ respectively, with $n, d_1,\dots,d_n \in \NN$, $x_1,\dots,x_n,y \in X$, and $a_1,\dots,a_n,b \in \boldsymbol{k}$ for sets $X$ and $\kk$ satisfying relatively simple partial ordering assumptions detailed in Remark \ref{order props}.
Recall that $\log E'(x)$ is a purely formal symbol at this stage, but we include it here so that when we define an exponential function in Section 4, there is an element whose image we can define to be $E'(x)$.
The main result of Section 3 is Lemma \ref{rho_0 is injective}, in which we show the partial ordering assumptions on $X$ and $\kk$ determine a unique order on the group ring $\kk[G_{X}\Lambda_{X}]$.
The intuition behind the proof is that we can ``rewrite" each sum in terms of the alternate ``basis" $E_0, E_1,\dots$ to make it easy to see what its sign must be. 
We formalize the idea of ``rewriting" by defining an embedding from $\kk[G_{X}\Lambda_{X}]$ into a Hahn series field built from monomials of the form $E_{d_1}(x_1-1)^{a_1} \cdots E_{d_n}(x_n-1)^{a_n}$.
The computations in Section 2 give intuition for how the embedding must be defined.

In Section 4, we adapt the construction of the logarithmic-exponential series in \cite{LEseries}.
Instead of starting from powers of $x$, we would like to start with the monomials $G_{X}\Lambda_{X}$.
The structure resulting from the adapted construction is a model of $T_{\mathrm{an},\exp}(\log)$.
The key change from the original construction is to allow only \textit{finite} sums of certain kinds of monomials at a time.
This restriction allows us to use Lemma \ref{rho_0 is injective} to define an order if we assume the sets $X$ and $\kk$ used to define $G_{X}\Lambda_{X}$ satisfy the partial ordering assumptions from Remark \ref{order props}.
Section 4 is divided into three parts: First, we adapt the exponential part of the log-exp series construction to build countably many partial exponential rings.
Second, we adapt the logarithmic part of the log-exp series construction by defining embeddings between the partial exponential rings.
To define the embeddings, we again use intuition from the computations in Section 2 to ensure we impose the correct relations on these formal objects. 
The partial exponential rings and the embeddings between them form a directed system, so in the third part of this section, we take the direct limit of the system and show that it is a model of $T_{\an,\exp}(\log)$.

In Section 5, we build $M_F$ by starting with the usual logarithmic-exponential series and using the construction in Section 4 to gradually close off under $E$, its derivatives, and the inverse function symbols.
We can think of the $i$th stage of this construction as being closed under $i$ applications of $E$, its derivatives, and $L$ (but not necessarily the other inverses).
We build the $(i+1)$th stage from the $i$th stage as a direct limit of many Section-4-constructions, with inclusion maps between them. 
In each Section-4-construction, $X$ and $\kk$ are subsets of $H_i$ that satisfy the assumptions in Remark \ref{order props}---this will ensure the direct limit is closed under one more application of $E, E', E'', \dots$. 
We also include a base case to ensure the $(i+1)$th stage is closed under one more application of $L$. 
$M_F$ is the direct limit of the $i$th stage constructions for $i \in \NN$.
After constructing $M_F$, we check that it is closed under not only $E, E', E'',\dots$ and $L$ but also the other inverse functions $\left(E'\right)^{-1}, \left(E''\right)^{-1}, \dots$.
We then define a derivation $\partial : M_F \to M_F$ that works like differentiation with respect to the variable of the original log-exp series field, and which respects $E$.
Finally, we show that if $F = \RR$, then the Hardy field of germs at $+\infty$ of $\mathcal{L}_{\transexp}$-terms is order-isomorphic to a substructure of $M_{\RR}$.

\subsection{Background on transseries}

Transseries can be understood as a generalization of Laurent series that represent (often divergent) asymptotic expansions of real functions at $+\infty$.
They were initially studied independently by Dahn and G\"oring \cite{DahnGoringTransseries} in work on Tarski's problem on the real exponential field and Ecalle in work on the Dulac problem \cite{EcalleTransseries}.
We refer to \cite{LEseries} for the construction of the logarithmic-exponential (log-exp) series, which closely follows Dahn and G\"oring's original construction and which we follow closely in Section 4.
Elements of the field of log-exp series are infinite sums of ``log-exp monomials" arranged in decreasing order.
For example,
    $$2e^{e^x} - \frac{1}{2}xe^{2x} - x^{1/3} + 3\log x +1 + x^{-1} + x^{-2} + x^{-3} + \cdots  + xe^{-x}$$
is a log-exp series.

\subsubsection{Other transexponential transseries constructions}

Schmeling, in his 2001 PhD thesis \cite{Schmeling}, presents a variety of interesting results on transseries fields containing transfinite iterates $\exp_{\alpha}$ and $\log_{\alpha}$ of $\exp$ and $\log$ for $\alpha < \omega^{\omega}$.
For Schmeling, the $\omega$th iterate of $\exp$ satisfies the same difference equation 
    $$f(x+1) = \exp f(x)$$
that we take $E$ to satisfy.
Schmeling also gives a formal definition of what it means for a generalized power series field $C[[\mathfrak{M}]]$ to be a \textit{transseries} field and develops a theory of derivations and compositions for transseries fields and their transfinite exponential expansions. 

Van den Dries, van der Hoeven, and Kaplan extend Schmeling's thesis work in \cite{LogarithmicHyperseries} to build a field $\mathbb{L}$ of \textit{logarithmic hyperseries} that contains \textit{all} transfinite iterates of $\log$.
They also construct natural differentiation, integration, and composition operations on $\mathbb{L}$.

In \cite{HyperserialFields}, Bagayoko, van der Hoeven, and Kaplan generalize Schmeling's methods and exploit the properties of $\mathbb{L}$ to build a \textit{hyperserial field} $\mathbb{H}$ that contains all transfinite iterates of $\log$ \textit{and} $\exp$, and in \cite{HyperseriesDerivation} Bagayoko shows the derivation on $\mathbb{L}$ extends to the closure under hyperexponentials.

The hyperserial field $\mathbb{H}$ differs from the sublogarithmic-transexponential series constructed here in two ways:
First, in addition to being closed under all ordinal iterates of exp and log, $\mathbb{H}$ also allows many kinds of sums that are not in the field of sublogarithmic-transexponential series.
The notion of summability in the sublogarithmic-transexponential series is quite restrictive---in order to use the ordering result of Section 3 and to define the derivation in a natural way, certain kinds of finiteness are built in through the whole construction.
Second, it is not clear how to define an embedding of the field $\mathcal{H}$ of germs at $+\infty$ of unary $\mathcal{L}_{\transexp}$ terms directly into $\mathbb{H}$, while the sublogarithmic-transexponential series are built concretely to embed $\mathcal{H}$.
This embedding is used to show that $\mathcal{H}$ is ordered, from which Theorem \ref{main theorem} follows.

\section{Basic properties of \texorpdfstring{$E$}{E} and its derivatives}\label{section basic props of E}

In \cite{Kneser}, Kneser constructs a real analytic ``half-iterate" of $e^x$, i.e., a function $h$ such that $h(h(x)) = e^x$. 
Finding half-iterates of a function $\psi$ reduces to finding solutions to the Abel functional equation, which in its full generality is the following:
\begin{equation*}
	g(\psi(x)) =g(x)+c.
\end{equation*}
After constructing a real analytic solution $g$ to the Abel equation with $c = 1$ and $\psi=\exp$, Kneser defines $h(x) = g^{-1}\left(g(x)+\frac{1}{2}\right)$, so that
\begin{align*}
    h(h(x)) &= g^{-1}\left(g\left(g^{-1}\left(g(x)+\frac{1}{2}\right)\right)+\frac{1}{2}\right) \\
        &= g^{-1}\left(g(x)+1\right) \\
        &= g^{-1}(g(\exp x)) \\
        &=\exp x.
\end{align*}

However, we are concerned not with partial iterates of $\exp$, but with Kneser's solution $g$ to the Abel equation with $c = 1$ and $\psi=\exp$. 
If $g$ satisfies $g(\exp(x))=g(x)+1$, then $g^{-1}$ satisfies 
\begin{equation*}
	g^{-1}(x+1)=\exp(g^{-1}(x)).
\end{equation*} 
Any such $g^{-1}$ will eventually dominate any iterate of $\exp$.
We call such a function \textit{transexponential}.
Kneser's construction gives the following:
\begin{theorem}[{{\cite{Kneser}}}]
    The functional equation $g(\exp x) = g(x) + 1$ has a real analytic solution on $x>0$.
\end{theorem}
This functional equation actually has infinitely many real-analytic solutions, but we will use Kneser's, which may be taken to satisfy $g(1) = 0$.
We will call the solution given by Kneser's construction $L$.
We will call its compositional inverse $E$.
We will also use $E$ and $L$ to refer to the corresponding operators on germs of functions at $+\infty$ and as formal symbols.

We now give several simple identities and relations derived from the functional equation $E$ satisfies.
\begin{rmk}
    Whenever we write a relation, it should be understood that it holds for all large enough $x$.
\end{rmk}

\begin{enumerate}
    \item $E(x+1) > E(x)^a$ for any $a \in \mathbb{R}$, and $\displaystyle \lim_{x \to \infty} \frac{E(x)^a}{E(x+1)} = 0$.
    
    \item $E'(x+1) = e^{E(x)}E'(x) = E(x+1)E'(x)$ by the chain rule, so that $E'(x) > E(x)$. 
    Again $\displaystyle \lim_{x \to \infty} \frac{E(x)}{E'(x)} = 0$. 
	This also shows that $E^{(d_1)}(x) > E^{(d_2)}(x)$ if $d_1 > d_2$.
	
	\item Repeatedly differentiating both sides of $E(x+1) = \exp E(x)$ gives algebraic-difference-differential equations of all orders, which are of the form
	    $$E^{(d)}(x+1) = E(x+1)\sum_{k=1}^d\sum_{\bar{j}} \frac{d!}{j_1! \cdots j_{d-k+1}!}\left(\frac{E'(x)}{1!}\right)^{j_1} \cdots \left(\frac{E^{(d-k+1)}(x)}{(d-k+1)!}\right)^{j_{d-k+1}}$$
	where the second sum is taken over all sequences $j_1,\dots,j_{d-k+1}$ in $\NN$ satisfying
    \begin{align*}
        &j_1 + j_2 + \cdots + j_{d-k+1} = k \\
        &j_1 + 2j_2 + \cdots + (d-k+1)j_{d-k+1} = d.
    \end{align*}
	The double sum in the above equation is the formula for the $d$th complete Bell polynomial with arguments $E'(x),\dots, E^{(d)}(x)$, which we will denote by $B_d(x)$, i.e.,
        $$B_d(x) := \sum_{k=1}^d\sum_{\bar{j}} \frac{d!}{j_1! \cdots j_{d-k+1}!}\left(\frac{E'(x)}{1!}\right)^{j_1} \cdots \left(\frac{E^{(d-k+1)}(x)}{(d-k+1)!}\right)^{j_{d-k+1}}.$$
    Bell polynomials are used in the study of set partitions, though they arise here as an instance of Fa\`a di Bruno's formula, which generalizes the chain rule by computing higher derivatives of a composition of functions.
\end{enumerate}

The following result will help us derive more identities involving $E$ and its derivatives.

\begin{lemma}[{{\cite[Lemma 3.5]{Boshernitzan}}}] \label{B}
	Let $g(x), h(x) > 0$ be continuous, $\displaystyle \lim_{x \to \infty} h(x)=+\infty$, and assume that
	\begin{align*}
		h(x+1)&> 2h(x)  \\
		|g(x+1)-g(x)| &\le h(x+1) 
	\end{align*}
	for all large enough $x$.
	Then $|g(x)| < 3h(x)$ for all large enough $x$. 
\end{lemma}

\begin{proof}
	Using the second inequality $n$ times, we have 
	\begin{align*}
		g(x+n) &\le g(x+(n-1)) + h(x+n) \\
			&\le \left[g(x+(n-2)) + h(x+(n-1))\right] + h(x+n) \\
			&\vdots \\
			&\le  g(x)+ \sum_{k=0}^{n-1} h(x+(n-k)).
	\end{align*}
	Now using the first inequality $n$ times, $h(x+n) > 2h(x+(n-1)) > \dots > 2^n h(x)$, so $\displaystyle \frac{h(x)}{h(x+n)} < \frac{1}{2^n}$. 
	So we have 
	\begin{align*}
		\lim_{n \to \infty} \frac{g(x+n)}{h(x+n)} &\le \lim_{n \to \infty} \frac{g(x) + \sum_{k=0}^{n-1} h(x+(n-k))}{h(x+n)} \\
			&\le \lim_{n \to \infty} \left (\sum_{k=0}^{n-1} \frac{1}{2^k} + \frac{g(x)}{h(x+n)} \right ) \\
			&<3. \qedhere
	\end{align*}
\end{proof}

If we take $g(x) = \log(E'(x))$ and $h(x) = \log(E(x))$ then we can apply Lemma \ref{B} because:
\begin{align*}
	h(x+1) &= \log(E(x+1)) = \log(e^{E(x)}) = E(x)= e^{E(x-1)} \\
	    &> 2E(x-1) = 2h(x) \\
	|g(x+1)-g(x)| &= \log(E'(x+1)) - \log(E'(x)) \\
	    &= \log\left(\frac{E'(x+1)}{E'(x)}\right) \\
	    &= \log(E(x+1)) = h(x+1).
\end{align*}
Lemma \ref{B} gives that $\log(E'(x)) < 3\log(E(x))$, and exponentiating, we get $E'(x) < E(x)^3$.
From here, we can make a sequence of comparisons. 

\begin{lemma}\label{comps}
For all large enough $x$, we have
\begin{enumerate}
	\item $E^{(d)}(x)<E(x)E(x-1)^{3d}$
	\item $E^{(d)}(x)<E(x)^2$
	\item $E^{(d)}(x-1)^a<E(x)$ for all $a \in \mathbb{R}$
	\item $E^{(d)}(x-r)^a<E(x)$ for all $a \in \mathbb{R}$ and $r>0$
\end{enumerate}
\end{lemma}

\begin{proof}
	 The first part is proven in \cite{Boshernitzan}, using Lemma \ref{B}.
	 The second part follows from the first part and Lemma \ref{B}:
	 \begin{align*}
			 \lim_{x \to \infty} \frac{E^{(d)}(x+1)}{E(x+1)^2} 
			 	&= \lim_{x \to \infty}\frac{E(x+1) E(x)^{3d}}{E(x+1)^2} \\
			 	& \le \lim_{x \to \infty} \frac{E(x)^{3d}}{E(x+1)} \\
				&= \lim_{x \to \infty} \frac{E(x)^{3d}}{e^{E(x)}} \\
				&= 0.
	\end{align*}
	
	The third part uses Lemma \ref{B} and the fact that $\displaystyle \lim_{x \to \infty} \frac{E(x)^a}{E(x+1)} = 0$ for any $a \in \mathbb{R}$:
	\begin{align*}
		\lim_{x \to \infty}\frac{E^{(d)}(x-1)^a}{E(x)} 
		    &\le \lim_{x \to \infty} \frac{(E(x-1)E(x-2)^{3d})^a}{E(x)} \\
		    &\le \lim_{x \to \infty} \frac{E(x-1)^{a+1}}{E(x)} \\
		    &= 0.
	\end{align*}
	
	For the fourth part, we need only show that $\displaystyle \lim_{x \to \infty} \frac{E(x-r)^m}{E(x)} = 0$ for any $m \in \mathbb{N}$, since the second part tells us that $\displaystyle \lim_{x \to \infty} \frac{E^{(d)}(x-r)^a}{(E(x-r)^2)^a} = 0$.
	The proof is the same as for the third part, using the fact that partial iterates of the exponential also have faster-than-polynomial growth.
	First assume $r$ is rational.
	Write $e_{1/q}$ to denote a $1/q$th iterate of $\exp$ (meaning $\exp = \overbrace{e_{1/q} \circ \cdots \circ e_{1/q}}^{q \text{ times}}$)  and write $e_{p/q}$ to denote $p$ iterates of $e_{1/q}$, i.e., $e_{p/q} = \overbrace{e_{1/q} \circ \cdots \circ e_{1/q}}^{p \text{ times}}$.
	Then
	\begin{align*}
		\lim_{x \to \infty} \frac{E(x-r)^m}{E(x)} 
			& = \lim_{x \to \infty} \frac{e_{1-r}(E(x-1))^m}{e^{E(x-1)}} \\
			&= \lim_{y \to \infty} \frac{y^m}{e_r(y)} \\
			&\le \lim_{y \to \infty} \frac{y^m}{y^a} \hspace{1cm} \text{for any $a$} \\
			&=0 \hspace{2cm} \text{for $a>m$.}
	\end{align*}
	If $r$ is irrational, then there is a rational $r'$ such that $0 < r' < r$, so 
	\begin{align*}
	    \lim_{x \to \infty} \frac{E(x-r)^m}{E(x)} &\le \lim_{x \to \infty} \frac{E(x-r')^m}{E(x)} = 0.  \qedhere
	\end{align*}
\end{proof}

\subsection{Monomials with equivalent growth rates and logarithmic derivatives} \label{Boshernitzan's sequence}

Note that what is actually proven in each part of Lemma \ref{comps} is that the limit of the quotient of the smaller expression by the larger expression is 0. This means that in each of these comparisons, the larger expression is not only larger, but also has a faster growth rate as $x \to \infty$. 

However, it is possible for the quotient of simple expressions involving $E$ and its derivatives to approach 1.
For example,
\begin{align*}
    \lim_{x \to \infty} \frac{E(x)E''(x)}{E'(x)^2} 
        &= \lim_{x \to \infty}\frac{E(x)^2\big(E'(x-1)^2 + E''(x-1)\big)}{E(x)^2E'(x-1)^2} \\
        &= \lim_{x \to \infty} \left(1 + \frac{E''(x-1)}{E'(x-1)^2}\right) \\
        &= 1.
\end{align*}
by Lemma \ref{comps}, so $E(x)E''(x)$ and $E'(x)^2$ grow at roughly the same rate as $x$ approaches infinity.

Even though $E(x)E''(x)$ and $E'(x)^2$ have the same growth rate, it is easy to determine that $E(x)E''(x) > E'(x)^2:$ 
\begin{align*}
    E(x)E''(x) - E'(x)^2 &= E(x)^2\big(E'(x-1)^2 + E''(x-1)\big) - E(x)^2E'(x-1)^2 \\
        &= E(x)^2E''(x-1) \\
        &> 0.
\end{align*}
However, it can be quite complicated to compute the sign of an expression involving many terms with equivalent growth rates.
To remedy this, it will be useful to be able to rewrite expressions involving $E$ and its derivatives in terms of different functions that all have distinct growth rates.
To this end, we follow \cite{Boshernitzan} and introduce a useful sequence of functions.
\begin{defn}
    Let $E_0(x) = E(x)$ and let $E_{d+1}(x) = \frac{E_d'(x)}{E_d(x)}$ for $d \in \mathbb{N}$.
\end{defn}

We would now like to show that $E_0, E_1, E_2,\dots$ all have distinct growth rates. 
We will obtain this result using slightly different computations than \cite{Boshernitzan} in order to provide intuition for the purely formal construction in Sections 3-5.

\begin{lemma}\label{log derivative asymptotics}
    For  $d \ge 2$, we have $E_d(x) = E'(x-d) + R_d(x)$ where $R_d(x)$ is such that $\lim_{x \to \infty} \frac{R_d(x)}{E'(x-d-1)} = 1$.
\end{lemma}

\begin{proof}
    We know $E_0(x) = E(x)$ and $E_1(x) = \frac{E'(x)}{E(x)} = E'(x-1)$.
    Then, 
        $$E_2(x) = \frac{E''(x-1)}{E'(x-1)} = \frac{E'(x-2)^2 +E''(x-2)}{E'(x-2)} = E'(x-2) + \frac{E''(x-2)}{E'(x-2)}.$$
    Let $R_2(x) := \frac{E''(x-2)}{E'(x-2)}$, which is real analytic on $x>0$.
    The same calculation as above, shifted down by 1, shows that $R_2(x) = E'(x-3) + \frac{E''(x-3)}{E'(x-3)}$.
    So $\lim_{x \to \infty} \frac{R_2(x)}{E'(x-3)} = 1$.
    
    We proceed by induction.
    Suppose 
    \begin{enumerate}
        \item $R_{d}(x)$ has been defined so that $\lim_{x \to \infty} \frac{R_{d}(x)}{E'(x-d-1)} = 1$ and $R_{d}$ is real analytic at all large enough values of $x$.
        
        \item $E_{d}(x) = E'(x-d) + R_{d}(x)$.
    \end{enumerate}
    Then 
    \begin{align*}
        E_{d+1}(x) =& \frac{E''(x-d) + R_{d}'(x)}{E'(x-d) + R_{d}(x)} \\
            =& \frac{E(x-d)\big(E'(x-d-1)^2 + E''(x-d-1)\big) + R_{d}'(x)}{E'(x-d) + R_{d}(x)} \\
            &+ \frac{R_{d}(x)E'(x-d-1)-R_{d}(x)E'(x-d-1)}{E'(x-d) + R_{d}(x)}\\
            =& \frac{E(x-d)E'(x-d-1)^2 + R_{d}(x)E'(x-d-1)}{E(x-d)E'(x-d-1) + R_{d}(x)} \\
            &+ \frac{E(x-d)E''(x-d-1) + R_{d}'(x) - R_{d}(x)E'(x-d-1)}{E'(x-d) + R_{d}(x)} \\
            = &E'(x-d-1) + \frac{E(x-d)E''(x-d-1)-R_{d}(x)E'(x-d-1)+R_{d}'(x)}{E'(x-d) + R_{d}(x)}.
    \end{align*}
    Now we would like to define $R_{d+1}(x)$ to be
        $$\frac{E(x-d)E''(x-d-1)-R_{d}(x)E'(x-d-1)+R_{d}'(x)}{E'(x-d) + R_{d}(x)}$$
    so we must show this expression satisfies the induction hypothesis.
    The expression is real analytic at large $x$ because the denominator is large when $x$ is large.
    So it remains to show the expression has the same growth rate as $E'(x-d-2)$.
    We will show $\frac{E(x-d)E''(x-d-1)}{E'(x-d) + R_{d}(x)}$ has the same growth rate at $E'(x-d-2)$, and $\frac{-R_{d}(x)E'(x-d-1)+R_{d}'(x)}{E'(x-d) + R_{d}(x)}$ approaches 0 as $x \to \infty$.
    
    By induction, we know $\displaystyle\lim_{x \to \infty} \frac{R_{d}(x)}{E'(x-d-1)} = 1$.
    Since $R_d$ is differentiable at large $x$, $\displaystyle \lim_{x \to \infty} \frac{R_{d}'(x)}{E''(x-d-1)} = 1$.
    So
    \begin{align*}
        \lim_{x \to \infty} \frac{-R_{d}(x)E'(x-d-1)+R_{d}'(x)}{E'(x-d) + R_{d}(x)}
            = \lim_{x \to \infty} \frac{-R_{d}(x)E'(x-d-1)+R_{d}'(x)}{E(x-d)\left(E'(x-d-1) + \frac{R_{d}(x)}{E(x-d)}\right)}
            = 0
    \end{align*}
    since $E'(x-d-1)^n < E(x-d)$ for all $n \in \mathbb{N}$ by Lemma \ref{comps}.
    
    Now we can rewrite 
    \begin{align*}
        \frac{E(x-d)E''(x-d-1)}{E(x-d)E'(x-d-1) + R_{d}(x)} 
            &= \frac{E(x-d-1)\big(E'(x-d-2)^2 + E''(x-d-2)}{E(x-d-1)E'(x-d-2) + \frac{R_{d}(x)}{E(x-d)}} \\
            &= \frac{E'(x-d-2)^2 + E''(x-d-2)}{E'(x-d-2) + \frac{R_{d}(x)}{E(x-d)E(x-d-1)}}.
    \end{align*}
    So again by Lemma \ref{comps} 
    \begin{align*}
        \lim_{x \to \infty} \frac{\frac{E'(x-d-2)^2 + E''(x-d-2)}{E'(x-d-2) + \frac{R_{d}(x)}{E(x-d)E(x-d-1)}}}{E'(x-d-2)} &= \lim_{x \to \infty} \frac{E'(x-d-2)^2 + E''(x-d-2)}{E'(x-d-2)^2} = 1. \qedhere
    \end{align*}
\end{proof}

\begin{cor}
    For all large enough $x$ and all $d,n \in \mathbb{N}$, we have $E_d(x) > E_{d+1}(x)^n$.
\end{cor}
This corollary follows immediately from Lemmas \ref{comps} and \ref{log derivative asymptotics}, and it is also proved in \cite{Boshernitzan} via different computations.

For each $d \in \mathbb{N}$, we can rewrite $E^{(d)}$ as a polynomial with integer coefficients in $E_0,\dots, E_d$.
For example,
\begin{align*}
    E' &= E_0 \cdot E_1 \\
    E'' &= E_0 \cdot (E_1)^2 + E_0\cdot E_1 \cdot E_2.
\end{align*}
We get similar expressions for all derivatives of $E$ because 
\begin{align*}
    E_d'=E_dE_{d+1}
\end{align*}
for all $d \in \mathbb{N}$.

The importance of the above corollary is that it allows us to rewrite expressions involving the derivatives of $E$ in terms of $E_0, E_1, E_2,\dots$ to identify a dominant term.
For example, we showed earlier that $E(x)E''(x)$ and $E'(x)^2$ grow at the same rate as $x \to \infty$.
This is easy to see when we rewrite these expressions in terms of $E_0, E_1, E_2,\dots$
\begin{align*}
    E(x)E''(x) &= E_0(x)^2 E_1(x)^2 + E_0(x)^2 E_1(x) E_2(x) \\
    E'(x)^2 &= E_0(x)^2E_1(x)^2.
\end{align*}

We can also rewrite $E_d(x)$ for $d > 0$ in terms of $E_0(x-1),\dots,E_d(x-1)$.
For example, 
\begin{align*}
    E_1(x) &= E'(x-1) = E_0(x-1) E_1(x-1) \\
    E_2(x) &= \frac{E_1'(x)}{E_1(x)} \\
        &= \frac{E_0'(x-1)E_1(x-1) + E_0(x-1)E_1'(x-1)}{E_0(x-1)E_1(x-1)} \\
        &= \frac{E_0(x-1)E_1(x-1)^2 + E_0(x-1)E_1(x-1)E_2(x-1)}{E_0(x-1)E_1(x-1)} \\
        &= E_1(x-1)\left(1 + \frac{E_2(x-1)}{E_1(x-1)}\right).
\end{align*}
The computation for $d=3$ below illustrates the general pattern:
\begin{align*}
    E_3(x) &= \frac{E_2'(x)}{E_2(x)} \\
        &= \frac{E_1'(x-1) + E_2'(x-1)}{E_1(x-1) + E_2(x-1)} \\
        &= \frac{E_1(x-1)E_2(x-1) + E_2(x-1)E_3(x-1)}{E_1(x-1) + E_2(x-1)} \\
        &= E_2(x-1)\left(1 + \frac{-E_2(x-1) + E_3(x-1)}{E_1(x-1) + E_2(x-1)}\right).
\end{align*}
We now show how the computation for $d=3$ can be used to obtain similar expressions for all $d \ge 3$.

\begin{lemma}\label{log derivative shift}
    Let $\epsilon_2(x) = \frac{E_2(x-1)}{E_1(x-1)}$.
    For all $d \ge 2$, $E_{d+1}(x) = E_{d}(x-1)(1 + \epsilon_{d+1}(x))$, where $\epsilon_{d+1}(x) = \frac{\epsilon_d'(x)}{E_d(x-1)}\cdot \frac{1}{1+\epsilon_d(x)}$ and $\epsilon_{d+1}(x)$ is expressed in terms of $E_1(x-1),\dots,E_{d+1}(x-1)$.
\end{lemma}

\begin{proof}
    We will use induction.
    In the base case, 
    \begin{align*}
        \epsilon_3(x)    
            &=\frac{\left(\frac{E_2(x-1)}{E_1(x-1)}\right)'}{E_2(x-1)\left(1+\frac{E_2(x-1)}{E_1(x-1)}\right)} \\
            &= \frac{\frac{E_1(x-1)E_2(x-1)E_3(x-1) - E_1(x-1)E_2(x-1)^2}{E_1(x-1)^2}}{E_2(x-1)\left(1+\frac{E_2(x-1)}{E_1(x-1)}\right)} \\
            &=\frac{E_3(x-1)-E_2(x-1)}{E_1(x-1) + E_2(x-1)}
    \end{align*}
    and our earlier computation shows $E_3(x) = E_2(x)(1 + \epsilon_3(x))$.
    
    Suppose the hypotheses hold for $n=2,\dots,d$.
    Then
    \begin{align*}
        E_{d+1} &= \frac{E_d'(x)}{E_d(x)} \\
        &= \frac{E_{d-1}'(x-1)(1+\epsilon_d(x)) + E_{d-1}(x-1)\epsilon_d'(x)}{E_{d-1}(x-1)(1 + \epsilon_d(x))} \\
        &= \frac{E_{d-1}'(x-1)}{E_{d-1}(x-1)} + \frac{\epsilon_{d}'(x)}{1 + \epsilon_d(x)} \\
        &= E_d(x-1)\left(1 + \frac{\epsilon_{d}'(x)}{E_d(x-1)}\cdot\frac{1}{1 + \epsilon_d(x)}\right) \\
        &= E_d(x-1)(1 + \epsilon_{d+1}(x)).
    \end{align*}
    By induction, $\epsilon_{d}(x)$ can be expressed in terms of $E_1(x-1),\dots, E_{d}(x-1)$.
    Since for all $n \in \mathbb{N}$ we have 
        $$E_n'(x-1) = E_n(x-1)E_{n+1}(x-1)$$
    we can express $\epsilon_d'(x)$ in terms of $E_1(x-1),\dots, E_{d+1}(x-1)$.
\end{proof}

\begin{cor}\label{problematic terms approximation}
    For all $d \ge 1$, $\lim_{x \to \infty} \frac{\left(\frac{E_{d-1}(x-1)}{E_{d}(x-1)}\right)}{E(x-d)} = 1$.
\end{cor}

\begin{proof}
    By Lemma \ref{log derivative asymptotics}
    \begin{align*}
        \frac{E_{d-1}(x-1)}{E_{d}(x-1)} &= \frac{E'(x-d) + R_{d-1}(x-1)}{E'(x-d-1) + R_{d}(x-1)} \\
            &= \frac{E'(x-d) + R_{d-1}(x-1)}{E'(x-d-1)}\left(\frac{1}{1 + \frac{R_{d}(x-1)}{E'(x-d-1)}}\right).
    \end{align*}
    Also by Lemma \ref{log derivative asymptotics}, we know
    \begin{align*}
        \lim_{x \to \infty} \frac{R_{d-1}(x-1)}{E'(x-d-1)} = 1 \\
        \lim_{x \to \infty} \frac{R_{d}(x-1)}{E'(x-d-1)} = 0.
    \end{align*}
    So $\lim_{x \to \infty} \frac{E_{d-1}(x-1)}{E_{d}(x-1)} = E(x-d) + 1$ and the result follows.
\end{proof}

We will use these results, in which $E_0,E_1,\dots$ are functions, to give intuition for how to order purely formal sums of monomials in the remaining sections.

\section{Ordering the \texorpdfstring{$E$}{E}-sums}\label{section 3}

In the previous section, we used the difference equation for $E$ and the difference-differential equations for its derivatives to derive some simple relations among these functions. 
We also introduced the sequence of logarithmic derivatives $E_0, E_1, \dots$ as a different ``basis" for expressions involving $E$ and its derivatives, which allowed us to identify a dominant monomial in sums of asymptotic monomials.
In this section, we introduce the formal notions of $E$-sums and $E_*$-sums and impose relations upon them to match the properties derived in Section \ref{section basic props of E}. 
We also show that under some simple hypotheses, the $E$-sums can be uniquely ordered.

In the following definitions, we use the same symbols we used for functions in the previous section (e.g., $E^{(d)}(x), E_d(x)$, etc.) to represent \textit{purely formal} objects.
We use these symbols to highlight the connection to results of the previous section and for consistency of notation.

\begin{defn}[$E,m$-sums]
    Let $X$ be a set of variables, and let $\boldsymbol{k}$ be a field.
    Fix $m \in \mathbb{Z}$.
    \begin{enumerate}
        \item First, let $G_{X-m}$ be the multiplicative abelian group generated by expressions of the form $E^{(d)}(x-m)^a$ 
        for $d \in \mathbb{N}$, $x \in X$, and $a \in \boldsymbol{k}$, where we identify expressions of the form $E^{(d)}(x - m)^aE^{(d)}(x - m)^b$ with $E^{(d)}(x - m)^{a+b}$. 
        
        \item Second, let $\Lambda_{X-m}$ be the multiplicative abelian group generated by expressions of the form $\big(\log E'(x-m)\big)^a$ for $x \in X$ and $a \in \boldsymbol{k}$, where we identify expressions of the form $\big(\log E'(x-m)\big)^a \big(\log E'(x-m)\big)^b$ with $\big(\log E'(x-m)\big)^{a+b}$.
    \end{enumerate}
    We will call the group ring $\boldsymbol{k}[G_{X-m}]$ the $E,m$-sums.
    We will call the group ring $\boldsymbol{k}[G_{X-m}\Lambda_{X-m}]$ the log-$E,m$-sums.
\end{defn}

We will write $\log E'(x-m)^a$ instead of $\big(\log E'(x-m)\big)^a$ to avoid using too many parentheses. 
The parentheses should be implicitly understood as surrounding $\log E'(x-m)$.
If we intend an exponent to apply only to $E'(x-m)$, we will instead write $\log \big(E'(x-m)^a\big)$.

\begin{defn}[$E_*,m$-sums]
    Let $X$ be an ordered set of variables, and let $\boldsymbol{k}$ be an ordered field.
    Fix $m \in \mathbb{Z}$.
    Let $H_{X-m}$ be the multiplicative abelian group generated by expressions of the form $E_d(x-m)^a$ and $E_0(x-m+k)^a$ for $d,k \in \mathbb{N}$, $x \in X$, and $a \in \boldsymbol{k}$, where we identify expressions of the form $E_d(x-m)^aE_d(x-m)^b$ with $E_d(x-m)^{a+b}$ and $E_0(x-m+k)^aE_0(x-m+k)^b$ with $E_0(x-m+k)^{a+b}$.
    
    Define an order on $H_{X-m}$ as follows: We can write any $h \in H_{X-m}$ as
        $$h = \prod_{j=1}^p E_0(x_j-m+k)^{\alpha_{j,k}} \cdots E_0(x_j-m+0)^{\alpha_{j,0}} E_1(x_j-m)^{\beta_j}E_2(x_j-m)^{a_{j,2}}\cdots E_d(x_j-m)^{a_{j,d}}$$
    for some $p,k,d \in \mathbb{N}$ and $x_1 > \cdots > x_p$.
    Let 
    \begin{enumerate}
        \item $\overline{\alpha_{l}} = (\alpha_{1,l}, \alpha_{2,l},\dots,\alpha_{p,l})$ for $l = 0,\dots,k$
        
        \item $\overline{\beta} = (\beta_{1}, \beta_{2},\dots,\beta_{p})$
        
        \item $\overline{a_{l}} = (a_{1,l}, a_{2,l},\dots,a_{p,l})$ for $l = 2,\dots,d$
    \end{enumerate}
    Define $h > 1$ if and only if the first nonzero element of the following sequence is positive:
        $$\overline{\alpha_{k}}, \overline{\alpha_{k-1}}, \dots, \overline{\alpha_{0}}, \overline{\beta}, \overline{a_{2}}, \dots, \overline{a_{d}}.$$
        
    Let the $E_*,m$-sums be the Hahn series field $\boldsymbol{k}((H_{X-m}))$.
\end{defn}

At this point, it is not yet clear how we should define a derivation on $\boldsymbol{k}((H_{X-m}))$, but we can define a derivation on the subfield of elements whose monomials have only integer exponents.

\begin{defn}\label{Boshernitzan's sequence with integer powers}
    Let $X$ be an ordered set of variables, and let $\boldsymbol{k}$ be an ordered differential field.
    Let $H_{X-m}' \subset H_{X-m}$ be the subgroup generated by expressions of the form $E_d(x-m)^n$ for $d \in \mathbb{N}$, $x \in X$, $n \in \mathbb{Z}$.
    Then $\boldsymbol{k}((H_{X-m}'))$ is a subfield of $\boldsymbol{k}((H_{X-m}))$.
Define a derivation on $\boldsymbol{k}((H_{X-m}'))$ by
    $$\partial_m(E_d(x-m)) = E_d(x-m)E_{d+1}(x-m)$$
which comes from the definition $E_{d+1} = \frac{E_d'}{E_d}$ of the sequence of logarithmic derivatives from subsection \ref{Boshernitzan's sequence}.
\end{defn}

In this section, we will only discuss $\boldsymbol{k}[G_{X-m}]$ and $\boldsymbol{k}[G_{X-m}\Lambda_{X-m}]$ with $m = 0$, and $\boldsymbol{k}((H_{X-m}))$ with $m = 0,1$, though the proofs work the same for any consecutive pair $m = n,n+1 \in \mathbb{Z}$.

\begin{rmk}
The goal of this section is to show that the log-$E$-sums $\boldsymbol{k}[G_{X}\Lambda_{X}]$ can be ordered in a way that is consistent with how we intend to interpret them as germs of functions. 
We will first order $\boldsymbol{k}[G_{X}]$ by defining an injective ring homomorphism
    $$\sigma_0 : \boldsymbol{k}[G_{X}] \to \boldsymbol{k}((H_{X})).$$
Since $\boldsymbol{k}((H_{X}))$ is an ordered Hahn series field, $\sigma_0$ will induce an order on $\boldsymbol{k}[G_{X}]$.
We will then define an order-preserving embedding $\nu_0 : \sigma_0(\boldsymbol{k}[G_{X}]) \to \boldsymbol{k}((H_{X-1}))$ and use it to define an ``approximation map" 
    $$\rho_0 : \boldsymbol{k}[G_{X}\Lambda_{X}] \to \boldsymbol{k}((H_{X-1}))$$
so-called because, intuitively, it ``approximates" log-$E$-sums by elements of $\boldsymbol{k}((H_{X-1}))$ well enough to induce a unique ordering of $\boldsymbol{k}[G_{X}\Lambda_{X}]$ compatible with the order induced on $\boldsymbol{k}[G_{X}]$.
\end{rmk}

\subsection{The order on \texorpdfstring{$\boldsymbol{k}((H_{X}))$}{k((H0))} induces an order on \texorpdfstring{$\boldsymbol{k}[G_{X}]$}{k[G0]}}

We will define a homomorphism $\sigma_0 : \boldsymbol{k}[G_{X}] \to \boldsymbol{k}((H_{X}))$. 
Intuitively, we want $\sigma_0(s)$ to be $s \in \boldsymbol{k}[G_{X}]$ ``rewritten" using the expressions for the derivatives of $E$ as polynomials in $E_0, E_1, \dots$ from Section \ref{section basic props of E}.
We can begin by defining $\sigma_0(E^{(d)}(x))$ exactly according to this intuition:
    $$\sigma_0(E^{(d)}(x)) := E_0(x)E_1(x)^d + \cdots + E_0(x)\cdots E_d(x)$$
is the polynomial expression for $E^{(d)}(x)$ in terms of the logarithmic derivative sequence $E_0(x),\dots,E_d(x)$.

Next, we must extend the definition to generators of $G_{X}$ of the form $E^{(d)}(x)^a$.
If $a=n \in \mathbb{N}$, then we can define $\sigma_0(E^{(d)}(x)^n) := \sigma_0(E^{(d)}(x))^n$. 
If $a \not \in \mathbb{N}$, then we cannot define $\sigma_0(E^{(d)}(x)^a)$ to be $\sigma_0(E^{(d)}(x))^a$ because $\sigma_0(E^{(d)}(x))^a$ is not formally an element of $\boldsymbol{k}((H_{X}))$ (unless $d=0,1$).
To represent $\sigma_0\left(E^{(d)}(x)\right)^a$ as an element of $\boldsymbol{k}((H_{X}))$, we reason as follows:
\begin{align*}
    \sigma_0\left(E^{(d)}(x)\right)^a 
        &= \big(E_0(x)E_1(x)^d + \cdots + E_0(x)\cdots E_d(x)\big)^a 
        \\
        &= E_0(x)^aE_1(x)^{da}\left(1 + \frac{E_1(x)^{d-1}E_2(x) + \cdots + E_1(x)\cdots E_d(x)}{E_1(x)^d}\right)^a 
        \\
        &= E_0(x)^a E_1(x)^{da} \sum_{k=0}^{\infty} \binom{a}{k} \left(\frac{E_1(x)^{d-1}E_2(x) + \cdots + E_1(x)\cdots E_d(x)}{E_1(x)^d}\right)^k
\end{align*}
where the infinite sum in the final line comes from the Taylor series for $(1 + (\cdot))^a$ with $\binom{a}{k} = \frac{a(a-1)\cdots (a-k+1)}{k!}$.
The sum in the final line expands to a valid element of $\boldsymbol{k}((H_{X}))$, so we take this to be the definition of $\sigma_0\left(E^{(d)}(x)^a\right)$.
Note that this definition extends the definition of $\sigma_0\left(E^{(d)}(x)\right)^n$ for $n \in \mathbb{N}$ above.

To summarize, we define $\sigma_0 : \boldsymbol{k}[G_{X}] \to \boldsymbol{k}((H_{X}))$ by
\begin{enumerate}
    \item $\sigma_0(E(x)^a) = E_0(x)^a$ 
    
    \item $\sigma_0(E'(x)^a) = E_0(x)^aE_1(x)^a$
    
    \item $\sigma_0(E^{(d)}(x)^a) = E_0(x)^a E_1(x)^{da} \sum_{k=0}^{\infty} \binom{a}{k} \left(\frac{E_1(x)^{d-1}E_2(x) + \cdots + E_1(x)\cdots E_d(x)}{E_1(x)^d}\right)^k$ for $d > 1$
    
    \item Extend $\sigma_0$ to products and sums so that it is a $\boldsymbol{k}$-algebra homomorphism, i.e., for $g_1,\dots,g_n \in G_{X}$ and $c_1,\dots,c_n \in \boldsymbol{k}$, define
    \begin{enumerate}
        \item $\sigma_0(g_1 \cdots g_n) = \sigma_0(g_1) \cdots \sigma_0(g_n)$
        
        \item $\sigma_0(c_1g_1 + \cdots + c_ng_n) = c_1\sigma_0(g_1) + \cdots c_n\sigma_0(g_n)$
    \end{enumerate}
\end{enumerate}

We must check that $\sigma_0$ is well defined, i.e., that 
\begin{enumerate}
    \item for each generator $g$ of $G_{X}$, $\sigma_0(g)$ and $g$ satisfy the same relations, and
    
    \item for each $s \in \boldsymbol{k}[G_{X}]$, $\sigma_0(s)$ is indeed a sum in the Hahn series field $\boldsymbol{k}((H_{X}))$.
\end{enumerate}
First, the only relations among the generators of $G_{X}$ are
\begin{align*}
    E^{(d)}(x)^{a+b} &= E^{(d)}(x)^aE^{(d)}(x)^b.
\end{align*}
For $d > 1$, we can compute that
\begin{align*}
    \sigma_0\left(E^{(d)}(x)^{a+b}\right) 
        &= E_0(x)^{a+b} E_1(x)^{d(a+b)} \sum_{k=0}^{\infty} \binom{a+b}{k} \left(\frac{\sigma_0(E^{(d)}(x))}{E_0(x) E_1(x)^{d}}-1\right)^k \\
        &= E_0(x)^a E_1(x)^{da} \sum_{k=0}^{\infty} \binom{a}{k} \left(\frac{\sigma_0(E^{(d)}(x))}{E_0(x) E_1(x)^{d}}-1\right)^k \\
        & \cdot E_0(x)^b E_1(x)^{db} \sum_{k=0}^{\infty} \binom{b}{k} \left(\frac{\sigma_0(E^{(d)}(x))}{E_0(x) E_1(x)^{d}}-1\right)^k  \\
        &= \sigma_0\left(E^{(d)}(x)^a\right)\sigma_0\left(E^{(d)}(x)^b\right).
\end{align*}
Also, we immediately have
\begin{align*}
    \sigma_0\left(E(x)^{a+b}\right) 
        &= \sigma_0\left(E(x)^a\right)\sigma_0\left(E(x)^b\right)
        \\
    \sigma_0\left(E'(x)^{a+b}\right) 
        &= \sigma_0\left(E'(x)^a\right)\sigma_0\left(E'(x)^b\right).
\end{align*}

Second, $\sigma_0(g)$ is a valid sum in $\boldsymbol{k}((H_{X}))$ for each generator $g$ of $G_{X}$, and only finitely many generators of $G_{X}$ appear in any $s \in \boldsymbol{k}[G_{X}]$.
So $\sigma_0(s)$ is a valid sum in $\boldsymbol{k}((H_{X}))$.

\begin{rmk}\label{image of sigma_0}
    Let $s \in \boldsymbol{k}[G_{X}]$.
    Because $\sigma_0$ is defined using the Taylor expansion of $(1 + (\cdot))^a$, any monomial in $\mathrm{Supp} (\sigma_0(s))$ must be of the form
        $$\prod_{j=1}^p E_0(x_j)^{\alpha_{j}} E_1(x_j)^{\beta_{j}} E_2(x_j)^{n_{j,2}} \cdots E_d(x_j)^{n_{j,d}}$$
    with $\alpha_{j}, \beta_j \in \boldsymbol{k}$ and $n_{j,2},\dots,n_{j,d} \in \mathbb{N}$ for all $j=1,\dots,p$, for some $d,k \in \mathbb{N}$.
    The key observation is that the exponents of $E_2,\dots,E_d$ generators are not just any elements of $\boldsymbol{k}$---they must be \textit{natural numbers}.
\end{rmk}

\begin{rmk}[Ordering and separation assumptions]\label{order props}
We will now give assumptions from which we can induce a total ordering of $\boldsymbol{k}[G_{X}]$ using the order on $\boldsymbol{k}((H_{X}))$:
\begin{enumerate}
    \item $\boldsymbol{k}$ is an ordered field.
    
    \item $X$ is a subset of an ordered field $L$.
    
    \item We have the following partial order:
    For all $m \in \mathbb{N}$, all $x, y \in X$ with $x > y$, and all $a \in \boldsymbol{k}$
    \begin{enumerate}
        \item $E(x-m) > \boldsymbol{k}$ 
        \item $E(x-m) > E(y-m)^a$.
    \end{enumerate}
    
    \item There is a map $r : X \times X \to \mathbb{Q} \cap (0,1)$ such that for all $x,y \in X$ with $x > y$, we have $x-y < r(x,y)$.
\end{enumerate}
\end{rmk}

Notice that if we identify
    $$E_0(x-m) = E(x-m)$$
then the ordering already defined on each $\boldsymbol{k}((H_{X-m}))$ extends the partial order generated by these assumptions.
Also, if $X$ is a subset of an ordered field, then we can allow the subtraction symbol in the \textit{purely formal} notation ``$x-m$" to have algebraic meaning.

\begin{lemma}\label{series comp}
    Suppose $X$ and $\boldsymbol{k}$ satisfy the ordering and separation assumptions in Remark \ref{order props}.
    Then $\sigma_0$ is injective.
\end{lemma}

\begin{proof}
Let $0 \ne s = \sum_{i=1}^n c_ig_i \in \boldsymbol{k}[G_{X}]$ with all $c_i \ne 0$.
Then $\sigma_0(s)$ is a possibly infinite sum in $\boldsymbol{k}((H_{X}))$.
To show that $\sigma_0(s) \ne 0$, we will find a monomial with nonzero coefficient in $\sigma_0(s)$.

Enumerate all the variables that appear in $g_1,\dots,g_n$ as $x_1 > \cdots > x_p$.
Split each $g_i = \prod_{j=1}^p g_i(x_j)$ into blocks for each of the $x_j$'s.
Let $d$ be largest such that for some $j$, $E^{(d)}(x_j)$ appears in some $g_i$.
By definition of $\sigma_0$, we can express
\begin{align*}
    \sigma_0(g_i(x_j)) &= \sigma_0(E(x_j)^{a_{i,j,0}})\sigma_0( E'(x_j)^{a_{i,j,1}}) \sigma_0( E''(x_j)^{a_{i,j,2}}) \cdots \sigma_0(E^{(d)}(x_j)^{a_{i,j,d}}) \\
        &= E_0(x_j)^{a_{i,j,0} + \cdots + a_{i,j,d}} E_1(x_j)^{a_{i,j,1} + 2a_{i,j,2} + \cdots + da_{i,j,d}} \\
        &\left(\sum_{k_2=0}^{\infty}\binom{a_{i,j,2}}{k_2}\left(\frac{E_2(x_j)}{E_1(x_j)}\right)^{k_2}\right) \cdots \left(\sum_{k_{d}=0}^{\infty}\binom{a_{i,j,d}}{k_d}\left(\cdots + \frac{E_2(x_j) \cdots E_{d}(x_j)}{E_1(x_j)^{d-1}}\right)^{k_{d}}\right)
\end{align*}
with $a_{i,j,0},\dots,a_{i,j,d} \in \boldsymbol{k}$.
When we fully expand out this product of sums, we have that for each $k_{d} \in \mathbb{N}$, the coefficient of the term
    $$E_0(x_j)^{a_{i,j,0} + \cdots + a_{i,j,d}} E_1(x_j)^{a_{i,j,1} + 2a_{i,j,2} + \cdots + da_{i,j,d}}\left(\frac{E_2(x_j) \cdots E_{d}(x_j)}{E_1(x_j)^{d-1}}\right)^{k_{d}}$$ 
is $\binom{a_{i,j,d}}{k_d}$.
This is because the only times $E_{d}(x_j)$ appears anywhere in $\sigma_0(g_i(x_j))$ come from instances of the monomial $\frac{E_2(x_j) \cdots E_{d}(x_j)}{E_1(x_j)^{d-1}}$.
So the coefficient of the term
    $$\left(\frac{E_2(x_j) \cdots E_{d}(x_j)}{E_1(x_j)^{d-1}}\right)^{k_{d}}\prod_{l=1}^p E_0(x_l)^{a_{i,l,0} + \cdots + a_{i,l,d}} E_1(x_l)^{a_{i,l,1} + 2a_{i,l,2} + \cdots + da_{i,l,d}}$$ 
in $\sigma_0(g_i) = \sigma_0(g_i(x_1))\cdots \sigma_0(g_i(x_p))$ is still $\binom{a_{i,j,d}}{k_d}$, since $E_{d}(x_j)$ does not appear in  $\sigma_0(g_i(x_l))$ for $l \ne j$.

Let $I \subset \{1,\dots,n\}$ be maximal such that for all $j = 1,\dots,p$ and all $i_1,i_2 \in I$, the exponents of $E_0(x_j)$ and $E_1(x_j)$ in $g_{i_1}$ are the same as in $g_{i_2}$, i.e., we have 
\begin{align*}
    a_{i_1,j,0} + \cdots + a_{i_1,j,d} &= a_{i_2,j,0} + \cdots + a_{i_2,j,d} \\
    a_{i_1,j,1} + 2a_{i_1,j,2} + \cdots + da_{i_1,j,d} &= a_{i_2,j,1} + 2a_{i_2,j,2} + \cdots + da_{i_2,j,d}.
\end{align*}
For the remainder of the argument, fix some $j \in \{1,\dots,p\}$.
Then the coefficient of 
    $$\left(\frac{E_2(x_j) \cdots E_{d}(x_j)}{E_1(x_j)^{d-1}}\right)^{k_{d}}\prod_{l=1}^p E_0(x_l)^{a_{i,l,0} + \cdots + a_{i,l,d}} E_1(x_l)^{a_{i,l,1} + \cdots + da_{i,l,d}}$$ 
in $\sigma_0(\sum_{i=1}^nc_ig_i)$ is 
    $$\sum_{i \in I}c_i\binom{a_{i,j,d}}{k_d}$$
Suppose the $a_{i,j,d}$ for $i \in I$ are distinct. 
If $\sum_{i \in I}c_i\binom{a_{i,j,d}}{k_d} = 0$ for $k_{d}=0,\dots,|I|-1$, then we would have $c_{i_1} = \dots = c_{i_{|I|}} = 0$, a contradiction. 
If $\sum_{i \in I}c_i\binom{a_{i,j,d}}{k_d} \ne 0$ for some $k_{d} < |I|$, then we are done because this shows $\sigma_0(s) \ne 0$.

So assume the $a_{i,j,d}$ for $i \in I$ are not all distinct.
Let $1 < q < |I|$ be the number of distinct values among $a_{i,j,d}$ for $i \in I$.
Let $(P_1,\dots,P_q)$ partition $I$ so that for any $i_1,i_2 \in P_l$, $\alpha_{l,j} := a_{i_1,j,d} = a_{i_2,j,d}$,
and $\alpha_{l_1,j} \ne \alpha_{l_2,j}$ for $l_1 \ne l_2$.
Consider the sums
    $$\sum_{i \in P_l}c_i\frac{g_i}{E^{(d)}(x_j)^{\alpha_{l,j}}}$$
for $l = 1,\dots,q$.
Each sum has strictly fewer terms than $s$, the largest number $r$ such that $E^{(r)}$ appears is less than $d$, and if $\alpha_{l,j} \ne 0$ then the sum's monomials have one less multiplicand than the corresponding monomials of $s$.

We will proceed by induction with the hypothesis that for each $l = 1,\dots, q$, we can find a leading term $b_lt_l$ of $\sigma_0\left(\sum_{i \in P_l}c_i\frac{g_i}{E^{(d)}(x_j)^{\alpha_{l,j}}}\right)$, where $0 \ne b_l \in \boldsymbol{k}$ and $t_l \in G_{X}$. 
Since $E^{(d)}(x_j)$ does not appear in $\sum_{i \in P_l}c_i\frac{g_i}{E^{(d)}(x_j)^{\alpha_{l,j}}}$, $E_d(x_j)$ does not appear in $t_l$.
Now let $I_0 \subset \{1,\dots,q\}$ be the set of indices $l$ at which
    $$E_0(x_j)^{\alpha_{l,j}} E_1(x_j)^{d\alpha_{l,j}} t_l$$
is maximized in $H_{X}$.

Consider the sum
    $$\sum_{l \in I_0} \sigma_0(E^{(d)}(x_j)^{\alpha_{l,j}})b_lt_l = \sum_{l \in I_0} E_0(x_j)^{\alpha_{l,j}} E_1(x_j)^{d\alpha_{l,j}} b_lt_l \sum_{k_d=0}^{\infty}\binom{\alpha_{l,j}}{k_d}\left(\dots + \frac{E_2(x_j) \cdots E_{d}(x_j)}{E_1(x_j)^{d-1}}\right)^{k_{d}}.$$
For each $k_d$, the following is the largest monomial of the sum with $k_d$ as the exponent of $E_{d}(x_j)$:
    $$E_0(x_j)^{\alpha_{l,j}} E_1(x_j)^{d\alpha_{l,j}} t_l\left(\frac{E_2(x_j) \cdots E_{d}(x_j)}{E_1(x_j)^{d-1}}\right)^{k_{d}}.$$
Its coefficient is
    $$\sum_{l \in I_0} \binom{\alpha_{l,j}}{k_d}b_l$$
since $E_d(x_j)$ only appears in instances of $\frac{E_2(x_j) \cdots E_{d}(x_j)}{E_1(x_j)^{d-1}}$.
If the coefficient is 0 for all $k_d=0,\dots,|I_0|-1$, then that forces $b_1 = \cdots = b_q = 0$, since the $\alpha_{l,j}$'s are all distinct.
But this contradicts the induction hypothesis.
So some $\sum_{l=1}^q\binom{\alpha_{l,j}}{k_d}b_l$ must be nonzero for $k_d < |I_0|$, which means we have found a monomial of $\sigma_0(s)$ with nonzero coefficient.
\end{proof}

Since $\sigma_0$ is injective, the ordering of $\boldsymbol{k}((H_{X}))$ induces an order on $\boldsymbol{k}[G_{X}]$ by defining $s > 0$ if and only if $\sigma_0(s) > 0$.

\subsection{Defining the order preserving embedding \texorpdfstring{$\nu_0$}{v0}}

In Lemma \ref{log derivative shift}, we showed that for the functions $E_d$, $E_{d+1}$, and $\epsilon_d$ with $d \ge 2$, we have 
    $$E_d(x) = E_{d-1}(x-1)(1 + \epsilon_d(x)).$$
The proof of Lemma \ref{log derivative shift} does not actually use any properties of $E_d$, $E_{d+1}$, and $\epsilon_d$ as functions.
Only the identities $E_0 = E$, $E_{d+1} = \frac{E_d'}{E_d}$ for $d \in \mathbb{N}$, and the difference-differential identities for the derivatives of $E$ are needed.
Thus, the conclusion of Lemma \ref{log derivative shift} still holds in any purely formal context in which the necessary identities have been imposed.

Our goal is to formally build an ordered field of series that embeds the germs at $+\infty$ of $\mathcal{L}_{\mathrm{transexp}}$-terms as an ordered differential field.
Thus, we must impose the conclusion of Lemma \ref{log derivative shift} if we hope to build something consistent with the identities used to prove this Lemma. 

In this subsection, we will first show how to represent $\epsilon_d(x)$ as an element of $\boldsymbol{k}((H_{X-1}))$ for each $d \in \mathbb{N}$.
We will also prove a lemma about the form this representation takes, which will be necessary in the next subsection.
Then we will define an order preserving embedding
    $$\nu_0 : \sigma_0(\boldsymbol{k}[G_{X}]) \to \boldsymbol{k}((H_{X-1})).$$
It is important that we restrict the domain of $\nu_0$ to $\sigma_0(\boldsymbol{k}[G_{X}]) \subsetneq \boldsymbol{k}((H_{X}))$.
Sums in the image of $\sigma_0$ can be infinite, but they are limited in two key ways:
\begin{enumerate}
    \item Only finitely many elements of $X$ appear in any sum.
    
    \item For each sum, there is some $d \in \mathbb{N}$ such that only $E_0,\dots,E_d$ may appear and $E_{d+1}, E_{d+2},\dots$ do not.
\end{enumerate}
Arbitrary elements of $\boldsymbol{k}((H_{X}))$ do not have these two properties, and this causes issues when trying to extend the definition of $\nu_0$ on $\sigma_0(\kk[G_{X}])$ to all of $\boldsymbol{k}((H_{X}))$ in the natural way.
Fortunately, we will only ever need to embed elements of $\sigma_0(\boldsymbol{k}[G_{X}])$ into $\boldsymbol{k}((H_{X-1}))$, so this is not a problem for us.

Recall from Definition \ref{Boshernitzan's sequence with integer powers} that $H_{X-1}' \subset H_{X-1}$ is the subgroup of monomials with integer exponents.
Fix $x \in X$.
We will represent each $\epsilon_d(x)$ by an element of $\boldsymbol{k}((H_{X-1}')) \subset \boldsymbol{k}((H_{X-1}))$, by induction.
In the base case, we can immediately represent $\epsilon_2(x)$ by $\frac{E_2(x-1)}{E_1(x-1)} \in \boldsymbol{k}((H_{X-1}))$.
We can also immediately identify the $n$th derivative of $\epsilon_2$ with $(\partial_1)^{\circ n}\left(\frac{E_2(x-1)}{E_1(x-1)}\right) \in \boldsymbol{k}((H_{X-1}'))$ for all $n \in \mathbb{N}$.

Now suppose $\epsilon_d(x)$ is represented by $s_d \in \boldsymbol{k}((H_{X-1}'))$.
We will use the formula
    $$\epsilon_{d+1}(x) = \frac{\epsilon_d'(x)}{E_d(x-1)} \cdot \frac{1}{1 + \epsilon_d(x)}$$
from Lemma \ref{log derivative shift}, replacing $\epsilon_d(x)$ by $s_d$, $\epsilon_d'(x)$ by $\partial_1(s_d)$, and $\frac{1}{1 + \epsilon_d(x)}$ by an infinite sum.
Following the formula, we represent $\epsilon_{d+1}(x)$ by
\begin{align*}
    \frac{\partial_1(s_d)}{E_d(x-1)}\sum_{j=0}^{\infty} \left(-s_d \right)^j \in \boldsymbol{k}((H_{X-1}')).
\end{align*}
We can then represent the $n$th derivative of $\epsilon_{d+1}$ by
\begin{align*}
    (\partial_1)^{\circ n}\left(\frac{\partial_1(s_d)}{E_d(x-1)}\sum_{j=0}^{\infty} \left(-s_d \right)^j\right) \in \boldsymbol{k}((H_{X-1}'))
\end{align*}
for all $n \in \mathbb{N}$.
In what follows, we will write $\epsilon_d(x)$ instead of writing out the sums as above in order to highlight the connection with the computations of Subsection \ref{Boshernitzan's sequence}.

\begin{lemma}\label{exponent of E_d in epsilon_d}
    Let $x \in X$ and $d \in \mathbb{N}$. Then we have the following:
    \begin{enumerate}
        \item The exponent of $E_d(x-1)$ in any monomial of $\epsilon_d(x)$ is either 0 or 1.
        
        \item $E_{d+k}(x-1)$ does not appear in $\epsilon_d(x)$ for any $k \ge 1$.
        
        \item The sum of the exponents of generators in any monomial of $\epsilon_d(x)$ is 0.
    \end{enumerate}
\end{lemma}

\begin{proof}
    All three claims are immediate for $d=2$.
    Suppose the results hold for $\epsilon_d(x)$.
    We will show they hold for $\epsilon_{d+1}(x)$.
    
    Since $\epsilon_{d+1}(x) = \frac{\epsilon_d'(x)}{E_d(x-1)} \sum_{j=1}^{\infty} (-\epsilon_d(x))^j$, every monomial of $\epsilon_{d+1}(x)$ is of the form
        $$\frac{M_1}{E_d(x-1)}M_2$$
    where $M_1$ is a monomial of $\epsilon_d'(x)$ and $M_2$ is a monomial of $\sum_{j=1}^{\infty} (-\epsilon_d(x))^j$.
    
    Since $E_{d+k}(x-1)$ does not appear in $\epsilon_d(x)$, it does not appear in $\sum_{j=1}^{\infty} (-\epsilon_d(x))^j$ either.
    Also, since the sum of exponents in every monomial of $\epsilon_d(x)$ is 0, the same is true for $\sum_{j=1}^{\infty} (-\epsilon_d(x))^j$.
    So we can ignore $M_2$ when proving (1), (2), and (3).
    
    Now we will consider the possible forms $M_1$ can take.
    Any monomial of $\epsilon_d'(x)$ arises as a monomial of $\partial_1(M)$ for some monomial $M$ of $\epsilon_d(x)$.
    By induction hypothesis, any monomial of $\epsilon_d(x)$ is of the form $E_1(x-1)^{n_1} \cdots E_{d}(x-1)^{n_{d}}$ with $n_1,\dots,n_{d-1} \in \mathbb{Z}$ and $n_d \in \{0,1\}$.
    Then
        $$\partial_1\big(E_1(x-1)^{n_1} \cdots E_{d}(x-1)^{n_{d}}\big) = \sum_{l=1}^{d} n_l \big(E_1(x-1)^{n_1} \cdots E_{d}(x-1)^{n_{d}}\big) \cdot E_{l+1}(x-1).$$
    In this sum, $E_{d+1}(x-1)$ can only appear in the final summand with exponent 0 or 1, which finishes (1).
    No $E_{d+k}(x-1)$ with $k \ge 2$ appears, which finishes (2).
    For (3), observe that each monomial in the sum has one new generator $E_{l+1}(x-1)$ with exponent 1.
    So the sum of exponents in any $M_1$ must be 1, which means the total sum of exponents in $\frac{M_1}{E_d(x-1)}M_2$ is 0.
\end{proof}

We now define an order-preserving field embedding $\nu_0 : \sigma_0(\boldsymbol{k}[G_{X}]) \to \boldsymbol{k}((H_{X-1}))$ with the intuition of rewriting elements of $\sigma_0(\boldsymbol{k}[G_{X}]) \subset \boldsymbol{k}((H_{X}))$ using the difference equations
    $$E_d(x) = E_{d-1}(x-1)(1 + \epsilon_d(x))$$
from Lemma \ref{log derivative shift}.
Define $\nu_0 : \sigma_0(\boldsymbol{k}[G_{X}]) \to \boldsymbol{k}((H_{X-1}))$ as follows for all $a \in \boldsymbol{k}$ and $n \in \mathbb{N}$:
\begin{enumerate}
    \item $\nu_0(E_0(x)^a) = E_0(x)^a$ 
    
    \item $\nu_0(E_1(x)^a) = E_{0}(x-1)^a E_{1}(x-1)^a$ 
    
    \item $\nu_0(E_d(x)^n) = E_{d-1}(x-1)^n(1 + \epsilon_d(x))^n$  for $d \ge 2$. 
    
    \item Extend $\nu_0$ to products and sums so that it is a $\boldsymbol{k}$-algebra homomorphism. 
    For generators $g_1,\dots,g_l$ of $H_{X}$, let 
        $$\nu_0(g_1 \cdots g_l) = \nu_0(g_1) \cdots \nu_0(g_l).$$
    For $\sum_{g \in H_{X}} c_g g \in \sigma_0(\boldsymbol{k}[G_{X}])$, let
        $$\nu_0\left(\sum_{g \in H_{X}} c_g g\right) = \sum_{g \in H_{X}} c_g \nu_0(g).$$
\end{enumerate}

We can check with an easy computation that for each $g \in H_{X}$, $\nu_0(g)$ satisfies the same relations as $g$:
The only relations among the elements of $H_{X}$ are 
    $$E_d(x)^aE_d(x)^b = E_d(x)^{a+b}$$
for all $d \in \mathbb{N}$, $x \in X$, and $a,b \in \boldsymbol{k}$.
We must check that $\nu_0(E_d(x)^a)\nu_0(E_d(x)^b) = \nu_0(E_d(x)^{a+b})$.
This is clear for $d = 0,1$, so suppose $d > 1$.
Then
\begin{align*}
    \nu_0(E_d(x)^n)\nu_0(E_d(x)^m) &= 
        E_{d-1}(x-1)^n \sum_{j=0}^{n} \binom{n}{j} \epsilon_d(x)^j 
        \cdot E_{d-1}(x-1)^m \sum_{j=0}^{m} \binom{m}{j} \epsilon_d(x)^j  \\
        &= E_{d-1}(x-1)^{n+m} \sum_{j=0}^{n+m} \binom{n+m}{j} \epsilon_d(x)^j.
\end{align*}
So $\nu_0(E_d(x)^n)$ satisfies the same relations as $E_d(x)^n$.

We now check that $\nu_0$ does indeed map to $\boldsymbol{k}((H_{X-1}))$.

\begin{lemma} \label{nu well defined}
    For every $s \in \boldsymbol{k}[G_{X}]$, $\nu_0(\sigma_0(s))$ is an element of the Hahn series field $\boldsymbol{k}((H_{X-1}))$.
\end{lemma}

\begin{proof}
Let $\sigma_0(s) = \sum_{g \in H_{X}} c_g g$.
We will show that $\nu_0\left(\sum_{g \in H_{X}} c_g g\right)$ is a valid sum in $\boldsymbol{k}((H_{X-1}))$.
Since each $g \in H_{X}$ is a finite product of generators, $\nu_0$ is well defined on monomials.
So it suffices to check that $\big(c_g \nu_0(g) : g \in H_{X}\big)$ is summable, i.e., that
\begin{enumerate}
    \item For each $h \in H_{X-1}$ there are only finitely many $g \in H_{X}$ such that $c_g \ne 0$ and $h \in \mathrm{Supp} (\nu_0(g))$.
    
    \item $\displaystyle \bigcup_{g \in \mathrm{Supp} (\sigma_0(s))} \mathrm{Supp} (\nu_0(g))$ is reverse well-ordered in $H_{X-1}$.
\end{enumerate}

First we introduce some notation and setup.
Let $x_1 > x_2 > \dots > x_p$ be the elements of $X$ that appear in $s$, and let $d$ be order of the largest derivative that appears in $s$.
Enumerate $\mathrm{Supp}(\sigma_0(s))$ as $(g_{i} : i < \delta)$, a reverse well-ordered sequence in $H_X$.
By Remark \ref{image of sigma_0}, we can write
\begin{align*}
    g_{i} = \prod_{j=1}^p 
    E_0(x_j)^{\alpha_{i,j}} \cdot E_1(x_j)^{\beta_{i,j}} E_2(x_j)^{n_{i,j,2}} \cdots E_d(x_j)^{n_{i,j,d}}
\end{align*}
with $\alpha_{i,j},\beta_{i,j} \in \boldsymbol{k}$ and $n_{i,j,2},\dots,n_{i,j,d} \in \mathbb{N}$ for $j=1,\dots,p$.
Then 
\begin{align*}
    \nu_0(g_i) =\prod_{j=1}^p \Big( 
        &E_0(x_j)^{\alpha_{i,j}} \cdot E_0(x_j-1)^{\beta_{i,j}}E_1(x_j-1)^{\beta_{i,j}} 
        \\
        &\cdot E_1(x_j-1)^{n_{i,j,2}}(1 + \epsilon_2(x_j))^{n_{i,j,2}} \cdots E_{d-1}(x_j-1)^{n_{i,j,d}}(1 + \epsilon_d(x_j))^{n_{i,j,d}}\Big) 
        \\
        = \prod_{j=1}^p \Big( 
        &E_0(x_j)^{\alpha_{i,j}} \cdot E_0(x_j-1)^{\beta_{i,j}} E_1(x_j-1)^{\beta_{i,j} + n_{i,j,2}} E_2(x_j-1)^{n_{i,j,3}} \cdots E_{d-1}(x_j-1)^{n_{i,j,d}}
        \\
        &\cdot (1 + \epsilon_2(x_j))^{n_{i,j,2}} \cdots (1 + \epsilon_{d-1}(x_j))^{n_{i,j,d}} \Big).
\end{align*}
We will use the following observations about the final line $(1 + \epsilon_2(x_j))^{n_{i,j,2}} \cdots (1 + \epsilon_{d-1}(x_j))^{n_{i,j,d}}$ above, which follow from the fact that $\epsilon_2,\dots,\epsilon_d \in \boldsymbol{k}((H_{X-1}'))$:
\begin{enumerate}
    \item [(A)] All exponents in all monomials are integers.
    
    \item [(B)] Only $E_1(x_j-1), \dots, E_{d}(x_j-1)$ may appear for $j=1,\dots,p$. No $E_0$ generator may appear.
    Thus, the exponents of the $E_0$ generators are fixed across all monomials of $\nu_0(g_i)$.
\end{enumerate}

To finish the setup, we introduce notation for tuples of exponents in $g_i$ for $i < \delta$.
Let
\begin{itemize}
        \item $\overline{\alpha_{i}} = (\alpha_{i,1}, \alpha_{i,2},\dots,\alpha_{i,p})$
        
        \item $\overline{\beta_i} = (\beta_{i,1}, \beta_{i,2},\dots,\beta_{i,p})$
        
        \item $\overline{n_{i,l}} = (n_{i,1,l}, n_{i,2,l},\dots,n_{i,p,l})$ for $l = 2,\dots,d$.
    \end{itemize}

Now we will prove (1).
It suffices to show that if $h \in \mathrm{Supp} (\nu_0(g_0))$, then $h \in \mathrm{Supp} (\nu_0(g_{i}))$ for only finitely many $i < \delta$.

So suppose $h \in \mathrm{Supp} (\nu_0(g_0))$.
Then by Observation (B), $h$ must be of the form
\begin{align*}
    h = \prod_{j=1}^p 
        &E_0(x_j)^{\alpha_{0,j}} E_0(x_j-1)^{\beta_{0,j}} E_1(x_j-1)^{\beta_{0,j} + N_{j,1}} E_2(x_j-1)^{N_{j,2}} \cdots E_d(x_j-1)^{N_{j,d}}
\end{align*}
with $(N_{j,1},\dots,N_{j,d}) \in \mathbb{Z}^d$ by Observation (A).

Similarly, if $h \in \mathrm{Supp} (\nu_0(g_{i}))$ for $i > 0$, then again by Observations (A) and (B), $h$ must be of the form
\begin{align*}
    h = \prod_{j=1}^p 
        &E_0(x_j)^{\alpha_{i,j}} E_0(x_j-1)^{\beta_{i,j}} E_1(x_j-1)^{\beta_{i,j} + M_{j,1}} E_2(x_j-1)^{M_{j,2}} \cdots E_d(x_j-1)^{M_{j,d}}
\end{align*}
with $(M_{j,1},\dots,M_{j,d}) \in \mathbb{Z}^d$.
So we must have $\alpha_{i,j} = \alpha_{0,j}$ and $\beta_{i,j} = \beta_{0,j}$ for $j=1,\dots,p$.

By the way the order on $H_{X-1}$ is defined, since $g_{i} < g_0$, we must have
    $$(\overline{\alpha_{i}},\overline{\beta_{i}},\overline{n_{i,2}},\overline{n_{i,3}},\dots,\overline{n_{i,d}}) < (\overline{\alpha_{0}},\overline{\beta_{0}},\overline{n_{0,2}},\overline{n_{0,3}},\dots,\overline{n_{0,d}})$$
in the lexicographic order on $\boldsymbol{k}^{2p} \times \mathbb{N}^{p(d-1)}$.
Since the first parts of each of these sequences are equal, we must have
    $$(\overline{n_{i,2}},\overline{n_{i,3}},\dots,\overline{n_{i,d}}) < (\overline{n_{0,2}},\overline{n_{0,3}},\dots,\overline{n_{0,d}})$$
in the lexicographic order on $\mathbb{N}^{p(d-1)}$.
Any reverse well-ordered sequence in $\mathbb{N}^{p(d-1)}$ with the lexicographic order must be finite, so there can be only finitely many $i < \delta$ with $(\overline{n_{i,2}},\overline{n_{i,3}},\dots,\overline{n_{i,d}}) < (\overline{n_{0,2}},\overline{n_{0,3}},\dots,\overline{n_{0,d}})$, i.e., only finitely many $i$ with $h \in \mathrm{Supp} (\nu_0(g_{i}))$.
This finishes the proof of (1).

For (2), let $\varnothing \ne B \subset \bigcup_{g \in \Supp(\sigma_0(s))} \mathrm{Supp} (\nu_0(g))$.
We will show $B$ has a greatest element by building an increasing sequence $b_0, b_1, b_2, \dots$ in $B$ and showing it must terminate with a largest element.
In the base case, let $B_0 = B \ne \varnothing$.
Given $\varnothing \ne B_i \subset B$, we define $\gamma_i \in H_X$, $b_{i} \in B$, and $B_{i+1} \subset B$ as follows:
Since $\mathrm{Supp} (\sigma_0(s))$ is reverse well-ordered and $B_i \ne \varnothing$, there must be a greatest element $g \in \mathrm{Supp} (\sigma_0(s))$ such that $B_i \cap \nu_0(g) \ne \varnothing$.
\begin{enumerate}
    \item Let $\gamma_{i} := g$.
    
    \item Let $b_i := \max \big(B_i \cap \mathrm{Supp} (\nu_0(\gamma_i))\big)$, which exists since $\mathrm{Supp}(\nu_0(\gamma_i))$ is reverse well-ordered.
    
    \item Let $B_{i+1} := \{b \in B : b > b_i\}$.
\end{enumerate}
If $B_{i+1} \ne \varnothing$, we continue.
Note that $b_i \in B_i \setminus B_{i+1}$, so $B_{i+1} \subsetneq B_i$.
So when we continue, we will have $\gamma_{i+1} < \gamma_i$ and $b_{i+1} > b_i$.
If $B_{i+1} = \varnothing$, then $b_i$ is the largest element of $B$ and the sequence terminates.

We must show that the sequence $b_0, b_1, \dots$ terminates.
Again by Remark \ref{image of sigma_0}, write 
\begin{align*}
    \gamma_{i} = \prod_{j=1}^{p} E_0(x_{j})^{\alpha_{i,j}} \cdot E_1(x_{j})^{\beta_{i,j}} \cdot E_2(x_{j})^{n_{i,j,1}} \cdots E_{d}(x_{j})^{n_{i,j,d}}
\end{align*}
with $\alpha_{i,j,k}, \dots, \alpha_{i,j,0}, \beta_{i,j} \in \boldsymbol{k}$ and $n_{i,j,2},\dots,n_{i,j,d} \in \mathbb{N}$ for $j=1,\dots,p$.
Then, again by Observations (A) and (B), we can write $b_i$ as
\begin{align*}
    b_i = \prod_{j=1}^p 
        &E_0(x_j)^{\alpha_{i,j}} E_0(x_j-1)^{\beta_{i,j}} 
        \cdot E_1(x_j-1)^{\beta_{i,j} + N_{i,j,1}} E_2(x_j-1)^{N_{i,j,2}} \cdots E_d(x_j-1)^{N_{i,j,d}}
\end{align*}
with $N_{i,j,1},\dots,N_{i,j,d} \in \mathbb{Z}$.

We will now show that the inequalities $\gamma_i < \gamma_0$ and $b_i > b_0$ for $i > 0$ imply that the sequence $b_0,b_1,\dots$ terminates.
First, since $\gamma_i < \gamma_0$, we must have that 
    $$(\overline{\alpha_{i}}, \overline{\beta_i}) \le (\overline{\alpha_{0}}, \overline{\beta_0})$$
in the lexicographic order on $\boldsymbol{k}^{2p}$.
Since $b_i > b_0$, we must have 
    $$(\overline{\alpha_{i}}, \overline{\beta_i}) \ge (\overline{\alpha_{0}}, \overline{\beta_0})$$
in the lexicographic order on $\boldsymbol{k}^{2p}$.
Thus $\alpha_{i,j} = \alpha_{0,j}$ and $\beta_{i,j} = \beta_{0,j}$ for all $j=1,\dots,p$.

So once again, $\gamma_i < \gamma_0$ implies
    $$(\overline{n_{i,2}},\overline{n_{i,3}},\dots,\overline{n_{i,d}}) < (\overline{n_{0,2}},\overline{n_{0,3}},\dots,\overline{n_{0,d}})$$
in the lexicographic order on $\mathbb{N}^{p(d-1)}$.
Since $\mathrm{Supp} (\sigma_0(s))$ is reverse well-ordered, and any reverse well-ordered sequence in $\mathbb{N}^{p(d-1)}$ with the lexicographic order must be finite, the sequence $\gamma_1,\gamma_2,\dots$ must terminate after finitely many steps.
Thus, the sequence $b_0, b_1,\dots$ terminates too, and its final element is the largest element of $B$.
This completes the proof of (2).
\end{proof}

Thus the map $\nu_0 : \sigma_0(\boldsymbol{k}[G_{X}]) \to \boldsymbol{k}((H_{X-1}))$ is well defined.

The next corollary and remark follow from computations in the proof of Lemma \ref{nu well defined}, but we state them separately for clarity and completeness.

\begin{cor}\label{g_1 > g_2}
    If $s \in \sigma_0(\boldsymbol{k}[G_{X}])$ and $g_1 > g_2$ are two monomials of $s$, then $\mathrm{Lm}(\nu_0(g_1)) > \mathrm{Lm}(\nu_0(g_2))$.
\end{cor}

\begin{proof} 
Let $x_1 > x_2 > \dots > x_p$ be the elements of $X$ that appear in $g_1$ or $g_2$, and let $d$ be largest such that some $E_d$ generator appears in $g_1$ or $g_2$.
Then we can write
    $$g_i = \prod_{j=1}^p E_0(x_j)^{\alpha_{i,j}} \cdot E_1(x_j)^{\beta_{i,j}}\cdot E_2(x_j)^{n_{i,j,2}} \cdots E_d(x_j)^{n_{i,j,d}}$$
with $\alpha_{i,j},\beta_{i,j} \in \boldsymbol{k}$ and $n_{i,j,2},\dots,n_{i,j,d} \in \mathbb{N}$ for $j=1,\dots,p$ and $i = 1,2$.
Then
\begin{align*}
    \mathrm{Lm}(\nu_0(g_i)) = \prod_{j=1}^p 
        E_0(x_j)^{\alpha_{i,j}} E_0(x_j-1)^{\beta_{i,j}}
        E_1(x_j-1)^{\beta_{i,j}+n_{i,j,2}} \cdots E_{d-1}(x_j-1)^{n_{i,j,d}}.
\end{align*}
Since $g_1 > g_2$, we must have
    $$(\overline{\alpha_{1}}, \overline{\beta_1}, \overline{n_{1,2}}, \overline{n_{1,3}}, \dots, \overline{n_{1,d}}) > (\overline{\alpha_{2}}, \overline{\beta_2}, \overline{n_{2,2}}, \overline{n_{2,3}}, \dots, \overline{n_{2,d}})$$
in the lexicographic order on $\boldsymbol{k}^{p(d+1)}$.
So we also have $\mathrm{Lm}(\nu_0(g_1)) > \mathrm{Lm}(\nu_0(g_2))$.
\end{proof}

\begin{rmk}\label{comparing E_0 exponents}
    Suppose $g_1$ and $g_2$ are as in Corollary \ref{g_1 > g_2}, and also that
        $$(\overline{\alpha_{1}}, \overline{\beta_1}) > (\overline{\alpha_{2}}, \overline{\beta_2})$$
    in the lexicographic order on $\boldsymbol{k}^{2p}$.
    By observation (B) from the proof of Lemma \ref{nu well defined}, the exponents of the $E_0$ generators are fixed across all monomials of $\nu_0(g_i)$ for $i=1,2$.
    Thus, every monomial of $\nu_0(g_1)$ is greater than every monomial of $\nu_0(g_2)$, i.e.,
        $$\mathrm{Supp} (\nu_0(g_1)) > \mathrm{Supp} (\nu_0(g_2)).$$
\end{rmk}

\begin{cor}
    $\nu_0$ is order-preserving and thus injective.
\end{cor}

\begin{proof}
    Let $s \in \sigma_0(\boldsymbol{k}[G_{X}])$ and let $g = \mathrm{Lm}(s)$.
    If $g_1 \ne g$ is any other monomial in $\Supp \sigma_0(\kk[G_{X}])$, then by Lemma \ref{g_1 > g_2}, we know $\mathrm{Lm}(\nu_0(g_1)) > \mathrm{Lm}(\nu_0(g_2))$.
    Thus we must have 
        $$\mathrm{Lm}(\nu_0(s)) = \mathrm{Lm}(\nu_0(g_1))$$
    and
        $$\mathrm{Lt}(\nu_0(s)) = \mathrm{Lc}(s) \cdot \mathrm{Lt}(\nu_0(g_1)).$$
    So $\nu_0(s) > 0$ if and only if $\mathrm{Lc}(s) > 0$ if and only if $s > 0$.
\end{proof}

\subsection{Showing \texorpdfstring{$\boldsymbol{k}[G_{X}\Lambda_{X}]$}{k[G0L0]} is ordered} \label{order log sums}

To define an order on $\boldsymbol{k}[G_{X}]$, we embedded it into $\boldsymbol{k}((H_{X}))$, which induced an order on $\boldsymbol{k}[G_{X}]$.
The embedding was defined based on intuition from ``rewriting" expressions involving the derivatives of $E$ as polynomials in $E_0,E_1,\dots$.
We would like to replicate this idea to define an order on $\boldsymbol{k}[G_{X}\Lambda_{X}]$, but the generators $\log E'(x)^a$ of $\Lambda_{X}$ are an obstacle.
Recall that the difference-differential equation for $E'$ is $E'(x) = E(x) E'(x-1)$.
Taking log of both sides, we get
    $$\log E'(x) = E(x-1) + \log E'(x-1).$$
There is not a clear (to us) way of ``rewriting" $\log E'(x)$ in terms of $E_0(x-1),E_1(x-1),\dots$.
However, we can ``approximate" $\log E'(x)$ by
    $$\log E'(x) \approx E_0(x-1) + E_1(x-1).$$
We get this approximation by identifying $E(x-1) = E_0(x-1)$, and since $E'(x-1) < E(x-1)^2$ implies
    $$\log E'(x-1) < 2E(x-2) < E'(x-2) = E_1(x-1).$$

In this subsection, we will define an embedding
    $$\rho_0 : \boldsymbol{k}[G_{X}\Lambda_{X}] \to \boldsymbol{k}((H_{X-1}))$$
so that 
\begin{enumerate}
    \item $\rho_0$ is ``true" on elements of $\boldsymbol{k}[G_{X}]$, meaning that if $s \in \boldsymbol{k}[G_{X}]$, then $\rho_0(s) = \nu_0(\sigma_0(s))$ 
    
    \item $\rho_0(\log E'(x)) = E_0(x-1) + E_1(x-1).$
\end{enumerate}
We can extend $\rho_0$ to generators of the form $\log E'(x)^b$ of $\Lambda_{X}$ using the Taylor series for $(1 + (\cdot))^b$.
We will call $\rho_0$ the ``approximation map" because its effect will be to approximate every $t \in \boldsymbol{k}[G_{X}\Lambda_{X}]$ well enough by elements of $\boldsymbol{k}((H_{X-1}))$ that the sign of $t$ is determined by the ordering and separation assumptions on $\boldsymbol{k}$ and $X$ and the order on $\boldsymbol{k}((H_{X-1}))$.

It is reasonable to approximate $\log E'(x)$ by $E_0(x-1) + E_1(x-1)$ because the approximation does not affect the leading monomial $E_0(x-1)$.
We could work to find an element of $\boldsymbol{k}((H_{X-1}))$ that better approximates $\log E'(x)$, but we will show that our approximation is close enough to determine how $\boldsymbol{k}[G_{X}\Lambda_{X}]$ should be ordered.

To summarize, we define the ``approximation map" $\rho_0 : \boldsymbol{k}[G_{X}\Lambda_{X}] \to \boldsymbol{k}((H_{X-1}))$ by 
\begin{enumerate}
    \item $\displaystyle \rho_0(s) = \nu_0(\sigma_0(s))$ for $s \in \boldsymbol{k}[G_{X}]$
    
    \item $\displaystyle \rho_0(\log E'(x)^b) = E_0(x-1)^{b}\sum_{k=0}^{\infty}\binom{b}{k} \left(\frac{E_1(x-1)}{E_0(x-1)}\right)^k$ 
    
    \item Extend $\rho_0$ to products and sums so that it is a $\boldsymbol{k}$-algebra homomorphism, i.e., for $g_1,\dots,g_n \in G_{X}\Lambda_{X}$ and $c_1,\dots,c_n \in \boldsymbol{k}$, define
    \begin{enumerate}
        \item $\rho_0(g_1 \cdots g_n) = \rho_0(g_1) \cdots \rho_0(g_n)$
        
        \item $\rho_0(c_1g_1 + \cdots + c_ng_n) = c_1\rho_0(g_1) + \cdots c_n\rho_0(g_n)$.
    \end{enumerate}
\end{enumerate}
Note that $\rho_0$ is well defined for the same reasons $\sigma_0$ is well defined: it respects the relations among generators of $G_{X}$ and $\Lambda_{X}$, and for each $s \in \boldsymbol{k}[G_{X}\Lambda_{X}]$, $\rho_0(s)$ is a valid sum in $\boldsymbol{k}((H_{X-1}))$.

\begin{rmk}\label{initial subsum of rho(s)}
    Given $s \in \boldsymbol{k}[G_{X} \Lambda_{X}]$, we will define a sum $\mathrm{Init}(s)$, which we intend to be an initial subsum of $s$. 
    In Lemma \ref{rho_0 is injective}, we will show that $\mathrm{Init}(s) \ne 0$, and therefore
        $$\mathrm{Lm}\big(\rho_0(s)\big) = \mathrm{Lm}\big(\mathrm{Init}(s)\big).$$
    From this, we will conclude that $\rho_0$ is injective.
    
    We can write any $s \in \boldsymbol{k}[G_{X}\Lambda_{X}] = \boldsymbol{k}[G_{X}][\Lambda_{X}]$ as
        $$s = s_1\ell_1 + \cdots + s_n\ell_n$$
    with $s_i \in \boldsymbol{k}[G_{X}]$ and $\ell_i \in \Lambda_{X}$ for $i = 1,\dots,n$ and $\ell_{i_1} \ne \ell_{i_2}$ for $i_1 \ne i_2$.
    So let $s = s_1\ell_1 + \cdots + s_n\ell_n \in \boldsymbol{k}[G_{X}\Lambda_{X}]$ with $s_i \in \boldsymbol{k}[G_{X}]$ and $\ell_i \in \Lambda_{X}$.
    Let $x_1 > x_2 > \dots > x_p$ list the elements of $X$ appearing in $s$.
    Let $d$ be the largest derivative appearing in any $s_i$.
    
    We will define $\mathrm{Init}(s)$ by looking at $\rho_0(s)$ and figuring out what form its largest monomials could take.
    First, we look at the images of the $s_i$'s under $\rho_0$.
    By definition, $\rho_0(s_i) = \nu_0(\sigma_0(s_i))$.
    By Remark \ref{image of sigma_0}, for each $i=1,\dots,n$ write
        $$\mathrm{Lm}(\sigma_0(s_i)) = \prod_{j=1}^p E_0(x_j)^{\alpha_{i,j}} \cdot E_1(x_j)^{\beta_{i,j}} \cdot E_2(x_j)^{n_{i,j,2}} \cdots E_d(x_j)^{n_{i,j,d}}$$
    with $\alpha_{i,j}, \beta_{i,j} \in \boldsymbol{k}$ and $n_{i,j,2}, \dots, n_{i,j,d} \in \mathbb{N}$ for each $j=1,\dots,p$. 
    Let $t_i$ be the initial subsum of $\sigma_0(s_i)$ with monomials of the form
        $$\prod_{j=1}^p E_0(x_j)^{\alpha_{i,j}} \cdot E_1(x_j)^{\beta_{i,j}} \cdot E_2(x_j)^{m_{i,j,2}} \cdots E_d(x_j)^{m_{i,j,d}}$$
    with $m_{i,j,2},\dots,m_{i,j,d} \in \mathbb{N}$ and all other exponents the same as in $\mathrm{Lm}(\sigma_0(s_i))$.
    Then $t_i$ is a finite sum because any reverse well-ordered sequence of tuples of exponents in $\mathbb{N}^{p(d-1)}$ must be finite.
    
    For any monomial $M$ of $t_i$, we have
    \begin{align*}
        \nu_0(M) = \prod_{j=1}^p
            &E_0(x_j)^{\alpha_{i,j}} E_0(x_j-1)^{\beta_{i,j}} E_1(x_j-1)^{\beta_{i,j} + m_{i,j,2}} 
            \\
            &E_2(x_j-1)^{m_{i,j,3}} \cdots E_{d-1}(x_j-1)^{m_{i,j,d}} (1 + \epsilon_2(x_j))^{m_{i,j,2}} \cdots (1 + \epsilon_{d-1}(x_j))^{m_{i,j,d}}.
    \end{align*}
    If $M_1 \in \mathrm{Supp} (\sigma_0(s_i)) \setminus \mathrm{Supp} (t_i)$, then by Remark \ref{comparing E_0 exponents}, $\mathrm{Supp} (\nu_0(M)) > \mathrm{Supp} (\nu_0(M_1))$, i.e., every monomial of $\nu_0(M)$ is greater than every monomial of $\nu_0(M_1)$.
    
    Now we look at the images of the $\ell_i$'s under $\rho_0$.
    Let $\ell_i = \prod_{j=1}^p\log E'(x_j)^{b_{i,j}}$ for each $i=1,\dots,n$.
    Then 
    \begin{align*}
        \rho_0(\ell_i) = \prod_{j=1}^p E_0(x_j-1)^{b_{i,j}}\sum_{k=0}^{\infty}\binom{b_{i,j}}{k} \left(\frac{E_1(x_j-1)}{E_0(x_j-1)}\right)^k.
    \end{align*}
    Observe that the leading monomial of $\rho_0(\ell_i)$ is $\prod_{j=1}^p E_0(x_j-1)^{b_{i,j}}$.
    If 
        $$M_2 \in \mathrm{Supp} (\rho_0(\ell_i)) \setminus \{\mathrm{Lm}(\rho_0(\ell_i))\}$$
    then the exponent of $E_0(x_j-1)$ in $M_2$ is strictly less than $b_{i,j}$, for some $j = 1,\dots,p$.
    So 
        $$\mathrm{Lm}(\rho_0(\ell_i)) > \mathrm{Supp} (M_2).$$
    
    Altogether, this shows that 
        $$\mathrm{Supp} \big(\nu_0(t_i)\cdot \mathrm{Lm}(\rho_0(\ell_i))\big) > \mathrm{Supp}\big( \rho_0(s_i\ell_i) - \nu_0(t_i)\cdot \mathrm{Lm}(\rho_0(\ell_i))\big)$$
    i.e., $\nu_0(t_i)\cdot \mathrm{Lm}(\rho_0(\ell_i))$ is an initial subsum of $\rho_0(s_i\ell_i)$.
    
    Now we select the indices $i$ for which the initial subsum $\nu_0(t_i)\cdot \mathrm{Lm}(\rho_0(\ell_i))$ of $\rho_0(s_i\ell_i)$ may possibly contribute to the leading monomial of $\rho_0(s)$.
    Define $I_0 \subset \{1,\dots,n\}$ to be the set of indices such that 
        $$(\overline{\alpha_{i}},\overline{\beta_{i}} + \overline{b_{i}})$$
    is maximal in the lexicographic order on $\boldsymbol{k}^{2p}$ (where the addition $\overline{\beta_{i}} + \overline{b_{i}}$ is coordinatewise).
    If $i_0 \in I_0$ and $i_1 \not \in I_0$, then every monomial of $\nu_0(t_{i_0}) \cdot \mathrm{Lm}(\rho_0(\ell_{i_0}))$ is greater than every monomial of $\nu_0(t_{i_1}) \cdot \mathrm{Lm}(\rho_0(\ell_{i_1}))$ by Remark \ref{comparing E_0 exponents}.
    
    Finally, we define 
    \begin{align*}
        \mathrm{Init}(s) := \sum_{i \in I_0}\nu_0(t_i) \cdot \mathrm{Lm}(\rho_0(\ell_i)).
    \end{align*}
    Observe that if $s \in \boldsymbol{k}[G_{X}\Lambda_{X}]$ and $\log E'(x)^b$ appears in $s$, then $\mathrm{Init}(s)$ is only (possibly) affected by the leading monomial of $\rho_0(\log E'(x)^b)$.
\end{rmk}

\begin{lemma}\label{rho_0 is injective}
    Suppose $X$, and $\boldsymbol{k}$ satisfy the order and separation assumptions in Remark \ref{order props}. 
    Then $\rho_0$ is injective.
\end{lemma}

\begin{proof}
Let $s \in \boldsymbol{k}[G_{X}\Lambda_{X}]$.
If $\mathrm{Init}(s) \ne 0$, then it follows from the definition of $\mathrm{Init}(s)$ in Remark \ref{initial subsum of rho(s)} that $\mathrm{Init}(s)$ is an initial subsum of $\rho_0(s)$.
Thus
    $$\mathrm{Lm}\big(\rho_0(s)\big) = \mathrm{Lm}\big(\mathrm{Init}(s)\big) \ne 0$$
which shows $\rho_0$ is injective.
Like the proof of Lemma \ref{series comp}, we will find a monomial with nonzero coefficient in $\mathrm{Init}(s)$.

Write $s = s_1\ell_1 + \cdots + s_n\ell_n \in \kk[G_{X}\Lambda_{X}]$. 
Following the setup of Remark \ref{initial subsum of rho(s)}, we may assume without loss of generality that $\alpha_{i,j} = \beta_{i,j} + b_{i,j} = 0$ for all $i \in I_0$, $j = 1,\dots,p$, by rescaling or factoring out common terms.
Let $q_i$ be the (finite) number of terms in $t_i$, and let $c_{i,l} \ne 0$ be the coefficient of the $l$th monomial of $t_i$.
Then we can write 
\begin{multline*}
    \mathrm{Init}(s) = \sum_{i \in I_0}\sum_{l=1}^{q_i} c_{i,l} \prod_{j=1}^p 
    E_1(x_j-1)^{\beta_{i,j} + m_{i,j,l,2}} (1 + \epsilon_2(x_j))^{m_{i,j,l,2}} \\
    E_2(x_j-1)^{m_{i,j,l,3}} (1 + \epsilon_3(x_j))^{m_{i,j,l,3}} \cdots E_{d-1}(x_j-1)^{m_{i,j,l,d}} (1 + \epsilon_d(x_j))^{m_{i,j,l,d}}.
\end{multline*}
Note that $\beta_{i,j}$ depends on $i$ and $j$, but not on $l$.
Let 
\begin{multline*}
    \mathrm{Init}(s)_{i,l}:= \prod_{j=1}^p 
    E_1(x_j-1)^{\beta_{i,j} + m_{i,j,l,2}} (1 + \epsilon_2(x_j))^{m_{i,j,l,2}} \\
    E_2(x_j-1)^{m_{i,j,l,3}} (1 + \epsilon_3(x_j))^{m_{i,j,l,3}} \cdots E_{d-1}(x_j-1)^{m_{i,j,l,d}} (1 + \epsilon_d(x_j))^{m_{i,j,l,d}}.
\end{multline*}
If $\mathrm{Supp}\left(\mathrm{Init}(s)_{i_1,l_1}\right) \cap \mathrm{Supp}\left(\mathrm{Init}(s)_{i_2,l_2}\right) \ne \varnothing$, then $\beta_{i_1,j}$ and $\beta_{i_2,j}$ must be in the same $\mathbb{Z}$-orbit for all $j=1,\dots,p$.
So without loss of generality, we may assume $\beta_{i,j} \in \mathbb{N}$ for all $i \in I_0$, $j=1,\dots,p$.

We will now proceed by induction on $d$.
If $d = 1$, then $\mathrm{Init}(s)$ has just one term, and its coefficient is nonzero, so we are done.
The base case of our induction will be $d=2$.

If $d = 2$, then
    $$\mathrm{Init}(s) = \sum_{i \in I_0}\sum_{l=1}^{q_i} c_{i,l} \prod_{j=1}^p 
    E_1(x_j-1)^{\beta_{i,j} + m_{i,j,l,2}} (1 + \epsilon_2(x_j))^{m_{i,j,l,2}}.$$
Recall that $\epsilon_2(x_j) = \frac{E_2(x_j-1)}{E_1(x_j-1)}$.
By Lemma \ref{exponent of E_d in epsilon_d}, for any $i,l$, the sum of the exponents of $E_1(x_j-1)$ and $E_2(x_j-1)$ in any monomial of $\mathrm{Init}(s)_{i,l}$ is $\beta_{i,j} + m_{i,j,l,2}$.
If this value differs for some pairs $i_1,l_1$ and $i_2,l_2$ and some $j=1,\dots,p$, then
    $$\mathrm{Supp} (\mathrm{Init}(s)_{i_1,l_1}) \cap \mathrm{Supp} (\mathrm{Init}(s)_{i_2,l_2}) = \varnothing.$$
In particular, since $\beta_{i,j}$ does not depend on $l$, we have
    $$\mathrm{Supp} (\mathrm{Init}(s)_{i,l_1}) \cap \mathrm{Supp} (\mathrm{Init}(s)_{i,l_2}) = \varnothing$$
for all $i \in I_0$ and $l_1 \ne l_2$.

So it will suffice to find a monomial with a nonzero coefficient among the sum
    $$S := \sum_{i \in I_1} c_{i,1} \prod_{j=1}^p E_1(x_j-1)^{\beta_{i,j} + m_{i,j,1,2}} (1 + \epsilon_2(x_j))^{m_{i,j,1,2}}$$
where $I_1 \subset I_0$ is set of indices at which
    $$(\beta_{i,1} + m_{i,1,1,2}, \dots, \beta_{i,p} + m_{i,p,1,2})$$
is maximized in the lexicographic order on $\mathbb{N}^p$.
Let $N_j := \beta_{i,j} + m_{i,j,l,2}$ for any $i \in I_1$.
Without loss of generality, we may assume that for each $j=1,\dots,p$, there are some $i_1,i_2 \in I_1$ such that $\beta_{i_1,j} \ne \beta_{i_2,j}$, since if $\beta_{i,j}$ is fixed across all $i \in I_1$ for some $j$, then so is $m_{i,j,1,2}$, and we can just factor out
    $$E_1(x_j-1)^{\beta_{i,j} + m_{i,j,1,2}} (1 + \epsilon_2(x_j))^{m_{i,j,1,2}}$$
from each term in $S$ without affecting the sign of $S$.

To simplify notation for the rest of the $d=2$ case, we drop the second index from the coefficients $c_{i,1}$ and the last two indices from the exponents $m_{i,j,1,2}$.

If $|I_1| = 1$, then we are done, since $S$ would be a sum of a single monomial with a nonzero coefficient.
So suppose $|I_1| > 1$.
We will now find a monomial of $S$ with nonzero coefficient.
We will do this by inductively eliminating terms from $S$ until we are left with a single monomial.
\begin{enumerate}
    \item Let $j_1$ be the smallest index for which there are $i_1,i_2 \in I_1$ such that $m_{i_1,j_1} \ne m_{i_2,j_1}$.
    Such an index exists because if $i_1 \ne i_2$, then there is some $j$ such that $\beta_{i_1,j} \ne \beta_{i_2,j}$ (since we assumed $g_{i_1} \ne g_{i_2}$ for $i_1 \ne i_2$), and thus $m_{i_1,j} \ne m_{i_2,j}$.
    
    \item Let $n_1 = \max_{i \in I_1} (m_{i,j_1})$.
\end{enumerate}
Then for $k = 0,\dots,n_1$, the coefficient of 
    $$\left(\prod_{j=1}^p E_1(x_j-1)^{N_j}\right) \cdot \epsilon_2(x_{j_1})^k$$
is
    $$\sum_{i \in I_1} c_i \binom{m_{i,j_1}}{k}.$$
If the $m_{i,j_1}$ are all distinct for $i \in I_1$, then we are done because these coefficients cannot all be zero, as the $c_i$'s are nonzero.

If the $m_{i,j_1}$ are not all distinct, then we proceed by induction:
Given $I_{l}$, $j_l$, and $n_l$ with $m_{i,j_l} \ne n_l$ for some $i \in I_l$, we define $I_{l+1}$, $j_{l+1}$, and $n_{l+1}$ as follows:
\begin{enumerate}
    \item Let $I_{l+1} := \{i \in I_l : m_{i,j_l} = n_l\}$. 
    Then $|I_{l+1}| < |I_l|$ by our choice of $j_l$.
    
    \item Let $j_{l+1}$ be the smallest index for which there are $i_1,i_2 \in I_{l+1}$ such that $m_{i_1,j_{l+1}} \ne m_{i_2,j_{l+1}}$, which exists by the same reasoning as above.
    
    \item Let $n_{l+1} = \max_{i \in I_{l+1}} (m_{i,j_{l+1}})$.
\end{enumerate}
Then for each $k = 1,\dots,n_{l+1}$, the coefficient of 
    $$\left(\prod_{j=1}^p E_1(x_j-1)^{N_j}\right) \cdot \Big(\epsilon_2(x_{j_1})^{n_1} \cdots \epsilon_2(x_{j_l})^{n_l}\Big) \cdot \epsilon_2(x_{j_{l+1}})^k$$
is
    $$\sum_{i \in I_{l+1}} c_i\binom{n_1}{n_1} \cdots \binom{n_l}{n_l} \cdot \binom{m_{i,j_{l+1}}}{k} = \sum_{i \in I_{l+1}} c_i \binom{m_{i,j_{l+1}}}{k}.$$
If $|I_{l+1}| = 1$, or if the values of $m_{i,j_{l+1}}$ for $i \in I_{l+1}$ are all distinct, then we find a nonzero coefficient.
Since $|I_{l+1}| < |I_{l}| < \cdots < |I_0| \le n$, this process must terminate.
So we find a term of $S$ with nonzero coefficient, which means $\mathrm{Init}(s) \ne 0$.
This proves the Lemma in the $d=2$ case.

Now, we will show that if we can find a nonzero term of $\Init(s)$ for any $s$ with largest order of derivative $d-1$, then we can also find a nonzero term of $\Init(s)$ if the largest order of derivative in $s$ is $d$.
Let 
    $$(\mu_{1},\dots,\mu_{p}) := \max \{(m_{i,1,l,d}, \dots, m_{i,p,l,d}) : i \in I_0, l \in \mathbb{N}\}$$
in the lexicographic order on $\mathbb{N}^p$.
Let $\mathcal{I}$ be the set of pairs $(i,l)$ at which this maximum is achieved, i.e., at which $(m_{i,1,l,d}, \dots, m_{i,p,l,d}) = (\mu_{1},\dots,\mu_{p})$.

Recall that 
\begin{multline*}
    \mathrm{Init}(s)_{i,l} = \prod_{j=1}^p 
    E_1(x_j-1)^{\beta_{i,j} + m_{i,j,l,2}} (1 + \epsilon_2(x_j))^{m_{i,j,l,2}} 
    \\
    E_2(x_j-1)^{m_{i,j,l,3}} (1 + \epsilon_3(x_j))^{m_{i,j,l,3}} \cdots  E_{d-1}(x_j-1)^{m_{i,j,l,d}} (1+\epsilon_d(x_j))^{m_{i,j,l,d}}.
\end{multline*}
By Lemma \ref{exponent of E_d in epsilon_d}, the only place $E_d(x_j-1)$ appears in this expression is in $\epsilon_d(x_j)$. 
Furthermore, $E_d(x_j-1)$ can only ever appear with exponent 1 in monomials of $\epsilon_d(x_j)$.
Let $\epsilon_{d*}(x_j)$ be the subsum of $\epsilon_d(x_j)$ with monomials in which $E_d(x_j-1)$ appears.

Let $S_d$ be the subsum of $\mathrm{Init}(s)$ such that for every monomial in $\mathrm{Supp} (S_d)$ and for each $j=1,\dots,p$, the exponent of $E_d(x_j-1)$ is $\mu_j$.
Then we can write
\begin{align*}
    S_d = \sum_{(i,l) \in \mathcal{I}} c_{i,l} \prod_{j=1}^p 
    &E_1(x_j-1)^{\beta_{i,j} + m_{i,j,l,2}} (1 + \epsilon_2(x_j))^{m_{i,j,l,2}} \\
    & E_2(x_j-1)^{m_{i,j,l,3}} (1 + \epsilon_3(x_j))^{m_{i,j,l,3}} \cdots E_{d-2}(x_j-1)^{m_{i,j,l,d-1}} (1+\epsilon_{d-1}(x_j))^{m_{i,j,l,d-1}} \\ &E_{d-1}(x_j-1)^{\mu_{j}} \epsilon_{d*}(x_j)^{\mu_{j}}.
\end{align*}
By definition, the last line $E_{d-1}(x_j-1)^{\mu_{j}} \cdot \epsilon_{d*}(x_j)^{\mu_{j}}$ above must be the same for all $(i,l) \in \mathcal{I}$, so we can factor it out to write
    $$S_d = \Bigg(\prod_{j=1}^p E_{d-1}(x_j-1)^{\mu_{j}} \cdot  \epsilon_{d*}(x_j)^{\mu_{j}}\Bigg) \cdot S_d'$$
where 
\begin{align*}
    S_d' = \sum_{(i,l) \in \mathcal{I}} c_{i,l} \prod_{j=1}^p 
    &E_1(x_j-1)^{\beta_{i,j} + m_{i,j,l,2}} (1 + \epsilon_2(x_j))^{m_{i,j,l,2}} \\
    & E_2(x_j-1)^{m_{i,j,l,3}} (1 + \epsilon_3(x_j))^{m_{i,j,l,3}} \cdots  E_{d-2}(x_j-1)^{m_{i,j,l,d-1}} (1+\epsilon_{d-1}(x_j))^{m_{i,j,l,d-1}}.
\end{align*}
By induction, we can find a monomial of $S_d'$ with nonzero coefficient, which gives a monomial of $\mathrm{Init}(s)$ with nonzero coefficient.
\end{proof}

\begin{rmk}\label{approx matches germs}
    Lemma \ref{rho_0 is injective} shows that $\rho_0$ is injective, so the order on $\boldsymbol{k}((H_{X-1}))$ induces an order on $\boldsymbol{k}[G_{X}\Lambda_{X}]$.
    However, we need the induced order to be compatible with how we intend to interpret elements of $\boldsymbol{k}[G_{X}\Lambda_{X}]$ as germs of functions.
    Because $\rho_0$ is defined using the intuition that we can ``approximate" $\log E'(x) = E(x-1) + \log E'(x-1)$ by $E_0(x-1) + E_1(x-1)$, it may not be immediately clear that the order induced by $\boldsymbol{k}((H_{X-1}))$ via $\rho_0$ is the one we want.
    
    If $s \in \boldsymbol{k}[G_{X}\Lambda_{X}]$ and $\log E'(x)^b$ appears in $s$, then by the way $\mathrm{Init}(s)$ is defined in Remark \ref{initial subsum of rho(s)}, the only term of $\rho_0\left(\log E'(x)^b\right)$ that has any possibility of contributing to $\mathrm{Init}(s)$ is $\mathrm{Lm}\left(\rho_0\left(\log E'(x)^b\right)\right)$. 
    The ``error" in the approximation $\rho_0\left(\log E'(x)^b\right)$ of $\log E'(x)^b$ appears in every term \textit{except the leading monomial}.
    So the sign of $\Init(s)$ is not distorted by the approximation.
\end{rmk}

\begin{cor}\label{log series comp}
    $\Lambda_{X}$ is ordered lexicographically on the exponents of its generators, i.e., if $x_1 > x_2 > \cdots > x_p \in X$ and $b_1,b_2, \dots,b_p \in \boldsymbol{k}^{\times}$, then
        $$\prod_{j=1}^p \log E'(x_j)^{b_j} > 1$$
    if and only if $b_1 > 0$. 
    So if $s \in \boldsymbol{k}[\Lambda_{X}]$, then $s > 0$ if and only if the coefficient of its largest monomial is positive.
\end{cor}

\begin{proof}
    Tracing through the proof of Lemma \ref{rho_0 is injective}, we first want to determine the sign of 
        $$\prod_{j=1}^p \log E'(x_j)^{b_j} - 1.$$
    If $b_1 > 0$, then $\mathrm{Init}(s) = \prod_{j=1}^p E_0(x_j-1)^{b_j} > 0$, and if $b_1 < 0$, then $\mathrm{Init}(s) = -1 < 0$.
    
    If $g_1 > g_2 > \dots > g_n \in \Lambda_{X}$ and $c_1,c_2,\dots,c_n \in \boldsymbol{k}^{\times}$, then we want to determine the sign of 
        $$s = \sum_{i=1}^n c_ig_i.$$
    In this case, $\mathrm{Init}(s) = c_1\mathrm{Lm}(\rho_0(g_1))$, so the sign of $s$ is determined by the sign of $c_1$.
\end{proof}

\section{A logarithmic-exponential series field constructed from \texorpdfstring{$E$}{E}-monomials}

In this section, we adapt the construction of the logarithmic-exponential series in \cite{LEseries} to build a logarithmic-exponential series field starting with monomials involving $E$ and its derivatives. 
The first part of the usual logarithmic-exponential series construction begins with the multiplicative group $x^{\boldsymbol{k}}$ of monomials at the zeroth stage and inductively adds new monomials for increasing levels of exponentiation to end up with an exponential field $\boldsymbol{k}((x^{-1}))^e$.
The second part of the construction uses an embedding $\varphi : \boldsymbol{k}((x^{-1}))^e \to \boldsymbol{k}((x^{-1}))^e$ such that every element in the image of $\varphi$ has a logarithm in $\boldsymbol{k}((x^{-1}))^e$.
If this approach is adapted naively to monomials of the form 
    $$\prod_{j=1}^p E(x_j)^{a_{j,0}} E'(x_j)^{a_{j,1}} \cdots E^{(d)}(x_j)^{a_{j,d}}$$
for $x_1 > \cdots > x_p \in X$ and $a_{j,l} \in \boldsymbol{k}$,
three problems arise.
\begin{enumerate}
    \item First, a logarithm of $E^{(d)}(x)$ for $d \ge 1$ would not arise naturally from the first part of the construction.
    
    \item Second, some monomials of the form above are ``small" relative to others.
    For example, as functions we have
        $$\exp\left(\frac{E'(x)^{1/2}}{E(x)^{1/2}}\right) = \exp\left(E'(x-1)^{1/2}\right) < E(x)$$
    since $E'(x-1)^{1/2} < E(x-1)$.
    So it does not make sense to add $\exp\left(\frac{E'(x)^{1/2}}{E(x)^{1/2}}\right)$ as a new monomial over a field of coefficients that contains $E(x)$.
    
    \item Third, it does not always make sense to take infinite sums of monomials that are ``close" in the sense that they differ by a factor of something ``small".
    For example, 
        $$\frac{E'(x)^2}{E(x)E''(x)} = \frac{1}{1 + \frac{E''(x-1)}{E'(x-1)^2}} < 1$$
    but it does not make sense to sum the infinite reverse well-ordered family 
        $$\left\{\left(\frac{E'(x)^2}{E(x)E''(x)}\right)^n : n \in \mathbb{N}\right\}$$
    because each of these monomials is approximately 1.
\end{enumerate}

To fix these issues, we include monomials for $\log E'$ and for exp of some ``small" infinite monomials at the very beginning of the first part of the construction.
Then at later stages in the Part 1 inductive construction, we add new monomials only for exp applied to ``large" monomials.
We also only ever allow finite sums of monomials whose quotient is ``small."
The result is that in our adaptation of the first part of the logarithmic-exponential series construction, we build a ring with a partially defined exponential function.
The second part of our construction is similar to the second part of the construction in \cite{LEseries} in that we build embeddings between countably many partial exponential rings constructed as in part one and then show that every element in each of these rings eventually has a multiplicative inverse, a logarithm, and an exponential under a finite sequence of embeddings.

\subsection{Part 1: Building partial exponential rings} \label{subsection 4.1}
Let $X$ and $\boldsymbol{k}$ satisfy the ordering and separation assumptions in Remark \ref{order props}, so that $\boldsymbol{k}[G_{X}\Lambda_{X}]$ is totally ordered by Lemma \ref{rho_0 is injective}.
Assume also that $\boldsymbol{k} \vDash T_{\mathrm{an}}(\exp,\log)$.

We start this subsection by pinning down the zeroth stage group of monomials, which we call $\Gamma_{0}$, from which we will build a partial exponential ring.
We want $\Gamma_{0}$ to include $G_{X}\Lambda_{X}$ along with new monomials to represent exp of the ``small" positive purely infinite elements of $\boldsymbol{k}[G_{X}\Lambda_{X}]$.
We illustrate the issue with another example:
\begin{ex}
    Treating $E(x)$ as a function, we can compute that
        $$\exp\left(\frac{E''(x)}{E'(x)}\right) > \exp\left(\frac{E'(x)}{E(x)}\right) = \exp\big(E'(x-1)\big) > E(x)^n$$
    for any $n \in \mathbb{N}$, since $E'(x-1) > nE(x-1)$.
    However, 
        $$\exp\left(\frac{E''(x)}{E'(x)} - \frac{E'(x)}{E(x)}\right) = \exp \left( \frac{E''(x-1)}{E'(x-1)}\right) < E(x).$$
    So even though the ``small" monomial issue does not arise with $\exp\left(\frac{E''(x)}{E'(x)}\right)$ or $\exp\left(\frac{E'(x)}{E(x)}\right)$ separately, it does not makes sense to take a group containing both of these expressions as monomials over a field of coefficients that contains $E(x)$.
\end{ex}

The examples point toward the following definition of which elements of $G_{X}$ are ``small" enough to cause a problem.
\begin{defn}\label{small defn}
    Let $g \in G_{X}$ and write
        $$g = \prod_{j=1}^pE(x_j)^{a_{j,0}} \cdot E'(x_j)^{a_{j,1}} \cdots E^{(d)}(x_j)^{a_{j,d}}$$
    with $p,d \in \mathbb{N}$, $x_1 > \cdots > x_p \in X$, and $a_{j,0},\dots,a_{j,d} \in \boldsymbol{k}$.
    Let $\xi_{x_j} := a_{j,0} + \cdots + a_{j,d}$, the sum of exponents of generators defined from $x_j$.
    If $x \in X$ does not appear in $g$, let $\xi_x := 0$.
    We say $g$ is small if $\xi_{x} = 0$ for all $x \in X$.
\end{defn}
Note that $1$ is small, and if $g_1,g_2 \in G_{X}$ are both small, then so is $g_1g_2$.
So the small elements of $G_{X}$ form a subgroup.

Since it does not always make sense to add new monomials for $\exp(g)$ with small infinite $g \in \boldsymbol{k}[G_{X}\Lambda_{X}]$ at later stages in the inductive Part 1 construction, we want to include such monomials from the very beginning.
But the examples above show that elements of the form $\exp(g)$ for small $g$ would intersperse with elements of $G_{X}\Lambda_{X}$, and that could make it challenging to figure out how to order even the monomials we start with.
To prevent the ordering on our initial monomial group from becoming too complicated, we will only include monomials $\exp(g)$ for certain relatively simple small infinite $g \in G_{X}\Lambda_{X}$.
Define
\begin{align*}
    T_{X} := \Bigg\langle &c \prod_{j=1}^p\log E'(x_j)^{a_{j,-1}}\frac{E''(x_j)^{a_{j,2}} \cdots E^{(d)}(x_j)^{a_{j,d}}}{E'(x_j)^{a_{j,2} + \cdots + a_{j,d}}} : c \in \boldsymbol{k}; p,d, a_{j,l} \in \mathbb{N}; x \in X; \\
        &\exists j,l \big(a_{j,l} \ne 0\big); \big(a_{j,-1} = 1 \big)\rightarrow \exists k \ne j, l \ne -1 \big(a_{k,l} \ne 0\big) \Bigg\rangle 
\end{align*}
where $T_{X}$ is generated additively. 
The final two conditions ensure that $\boldsymbol{k} \cap T_{X} = \{0\}$ and $c\log E'(x) \not \in T_{X}$ for any $c \in \boldsymbol{k}$, $x \in X$, since we already know how to exponentiate such elements---we assumed $\boldsymbol{k}$ is an exponential field, and we want to define exp of $c\log E'(x)$ to be $E'(x)^c$.

Since $X$ and $\boldsymbol{k}$ satisfy the ordering and separation assumptions in \ref{order props}, $T_{X}$ is totally ordered because it is a subgroup of $\boldsymbol{k}[G_{X}\Lambda_{X}]$.
Let $e_{T_{X}}(T_{X})$ be a multiplicative copy of $T_{X}$ via an order preserving isomorphism $e_{T_{X}} : T_{X} \to e_{T_{X}}(T_{X})$.

\begin{rmk}
    For consistency of notation, we will always express elements of $T_{X}$ as sums of the standard monomials of the group ring $\boldsymbol{k}[G_{X}\Lambda_{X}]$ instead of as sums of the generators of $T_{X}$ itself.
\end{rmk}

We can now define the group of monomials from which we will start the inductive construction of Part 1:
    $$\Gamma_{0} := G_{X} \Lambda_{X} e_{T_{X}}(T_{X}).$$

\begin{rmk}
    The 0 subscript in $\Gamma_{0}$ indicates the zeroth stage of the inductive Part 1 construction.
    Later on, we will need to specify the sets $X$ and $\kk$ from which $\Gamma_0$ is built. 
    When this is necessary, we will include $X$ and $\kk$ as additional subscripts and write $\Gamma_{X,\kk,0}$.
    For now, we suppress the additional subscripts for clarity.
\end{rmk}

Next, we will extend the order on $G_{X}\Lambda_{X}$ to $\Gamma_{0}$.
In fact, we will see that there is a unique way to order $\Gamma_{0}$ consistent with the order defined in Lemma \ref{series comp} and the way a logarithm must interact with $E$ and its derivatives as a consequence of $T_{\transexp}$.
We would like to define this order on $\Gamma_{0}$ by
\begin{align*}
    g \cdot \ell > e_{T_X}(t) \text{ if and only if } \log (g \cdot \ell) > t
\end{align*}
and show the comparison of $\log(g \cdot \ell)$ and $t$ is determined by Lemma \ref{series comp} and the difference equation for $E$.
However, we cannot yet make sense of ``$\log(g \cdot \ell)$". 
(Later in this section, we will define $\log(g \cdot \ell)$ using an infinite sum.)
Instead, we will define a temporary stand-in  
    $$\widehat{\log} : \rho_0(G_{X}\Lambda_{X}) \to \boldsymbol{k}[H_{X-1}]$$
to approximate the way the logarithm must later be defined.
We will then extend the order on $G_{X} \Lambda_{X}$ to an order on $\Gamma_{0}$ by 
\begin{align*}
    g \cdot \ell > e_{T_{X}}(t) \text{ if and only if } \widehat{\log} (\rho_0(g \cdot \ell)) > \rho_0(t)
\end{align*}
for $t \ne 0$.

How should we define $\widehat{\log}$ so that it approximates the behavior of log?
As functions, $\log$ and the derivatives of $E$ work as follows:
\begin{align*}
    \log \big(E^{(d)}(x)^a\big) 
        &= \log \big(E(x)^a B_d(x-1)\big) \\
        &= aE(x-1) + a\log B_d(x-1) \\
        &= aE(x-1) + ad\log E'(x-1) + a\log\left(\frac{B_d(x-1)}{E'(x-1)^d}\right) \\
        &= aE(x-1) + adE(x-2) + ad\log E'(x-2) + \cdots \\
    \log \big(\log E'(x)^b\big) 
        &= b\log\big(E(x-1) + \log E'(x-2)\big) \\
        &= bE(x-2) + b \log\left(1 + \frac{\log E'(x-2)}{E(x-1)}\right) \\
        &= bE(x-2) + \cdots
\end{align*}
where the ``$\cdots$" consists of infinitesimal expressions. 
The largest term above is $aE(x-1)$, and it appears in $\log \big(E^{(d)}(x)^a\big)$ but not in $\log \big(\log E'(x)^b\big)$.
We will define $\widehat{\log}$ so that it matches the computations above up to terms of the form $aE(x-1)$ and ignores anything smaller.

\begin{defn}
    If $x_1 > \cdots > x_p$ are the elements of $X$ appearing in $g \cdot \ell \in G_{X}\Lambda_{X}$ and $d$ is the largest derivative appearing in $g$, then express $\Lt(\rho_0(g \cdot \ell)) = \Lm(\rho_0(g \cdot \ell))$ as
        $$\prod_{j=1}^p E_0(x_j)^{\alpha_j}E_0(x_j-1)^{\beta_j}E_1(x_j-1)^{\beta_j+m_1}E_2(x_j-1)^{m_2} \cdots E_d(x_j-1)^{m_d}$$
    with $\alpha_j,\beta_j \in \kk$ and $m_1,\dots,m_d \in \NN$, for $j=1,\dots,p$.
    Define
        $$\widehat{\log}(\rho_0(g \cdot \ell)) := \alpha_1 E_0(x_1-1) + \cdots + \alpha_p E_0(x_p-1)$$
    which may be 0 if $\alpha_1 = \cdots = \alpha_p = 0$, i.e., if no generator of the form $E_0(x)$ appears in the leading monomial of $\rho_0(g \cdot \ell)$.
\end{defn}

\begin{rmk}\label{why log hat is good enough}
We will only ever use $\widehat{\log}$ to define an ordering of $\Gamma_{0}$.
To justify that $\widehat{\log}$ is a good enough approximation of log for this purpose, we will show that
\begin{enumerate}
    \item if $g$ is small and $t \ne 0$, then Lemma \ref{series comp} and the difference equation for $E$ force $e_{T_{X}}(t) > g \cdot \ell$;
    
    \item if $g$ is not small, then the ``dominant" terms of $\log(g \cdot \ell)$ and $\widehat{\log}(\rho_0(g \cdot \ell))$ are the same (we make this precise in Corollary \ref{log hat = Lm(log)}, after defining the full logarithm); and
    
    \item if $g$ is not small, then only the ``dominant" leading term of $\widehat{\log}(\rho_0(g \cdot \ell))$ matters in determining how it compares to $\rho_0(t)$.
\end{enumerate}
We will obtain (1) and (3) as consequences of Lemmas \ref{Leading monomial has E_1} and \ref{not of the form E_0(x-1)}.
\end{rmk}

\begin{lemma}\label{Leading monomial has E_1}
    Suppose $X$ and $\boldsymbol{k}$ satisfy the ordering and separation assumptions in Remark \ref{order props}. 
    If $t \in T_{X} \cap \boldsymbol{k}[G_{X}]$, then the leading term of $\sigma_0(t)$ has a positive integer power of $E_1(x)$ as a generator, for some $x \in X$.
\end{lemma}

\begin{proof}
    Let $t = \sum_{i=1}^n c_ig_i \in T_{X} \cap \boldsymbol{k}[G_{X}]$ with $c_1,\dots,c_n \ne 0$.
    Let $x_1 > x_2 > \dots > x_p$ list the elements of $X$ appearing in $t$, and let $d$ be the largest derivative appearing in any $g_i$.
    Then we can write
    \begin{align*}
        g_i = \prod_{j=1}^p&\frac{E''(x_j)^{a_{i,j,2}} \cdots E^{(d)}(x_j)^{a_{i,j,d}}}{E'(x_j)^{a_{i,j,2} + \cdots + a_{i,j,d}}} \\
        \sigma_0(g_i) = \prod_{j=1}^p & E_0(x_j)^{a_{i,j,2} + \cdots + a_{i,j,d} - (a_{i,j,2} + \cdots + a_{i,j,d})} E_1(x_j)^{2a_{i,j,2} + \cdots + da_{i,j,d} - (a_{i,j,2} + \cdots + a_{i,j,d})}  \\
            &\left(\sum_{k_2=0}^{a_{i,j,2}}\binom{a_{i,j,2}}{k}\left(\frac{E_2(x_j)}{E_1(x_j)}\right)^{k_2}\right) \cdots \left(\sum_{k_{d}=0}^{a_{i,j,d}}\binom{a_{i,j,d}}{k_d}\left(\cdots + \frac{E_2(x_j) \cdots E_{d}(x_j)}{E_1(x_j)^{d-1}}\right)^{k_{d}}\right).
    \end{align*}
    Note that $\sigma_0(g_i)$ is a finite sum because $a_{i,j,2},\dots,a_{i,j,d}$ are natural numbers for each $j=1,\dots,p$.
    In every monomial of $\sigma_0(g_i)$, the exponent of $E_0(x_j)$ is 0 and the exponent of $E_1(x_j)$ is a natural number for all $j=1,\dots,p$.
    
    If $d = 2$, then the leading monomial of $\sigma_0(t)$ is
        $$\prod_{j=1}^p E_1(x_j)^{2a_{i,j,2}}$$
    where $i$ is the single index at which $\overline{a_{i,2}}$ is maximized.
    So we may assume $d > 2$.
    
    We will trace through the proof of Lemma \ref{series comp} to show that we find a monomial with $E_1(x_{j_*})$ as a generator for some $j_*=1,\dots,p$.
    Let $I$ be as in Lemma \ref{series comp}, and let 
        $$A_j = a_{i,j,2} + 2a_{i,j,3} + \cdots + (d-1)a_{i,j,d}$$
    for all $i \in I$.
    Without loss of generality we may assume $I = \{1,\dots,n\}$.
    In Lemma \ref{series comp}, we proceeded with the argument using some fixed $j$. 
    Here we will use $j=1$.

    First suppose $a_{1,1,d},\dots,a_{n,1,d} \in \mathbb{N}$ are all distinct.
    By possibly renumbering $g_1,\dots,g_n$, we may assume $a_{1,1,d} < \dots < a_{n,1,d}$.
    Note that this means $a_{i,1,d} \ge i-1$.
    The coefficient $\sum_{i =1}^n c_i \binom{a_{i,1,d} }{k_d}$ of
        $$\left(\frac{E_2(x_1) \cdots E_{d}(x_1)}{E_1(x_1)^{d-1}}\right)^{k_{d}}\prod_{j=1}^p E_1(x_j)^{A_j}$$ 
    cannot be 0 for each $k_d = 0,\dots,n-1$ because then we would have $c_1 = \cdots = c_n = 0$.
    So some coefficient must be nonzero for some $k_* \le n - 1 \le a_{n,1,d}$.
    If $p > 1$ or if $(d-1)k_* < A_1$ then this monomial has some $E_1(x_{j_*})$ as a generator, and we are done.
    So suppose $p = 1$ and $(d-1)k_* = A_1$.
    
    Since $p = 1$, we remove the index $j$ counting $1,\dots,p$ in the following computations to lighten notation.
    
    Note that $(d-1)k_* = A = a_{i,2} + 2a_{i,3} + \cdots + (d-1)a_{i,d}$ for all $i = 1,\dots,n$.
    Since this holds for $i=n$, and since $a_{n,d} \ge n-1$, we must have the following:
    \begin{enumerate}
        \item $k_* = n-1 = a_{n,d}$ is the smallest value of $k_d$ such that $\sum_{i =1}^n c_i \binom{a_{i,d} }{k_d} \ne 0$, and $0 \ne \sum_{i = 1}^n c_i \binom{a_{i,d} }{k_d} = c_n$
        
        \item $a_{n,2} = 0 = a_{n,3} = \dots = a_{n,d-1} = 0$.
    \end{enumerate}
    So $A = (d-1)(n-1)$.
    Observe also that since $a_{1,d} < \dots < a_{n,d}$, we must have $a_{i,d} = i-1$.
    In particular $a_{n-1,d} = n-2$.
    Since 
        $$(d-1)(n-1) = A = a_{n-1,2} + \cdots + (d-2)a_{n-1,d-1} + (d-1)(n-2)$$
    we have that $a_{n-1,d-1} \le \frac{d-1}{d-2} \le 2$.
    
    We will use the fact that $0 \le a_{i,d-1} \le 2$ to show that the coefficient of the following monomial with $E_1(x)$ as a generator is nonzero:
        $$M := \left(\frac{E_2(x) \cdots E_{d}(x)}{E_1(x)^{d-1}}\right)^{n-2} \left(\frac{E_2(x) \cdots E_{d-1}(x)}{E_1(x)^{d-2}}\right) E_1(x)^{A}.$$ 
    Observe that $\frac{E_2(x) \cdots E_{d-1}(x)}{E_1(x)^{d-2}}$ appears only in $\left(\cdots + \frac{E_2(x) \cdots E_{l}(x)}{E_1(x)^{l-1}}\right)$ for $l=d-1,d$.
    Its coefficient is 1 if $l = d-1$, and its coefficient is $d$ if $l=d$.
    So the coefficient of $M$ is
        $$\sum_{i=1}^n c_i \Bigg(\binom{a_{i,d} }{n-2} \binom{a_{i,d-1}}{1} + \binom{a_{i,d} }{n-1} \binom{n-1}{1}d\Bigg).$$
    The first term of the summand comes from building $M$ as a product of the monomials 
    \begin{align*}
        \left(\frac{E_2(x) \cdots E_{d}(x)}{E_1(x)^{d-1}}\right)^{n-2} &\text{ from } \binom{a_{i,d} }{n-1}\left(\cdots + \frac{E_2(x) \cdots E_{d}(x)}{E_1(x)^{d-1}}\right)^{n-2} \\
        \left(\frac{E_2(x) \cdots E_{d-1}(x)}{E_1(x)^{d-2}} \right)
        &\text{ from } \binom{a_{i,d-1}}{1}\left(\cdots + \frac{E_2(x) \cdots E_{d-1}(x)}{E_1(x)^{d-2}}\right)^{1}.
    \end{align*}
    The second term of the summand comes from seeing $M$ as a monomial of
    \begin{align*}
        \binom{a_{i,d} }{n-1}\left(\cdots + d\frac{E_2(x) \cdots E_{d-1}(x)}{E_1(x)^{d-2}} + \cdots + \frac{E_2(x) \cdots E_{d}(x)}{E_1(x)^{d-1}}\right)^{n-1}.
    \end{align*}
    Since $a_{i,d} = i-1$, most of the terms in the sum are 0.
    We can calculate
    \begin{align*}
        \binom{a_{i,d} }{n-2} \binom{a_{i,d-1}}{1} + \binom{a_{i,d} }{n-1} \binom{n-1}{1}d &= 0 \text{ for } i=1,\dots,n-1 \\
        \binom{a_{n-1,d}}{n-2} \binom{a_{n-1,d-1}}{1} + \binom{a_{n-1,d}}{n-1} \binom{n-1}{1}d &\le 1 \cdot 2 + 0 \cdot(n-1)d = 2 \\
        \binom{a_{n,d}}{n-2} \binom{a_{n,d-1}}{1} + \binom{a_{n,d}}{n-1} \binom{n-1}{1}d &= (n-1)\cdot 0 + 1 \cdot (n-1)d = (n-1)d.
    \end{align*}
    Altogether, the coefficient of $M$ is $\delta \cdot c_{n-1} + (n-1)d \cdot c_n$ with $\delta \in \{0,1,2\}$.
    Since 
        $$0 = \sum_{i =1}^n c_i \binom{a_{i,d} }{n-2} = c_{n-1} + (n-1)c_n$$
    and $d > 2$, we must have $\delta \cdot c_{n-1} + (n-1)d \cdot c_n \ne 0$.
    This finishes the proof in the case that $a_{1,d},\dots,a_{n,d} \in \mathbb{N}$ are all distinct.
    
    Suppose $a_{i,1,d} \in \mathbb{N}$ are not all distinct for $i = 1,\dots,n$.
    We will modify the inductive argument of Lemma \ref{series comp} to find a monomial with some $E_1(x_j)$.
    
    Let $1 < q < n$ be the number of distinct values among $a_{i,1,d}$.
    Let $(P_1,\dots,P_n)$ partition $\{1,\dots,n\}$ so that for any $i_1,i_2 \in P_l$, $\alpha_{l} :=a_{i_1,1,d} = a_{i_2,1,d}$ and $\alpha_{l_1} \ne \alpha_{l_2}$ for $l_1 \ne l_2$.
    Consider the sums
        $$\sum_{i \in P_l} \frac{E_1(x_j)^{\alpha_l}}{E^{(d)}(x_j)^{\alpha_l}}g_i$$
    for $l=1,\dots,q$.
    Each of these sums is an element of $T_{X}$ with strictly fewer terms than $t$, and the monomials in each sum have one less multiplicand, since $E^{(d)}(x_1)$ does not appear.
    
    We proceed by induction.
    As in Lemma \ref{series comp}, assume that for each $l=1,\dots,q$ we can find a leading term $b_lt_l$ of $\sigma_0\left(\sum_{i \in P_l} \frac{E'(x_j)^{\alpha_l}}{E^{(d)}(x_j)^{\alpha_l}}g_i\right)$, where $0 \ne b_l \in \boldsymbol{k}$ and $t_l \in G_{X}$.
    Assume further that each $t_l$ has some $E_1(x_j)$ as a generator.
    By the same reasoning as in Lemma \ref{series comp}, the sum 
    \begin{align*}
        \sum_{l=1}^q \sigma_0\left(\frac{E^{(d)}(x_1)^{\alpha_l}}{E'(x)^{\alpha_l}}\right)b_lt_l = \sum_{l=1}^q E_1(x_1)^{(d-1)\alpha_l} b_lt_l \sum_{k_d = 0}^{\alpha_l} \binom{\alpha_l}{k_d} \left(\dots + \frac{E_2(x_1)\cdots E_d(x_1)}{E_1(x_1)^{d-1}}\right)^{k_d}
    \end{align*}
    has a nonzero monomial $M$.
    Our additional hypothesis guarantees that $M$ has some $E_1(x_j)$ as a generator.
    Since $\mathrm{Lm}\big(\sigma_0(t)\big) \ge M$, $\mathrm{Lm}\big(\sigma_0(t)\big)$ must have $E_1(x_{j'})$ as a generator for some $j' \le j$.
\end{proof}

\begin{lemma}\label{not of the form E_0(x-1)}
    Suppose $X$ and $\boldsymbol{k}$ satisfy the ordering and separation assumptions in Remark \ref{order props}.
    If $t \in T_{X} \cap \boldsymbol{k}[G_{X}]$, then no monomial of $\rho_0(t)$ is of the form $E_0(x-1)^n$ for any $x \in X$, $n \in \mathbb{N}$.
\end{lemma}

\begin{proof}
    As in Lemma \ref{Leading monomial has E_1}, if $t = \sum_{i=1}^n c_ig_i \in T_{X} \cap \boldsymbol{k}[G_{X}]$, then for some $p,d \in \mathbb{N}$ we can write
    \begin{align*}
        g_i = \prod_{j=1}^p&\frac{E''(x_j)^{a_{i,j,2}} \cdots E^{(d)}(x_j)^{a_{i,j,d}}}{E'(x_j)^{a_{i,j,2} + \cdots + a_{i,j,d}}} \\
        \sigma_0(g_i) = \prod_{j=1}^p & E_0(x_j)^{a_{i,j,2} + \cdots + a_{i,j,d} - (a_{i,j,2} + \cdots + a_{i,j,d})} E_1(x_j)^{2a_{i,j,2} + \cdots + da_{i,j,d} - (a_{i,j,2} + \cdots + a_{i,j,d})}  \\
            &\left(\sum_{k_2=0}^{a_{i,j,2}}\binom{a_{i,j,2}}{k}\left(\frac{E_2(x_j)}{E_1(x_j)}\right)^{k_2}\right) \cdots \left(\sum_{k_{d}=0}^{a_{i,j,d}}\binom{a_{i,j,d}}{k_d}\left(\cdots + \frac{E_2(x_j) \cdots E_{d}(x_j)}{E_1(x_j)^{d-1}}\right)^{k_{d}}\right)
    \end{align*}
    As noted above, $\sigma_0(g_i)$ is a finite sum because $a_{i,j,2},\dots,a_{i,j,d}$ are natural numbers for each $j=1,\dots,p$, and the exponent of $E_0(x_j)$ is 0 in every term of $\sigma_0(g_i)$.
    Additionally, in every term of $\sigma_0(g_i)$ the sum of the exponents of generators composed with $x_j$ is $A_j  = a_{i,j,2} + 2a_{i,j,3} + \cdots + a_{i,j,d}$.
    
    Every monomial of $\rho_0(t) = \nu_0(\sigma_0(t))$ is an element of $\mathrm{Supp} (\nu_0(M))$ for some $i=1,\dots,n$ and some monomial $M$ of $\sigma_0(g_i)$.
    We can write any monomial $M$ of $\sigma_0(g_i)$ as
    \begin{align*}
        M = \prod_{j=1}^p E_1(x_j)^{N_{j,1}} \cdots E_d(x_j)^{N_{j,d}}
    \end{align*}
    with $N_{j,1},\dots,N_{j,d} \in \mathbb{N}$ and 
        $$\sum_{k=1}^d N_{j,k} = A_j$$
    for each $j=1,\dots,p$.
    Now 
    \begin{align*}
        \nu_0(M) &= \prod_{j=1}^p E_0(x-1)^{N_{j,1}}E_1(x-1)^{N_{j,1}+N_{j,2}}(1+\epsilon_2(x_j))^{N_{j,2}} \cdots E_{d-1}(x_j)^{N_{j,d}}(1 + \epsilon_d(x_j))^{N_{j,d}}.
    \end{align*}
    By Lemma \ref{exponent of E_d in epsilon_d}, we know any $\epsilon_l(x_j)$ for $l=2,\dots,d$, $j=1,\dots,p$ contributes 0 to the sum of the exponents and involves only $E_1(x_j-1),\dots,E_{l+1}(x_j-1)$.
    So the sum of the exponents of $E_1(x_j-1),\dots,E_{d}(x_j-1)$ is always  $\sum_{k=1}^d N_{j,k} = A_j$ for each $j=1,\dots,p$.
    Since at least one of these sums must be nonzero, $\nu_0(M)$ cannot have any generators of the form $E_0(x-1)^n$.
\end{proof}

\begin{cor}\label{form of rho_0(t)}
    If $0 \ne t \in T_{X}$ then $\mathrm{Lm}(\rho_0(t))$ must be of the form $E_0(x-1)^k \cdot h$ with $k >0$ and $1 \ne h \in H_{X-1}$.
\end{cor}

\begin{proof}
    Let $0 \ne t \in T_{X}$ and write $t = s_1\ell_1 + \cdots + s_n\ell_n$ with $s_1,\dots,s_n \in T_{X} \cap \kk[G_{X}]$.
    By Lemma \ref{Leading monomial has E_1}, for each $i=1,\dots,n$ the leading monomial of $\sigma_0(s_i)$ has a positive integer power $k_{i}$ of $E_1(x)$ for some $x \in X$.
    By the definition of $\Init(t)$ from Remark \ref{initial subsum of rho(s)}, we may ignore all monomials of $s_i$ except those in which $E_1(x)$ appears with exponent $k_i$.
    Recall that $\nu_0(E_1(x)^{k_i} = E_0(x-1)^{k_1}E_1(x-1)^{k_1}$.
    So for every monomial $M$ of $\Init(t)$, there is some $x \in X$ and integer $k >0$ such that $E_0(x-1)^k$ appears in $M$.
    
    By Lemma \ref{not of the form E_0(x-1)}, no monomial of $\rho_0(t)$ is of the form $E_0(x-1)^n$ for any $x \in X$, $n \in \mathbb{N}$.
    So every monomial of $\rho_0(t)$ must be of the form $E_0(x-1)^k \cdot h$ for some $x \in X$, $0 < k \in \NN$, and $1 \ne h \in H_{X-1}$.
\end{proof}

Now we apply Corollary \ref{form of rho_0(t)} to justify (1) and (3) of Remark \ref{why log hat is good enough}.
Let $0 \ne t \in T_{X}$, $g \in G_{X}$, and $\ell \in \Lambda_{X}$.
By Corollary \ref{form of rho_0(t)}, we can write 
    $$\Lm(\rho_0(t)) = E_0(x-1)^k \cdot h$$
for some $x \in X$, $0< k \in \NN$, and $1 \ne h \in H_{X-1}$.

For (1), if $g$ is small, then no generator $E_0(y)$ appears in $\Lm(\rho_0(g \cdot \ell))$.
So $\Lm(\rho_0(g \cdot \ell))$ is of the form
    $$\prod_{j=1}^p E_0(x_j-1)^{a_{j,0}}E_1(x_j-1)^{n_{j,1}} \cdots E_d(x_j-1)^{n_{j,d}}.$$
In particular, $\rho_0(g \cdot \ell) < E_0(x_1-1)^{a_{j,0}+1}$.
The difference equation for $E$ tells us that $\log E(x_1-1)^{a_{j,0}+1} = (a_{j,0} + 1)E(x_1-2)$, which is less than $\Lm(\rho_0(t))$.
So if $g$ is small, then we must have $e_{T_{X}}(t) > g \cdot \ell$.

For (3), if $g$ is not small, then $\Lm\left(\widehat{\log}(\rho_0(g \cdot \ell))\right)$ must be of the form $E_0(x-1) \ne \Lm(\rho_0(t))$.

We can now extend the order on $G_{X} \Lambda_{X}$ to an order on $\Gamma_{0}$:
\begin{defn}\label{monomial comp advanced}
    Let $0 \ne t \in T_{X}$, $g \in G_{X}$, and $\ell \in \Lambda_{X}$.
    Define $g \cdot \ell > e_{T_{X}}(t)$ in $\Gamma_{0}$ if and only if $\widehat{\log} (\rho_0(g \cdot \ell)) > \rho_0(t)$ in $\kk((H_{X-1}))$.
\end{defn}

Now that we have identified and ordered our stage 0 group of monomials, $\Gamma_0$, we would next like to define a modification of the Hahn series field with infinite sums of elements of $\Gamma_{0}$ but which disallows infinite sums of monomials that are too ``close".
First, we extend the notion of smallness from $G_{X}$ to $\Gamma_{0}$.

\begin{defn}
    Let $\Gamma_{0}^{small} := \{g \cdot \ell : g \in G_{X} \text{ is small, } \ell \in \Lambda_{X} \}$. 
\end{defn}

$\Gamma_{0}^{small}$ is a subgroup of $\Gamma_{0}$ because the small elements of $G_{X}$ form a subgroup.
By Lemma \ref{rho_0 is injective}, $\kk[\Gamma_{0}^{small}]$ is ordered.

\begin{cor}\label{convex subgroup}
    $\Gamma_{0}^{small}$ is a convex subgroup of $\Gamma_{0}$.
\end{cor}

\begin{proof}
    Let $g_0 \in G_{X}$ be small and $\ell_0 \in \Lambda_{X}$ with $g_0 \cdot \ell_0 > 1$.
    Let $1 < g_1 \cdot \ell_1 \cdot e_{T_{X}}(t_1) \in \Gamma_{0}$, and suppose $g_1 \cdot \ell_1 \cdot e_{T_{X}}(t_1) < g_0 \cdot \ell_0$.
    We will show $g_1 \cdot \ell_1 \cdot e_{T_{X}}(t_1) \in \Gamma_{0}^{small}$, i.e., $t_1 = 0$ and $g_1$ is small.
    
    If $t_1 \ne 0$, then following the way $\Gamma_{0}$ is ordered via Definition \ref{monomial comp advanced}, $1 < g_1 \cdot \ell_1 \cdot e_{T_{X}}(t_1) < g_0 \cdot \ell_0$ means
        $$0 < \rho_0(t_1) + \widehat{\log} \big(\rho_0(g_1 \cdot \ell_1)\big) < \widehat{\log} \big(\rho_0(g_0 \cdot \ell_0)\big).$$
    We can write $\mathrm{Lm}\big(\rho_0(g_0 \cdot \ell_0)\big)$ as
    \begin{align*}
        \prod_{j=1}^p E_0(x_j)^{\alpha_{j}} E_0(x_j-1)^{\beta_j + b_j}E_1(x_j-1)^{\beta_j + N_{j,1}}.
    \end{align*}
    Since $g_0$ is small, we must have $\alpha_{j} = 0$ for all $j=1,\dots,p$, so $\widehat{\log} \big(\rho_0(g_0 \cdot \ell_0)\big) = 0$, a contradiction.
    So we must have $t_1 = 0$.
    
    We can write $\mathrm{Lm}\big(\rho_0(g_1 \cdot \ell_1)\big)$ as
    \begin{align*}
        \prod_{j=1}^{p'} E_0(y_j)^{\alpha_{j}'} E_0(y_j-1)^{\beta_j' + b_j'}E_1(y_j-1)^{\beta_j' + N_{j,1}'}
    \end{align*}
    where $\prod_{j=1}^{p'} E_0(y_j-1)^{b_j'}$ is the contribution of $\ell_1$ to this monomial and the other multiplicands arise from $g_1$.
    Since $1 < g_1 \cdot \ell_1 < g_0 \cdot \ell_0$ and $\alpha_j = 0$ for all $j=1,\dots,p$, we must have $\alpha_j' = 0$ for all $j=1,\dots,p'$ too, i.e., $g_1$ is small. \qedhere
\end{proof}

Since $\Gamma_{0}^{small}$ is a convex subgroup of $\Gamma_{0}$, the quotient group $\Gamma_{0}/\Gamma_{0}^{small}$ is an ordered group.

\begin{defn}\label{modified Hahn series}
    Let $K$ be an ordered field and let $\Gamma$ be an ordered multiplicative abelian group.
    Suppose $H$ is a convex subgroup of $\Gamma$ such that $K[H]$ is ordered.
    Define $K((\Gamma))_H$ to be the ring whose elements are sums of the form
        $$s = \sum_{M \in \Gamma} c_MM$$
    with $c_M \in K$, where
    \begin{enumerate}
        \item for each coset $w \in \Gamma/H$, $\{M \in w : c_M \ne 0\}$ is finite, and 
        
        \item $\Supp(s) = \{M \in \Gamma : c_M \ne 0\}$ is reverse well-ordered in $\Gamma$.
    \end{enumerate}
    For $s \in K((\Gamma))_H$, define the leading coset $\Lv(s)$ to be the coset of $\Lm(s)$ in $\Gamma/H$.
    Define $s > 0$ in $K((\Gamma))_H$ if 
        $$\sum_{M \in \Lv(s)} c_M \frac{M}{\Lm(s)} > 0$$
    in $K[H]$.
    For $w \in \Gamma/H$, we will write $s \vert_w$ to denote the finite subsum of $s$ with monomials in $w$.
\end{defn}

Observe that $K((\Gamma))_{H}$ is not necessarily a field because a sum with more than one element in its leading coset may not have a multiplicative inverse.

We could form the ordered ring $\kk((\Gamma_{0}))_{\Gamma_{0}^{small}}$, whose order extends the order on $\kk[G_{X}\Lambda_{X}]$ from Lemma \ref{rho_0 is injective}.
However, we must impose one additional restriction on our modified Hahn series that will be necessary in Subsection \ref{subsection derivation}, when we define a derivation on the field of sublogarithmic-transexponential series. Even though we have not yet discussed how to define a derivation, we include the following example now, to motivate this restriction.

\begin{ex}\label{problem with derivation}
    Let $x$ be the germ of the identity function, and let 
        $$X = \left\{E(x) + r + \frac{1}{E(x^r)} : r \in (0,1) \cap \mathbb{R}\right\}.$$
    Let $\boldsymbol{k}$ be any ordered exponential field such that $\boldsymbol{k}$ and $X$ satisfy the ordering and separation assumptions of Remark \ref{order props} (for example, $\boldsymbol{k} = \mathbb{R}$ works).
    Then we would want to define
    \begin{align*}
        \frac{d}{dx}E\left(E(x) + r + \frac{1}{E(x^r)}\right) = E'\left(E(x) + r + \frac{1}{E(x^r)}\right)\left(E'(x) - \frac{rx^{r-1}E'(x^r)}{E(x^r)^2}\right).
    \end{align*}
    Now, $E'\left(E(x) + r + \frac{1}{E(x^r)}\right) \in \kk((\Gamma_{0}))_{\Gamma_{0}^{small}}$ for all $r \in (0,1)$.
    However, due to the finitary nature of the full sublogarithmic-transexponential series construction in Subections \ref{subsection 5.1} and \ref{subsection 5.2}, there is no structure in which $\{E(x^r) : r \in C\}$ can all simultaneously appear in a sum for any infinite $C \subset (0,1)$. 
    This is because $(x^r)_{r \in (0,1)}$ all have different growth rates.
    So we do not want to  allow sums with generators built from infinitely many elements of $X$ because the ``derivative" of such an element would not be in the structure.
\end{ex}

Because of the issue that arises in Example \ref{problem with derivation}, any element $s$ of our modified Hahn series construction must only involve generators built from some finite subset $X' \subset X$.
A similar issue may arise if the derivatives of exponents and coefficients in some sum $s$ do not all lie in one structure.
So instead of using all of $\kk((\Gamma_{0}))_{\Gamma_{0}^{small}}$ in our construction, we will introduce a framework that allows us to restrict to a subring of sums built up in a finitary way.

\begin{rmk}\label{Modified Hahn series well defined}
    Suppose $X$ and $\kk$ satisfy the ordering and separation assumptions of Remark \ref{order props}.
    If $X' \subset X$ and $\kk_{\alpha}$ is a subfield of $\kk$, then $\Gamma_{X',\kk_{\alpha},0}$ is a subgroup of $\Gamma_{X,\kk,0}$ and $\Gamma_{X',\kk_{\alpha},0}^{small}$ is a subgroup of $\Gamma_{X,\kk,0}^{small}$.
    So we can form the ordered subring 
        $$\kk_{\alpha}((\Gamma_{X',\kk_{\alpha},0}))_{\Gamma_{X',\kk_{\alpha},0}^{small}}$$
    of $\kk((\Gamma_{0}))_{\Gamma_{0}^{small}}$.
    Note that the sign of $s \in \kk_{\alpha}((\Gamma_{X',\kk_{\alpha},0}))_{\Gamma_{X',\kk_{\alpha},0}^{small}}$ matches the sign of $s$ viewed as an element of $\kk((\Gamma_{0}))_{\Gamma_{0}^{small}}$.
\end{rmk}

Let $[X]^{<\omega}$ be the collection of finite subsets of $X$ and $(\kk_{\alpha})_{\alpha \in \mathcal{A}}$ a directed set of subfields of $\kk$.
By Remark \ref{Modified Hahn series well defined}, the collection of ordered subrings of $\kk((\Gamma_{0}))_{\Gamma_{0}^{small}}$ of the form
\begin{align*}
    \kk_{\alpha}((\Gamma_{X',\kk_{\alpha},0}))_{\Gamma_{X',\kk_{\alpha},0}^{small}}
\end{align*}
for $X' \in [X]^{< \omega}$ and $\alpha \in \mathcal{A}$ forms a directed set.
So the direct limit of this collection of subrings is its union, which we will call
\begin{align*}
    K_{X,0} 
        &:= \bigcup_{(X', \alpha) \in [X]^{<\omega} \times \mathcal{A}} \kk_{\alpha}((\Gamma_{X',\kk_{\alpha},0}))_{\Gamma_{X',\kk_{\alpha},0}^{small}}
\end{align*}
Again by Remark \ref{Modified Hahn series well defined}, $K_{X,0}$ is ordered and its order extends the order on $\kk[G_{X} \Lambda_{X}]$ from Lemma \ref{rho_0 is injective}.
Define
\begin{align*}
    A_{X,0} 
        &:= \{s \in K_{X,0} : \mathrm{Supp} (s) > \Gamma_{0}^{small}\} \\
    B_{X,0} 
        &:= \{s \in K_{X,0} : \forall M \in \mathrm{Supp} (s), \exists M_0 \in \Gamma_{0}^{small} \text{ with } M \le M_0\} \\
    B_{X,0}^* 
        &:= T_{X} \oplus \boldsymbol{k} \oplus \{s \in K_{X,0} : \mathrm{Supp} (s) < \Gamma_{X,0}^{small}\}.
\end{align*}
Define $e_{X,0} : B_{X,0}^* \to (K_{X,0})^{>0}$ by
    $$e_{X,0}(t + r + \epsilon) = e_{T_{X}}(t) \exp(r) \sum_{n=0}^{\infty} \frac{\epsilon^n}{n!}$$
for $t \in T_{X}$, $r \in \boldsymbol{k}$, and $\epsilon \in B_{X,0}$ with $\mathrm{Supp} (\epsilon) < \Gamma_{0}^{small}$.
Notice that the only part of $\boldsymbol{k}[\Gamma_{0}^{small}]$ on which it makes sense to define an exponential function is $T_{X} \oplus \boldsymbol{k}$.
For example, there is no element of $K_{X,0}$ to represent the exponential of the ``large" infinitesimal $\frac{E(x)}{E'(x)}$.

Additionally, we introduce notation for the directed family of subfields from which $K_{X,0}$ is built. 
Given $(X', \alpha) \in [X]^{<\omega} \times \mathcal{A}$, define 
    $$(K_{X,0})_{(X', \alpha)} := \kk_{\alpha}((\Gamma_{X',\kk_{\alpha},0}))_{\Gamma_{X',\kk_{\alpha},0}^{small}}.$$ 
We will write $(A_{X,0})_{(X', \alpha)} := \{s \in (K_{X,0})_{(X', \alpha)} : \Supp(s) > \Gamma_{0}^{small}\}$.

We now adapt the notion of the \textit{first extension} of a pre-exponential ordered field from \cite{LEseries}: 
An \textit{almost pre-exponential ordered ring} $\big(K,(K_{\beta})_{\beta \in \mathcal{B}},A,B,B^*,e\big)$ consists of an ordered ring $K$, a directed family $(K_{\beta})_{\beta \in \mathcal{B}}$ of subfields of $K$ such that $K = \bigcup_{\beta \in \mathcal{B}} K_{\beta}$, an additive subgroup $A$ of $K$, a convex subgroup $B$ of $K$ with $K = A \oplus B$, a subgroup $B^*$ of $B$, and a strictly increasing homomorphism $e : B^* \to (K)^{>0}$. 

Define the \textit{first extension} $(K',A',B',(B^*)',e')$ of an almost pre-exponential ordered ring $(K,A,B,B^*,e)$:
\begin{enumerate}
    \item For each $\beta \in \mathcal{B}$, let $A_{\beta} = A \cap K_{\beta}$.
    Take a multiplicative copy $e(A_{\beta})$ of the ordered additive abelian group $A_{\beta}$ with order-preserving isomorphism $e_{A_{\beta}} : A_{\beta} \to e(A_{\beta})$.
    
    \item Define $K' = \bigcup_{\beta \in \mathcal{B}}\big(K_{\beta}((e_{A_{\beta}}(A_{\beta})))\big)$, a subring of the usual Mal'cev-Neumann series ring $K((e_A(A)))$.
    
    \item Let $A' = \{s \in K' : \mathrm{Supp} (s) > 1\}$, and $B' = \{s \in K' : \mathrm{Supp} (s) \le 1\}$, so that $K' = A' \oplus B'$ and $B'$ is a convex subring of $K'$. 
    
    \item Let $(B^*)' = A \oplus B^* \oplus \mathfrak{m}(B')$, and extend $e$ to $e' : (B^*)' \to (K')^{>0}$ as follows: 
    Let $a \in A$, $b \in B^*$, and $\epsilon \in \mathfrak{m}(B')$.
    Since $K = \bigcup_{\beta \in \mathcal{B}} K_{\beta}$, there are $\beta_a$, $\beta_b$, and $\beta_{\epsilon}$ such that $a \in K_{\beta_a}$, $e(b) \in K_{\beta_b}$, and $\epsilon \in K_{\beta_{\epsilon}}((e_{A_{\beta_{\epsilon}}}(A_{\beta_{\epsilon}})))$.
    Since $(K_{\beta})_{\beta \in \mathcal{B}}$ is a directed family, there is some $\gamma \in \mathcal{B}$ such that $K_{\beta_a}, K_{\beta_b}, K_{\beta_{\epsilon}} \subset K_{\gamma}$.
    Define 
        $$e(a + b + \epsilon) = e_{A_{\gamma}}(a)e(b) \sum_{n=0}^{\infty} \frac{\epsilon^n}{n!}$$
    which is an element of $K_{\gamma}((e_{A_{\gamma}}(A_{\gamma})))^{>0} \subset (K')^{>0}$.
    
    \item Associate to $K'$ the directed family $(K_{\beta}')_{\beta \in \mathcal{B}}$ of subfields of $K'$, where $K_{\beta}' = K_{\beta}((e_{A_{\beta}}(A_{\beta})))$.
\end{enumerate}
Then $\big(K',(K_{\beta}')_{\beta \in \mathcal{B}},A',B',(B^*)',e'\big)$ is an almost pre-exponential ordered ring, and $e'$ is defined on all of $A \oplus B^* \subset K$.

Starting with $\big(K_{X,0},((K_{X,0})_{\beta})_{\beta \in [X]^{<\omega} \times \mathcal{A}}, A_{X,0}, B_{X,0}, B_{X,0}^*,e_{X,0}\big)$, define 
    $$(K_{X,n+1}, ((K_{X,n+1})_{\beta})_{\beta \in [X]^{<\omega} \times \mathcal{A}}, A_{X,n+1}, B_{X,n+1}, B_{X,n+1}^*,e_{X,n+1})$$
to be the first extension of $(K_{X,n}, ((K_{X,n})_{\beta})_{\beta \in [X]^{<\omega} \times \mathcal{A}}, A_{X,n}, B_{X,n}, B_{X,n}^*,e_{X,n})$.
Then 
\begin{align*}
    B_{X,n+1}^* &= A_{X,n} \oplus \cdots \oplus A_{X,0} \oplus B_{X,0}^* \oplus \mathfrak{m}(B_{X,1}) \oplus \cdots \oplus \mathfrak{m}(B_{X,n+1}).
\end{align*}
Let $K_{X} = \bigcup_{n \in \mathbb{N}} K_{X,n}$, $B_{X}^* = \bigcup_{n \in \mathbb{N}} B_{X,n}^*$, and let $e_{X} : B_{X}^* \to (K_{X})^{>0}$ be the common extension of all the $e_{X,n}$.

We would also like to view elements of $K_{X}$ as sums of monomials with coefficients in $\boldsymbol{k}$ instead of some $K_{X,n}$.
To do this, we define an increasing sequence of multiplicative subgroups of $K_{X}$ starting with $\Gamma_0$ and taking 
    $$\Gamma_{n+1} := \Gamma_{n}e_X(A_{X,n}).$$
Since $\Gamma_{n}$ is convex in $\Gamma_{n+1}$, by induction $\Gamma_{0}^{small}$ is a convex subgroup of $\Gamma_{n}$ for all $n \in \mathbb{N}$.
Analogously for each $X \in [X]^{<\omega}$ and $\alpha \in \mathcal{A}$, starting with $\Gamma_{X',\boldsymbol{k}_{\alpha},0}$ we also define 
    $$\Gamma_{X',\boldsymbol{k}_{\alpha},n+1} := \Gamma_{X',\boldsymbol{k}_{\alpha},n}e_X\left(A_{X,n} \cap (K_{X,n})_{(X',\alpha)}\right).$$
Again, since $\Gamma_{X',\boldsymbol{k}_{\alpha},n}$ is convex in $\Gamma_{X',\boldsymbol{k}_{\alpha},n+1}$, by induction $\Gamma_{X',\boldsymbol{k}_{\alpha},0}^{small}$ is a convex subgroup of $\Gamma_{X',\boldsymbol{k}_{\alpha},n}$ for all $n \in \NN$.

\begin{lemma}\label{direct product of monomial groups}
    Let $\boldsymbol{k}$ be an ordered field, $\Gamma$ an ordered group with convex subgroup $G_1$ such that $\Gamma$ is the internal direct product of $G_1$ and another subgroup $G_2$.
    Let $H$ be a convex subgroup of $G_1$.
    Then $\boldsymbol{k}((G_1))_{H}((G_2)) \cong \boldsymbol{k}((\Gamma))_{H}$.
\end{lemma}
\begin{proof}
    The isomorphism is given by
    \begin{align*}
        \sum_{M \in G} c_MM \mapsto \sum_{M_2 \in G_2} &\left(\sum_{M_1 \in G_1} c_{M_1M_2}M_1\right)M_2. \qedhere
    \end{align*}
\end{proof}

We know $(K_{X,0})_{X',\alpha} =  \boldsymbol{k}_{\alpha}((\Gamma_{X', \boldsymbol{k}_{\alpha},0}))_{\Gamma_{X', \boldsymbol{k}_{\alpha},0}^{small}}$.
Assume $(K_{X,n})_{X',\alpha} = \boldsymbol{k}_{\alpha}((\Gamma_{X', \boldsymbol{k}_{\alpha},n}))_{\Gamma_{X', \boldsymbol{k}_{\alpha},0}^{small}}$.
Then using Lemma \ref{direct product of monomial groups},
\begin{align*}
    (K_{X,n+1})_{(X',\alpha)} 
        &= (K_{X,n})_{X',\alpha}\Big(\Big(e_{X}\big(A_{X,n} \cap (K_{X,n})_{(X',\alpha)}\big)\Big)\Big)
        \\
        &= \boldsymbol{k}_{\alpha}((\Gamma_{X', \boldsymbol{k}_{\alpha},n}))_{\Gamma_{X', \boldsymbol{k}_{\alpha},0}^{small}} \Big(\Big(e_{X}\big(A_{X,n} \cap (K_{X,n})_{(X',\alpha)}\big)\Big)\Big) \\
        &= \boldsymbol{k}_{\alpha}\Big(\Big(\Gamma_{X', \boldsymbol{k}_{\alpha},n} \cdot e_{X}\big(A_{X,n} \cap (K_{X,n})_{(X',\alpha)}\big)\Big)\Big)_{\Gamma_{X', \boldsymbol{k}_{\alpha},0}^{small}}
        \\
        &= \boldsymbol{k}_{\alpha}((\Gamma_{X',\boldsymbol{k}_{\alpha},n+1}))_{\Gamma_{X', \boldsymbol{k}_{\alpha},0}^{small}}.
\end{align*}
So by induction,
\begin{align*}
    K_{X,n+1} 
        &= \bigcup_{(X',\alpha) \in [X]^{<\omega} \times \mathcal{A}}(K_{X,n+1})_{(X',\alpha)}
        \\
        &= \bigcup_{(X',\alpha) \in [X]^{<\omega} \times \mathcal{A}} \boldsymbol{k}_{\alpha}((\Gamma_{X',\boldsymbol{k}_{\alpha},n+1}))_{\Gamma_{X', \boldsymbol{k}_{\alpha},0}^{small}}
\end{align*}
for all $n \in \NN$.

\subsection{Part 2: Building a logarithmic-exponential field}\label{Subsection 4.2}

Suppose $X$ and $\boldsymbol{k}$ satisfy the ordering and separation assumptions of Remark \ref{order props}.
Since $X$ is assumed to be a subset of an ordered field, the notation
    $$X-m := \{x-m : x \in X\}$$
has algebraic meaning.
Note that if $X$ and $\boldsymbol{k}$ satisfy the ordering and separation assumptions, then so do $X-m$ and $\boldsymbol{k}$.

In the previous subsection, we showed how to build a partial exponential ring $K_X$ starting with $E$-monomials built from $X$ and $\boldsymbol{k}$.
In this subsection, we would like to construct a logarithmic-exponential field, starting from $K_X$.
We motivate the construction in this subsection with several examples.
\begin{ex}\label{ex1}
    The difference equation for $E$ tells us that the logarithm of $E(x)$ should be $E(x-1)$.
    There is no element $E(x-1)$ in $K_X$, but if we build $K_{X-1}$ using $X-1$ instead of $X$, then $E(x-1)$ and $e_{X-1}\big(E(x-1)\big)$ are both elements of $K_{X-1}$, and we can identify them with $\log E(x)$ and $E(x)$ respectively.
\end{ex}

\begin{ex}\label{ex3}
    The natural way to represent the multiplicative inverse of $E'(x) + E(x)$ is by the following computation:
    \begin{align*}
        \frac{1}{E'(x) + E(x)} = \frac{1}{E'(x)\left(1 + \frac{E(x)}{E'(x)}\right)} = \frac{1}{E'(x)} \sum_{n=0}^{\infty}\left(-\frac{E(x)}{E'(x)}\right)^n.
    \end{align*}
    Since $E'(x)$ and $E(x)$ are in the same coset of $\Gamma_{X,0}^{small}$, the infinite sum $\sum_{n=0}^{\infty}\left(-\frac{E(x)}{E'(x)}\right)^n$ is not an element of $K_X$.
    However, if we identify $E(x)$ with $e_{X-1}\big(E(x-1)\big)$,  $E'(x)$ with $e_{X-1}\big(E(x-1)\big)E'(x-1)$, and $\frac{E(x)}{E'(x)}$ with $\frac{1}{E'(x-1)}$, then we have a multiplicative inverse 
        $$\frac{1}{e_{X-1}\big(E(x-1)\big)E'(x-1)} \sum_{n=0}^{\infty} \left(-\frac{1}{E'(x-1)}\right)^n \in K_{X-1}$$
    of $e_{X-1}\big(E(x-1)\big)E'(x-1) + e_{X-1}\big(E(x-1)\big)$.
\end{ex}

\begin{ex}\label{ex2}
    Reasoning as in Example \ref{ex3}, $\frac{E(x)}{E'(x)}$ does not have an exponential in $K_{X}$.
    However, if we identify $\frac{E(x)}{E'(x)}$ with $\frac{1}{E'(x-1)} \in K_{X-1}$, then
        $$e_{X-1}\left(\frac{1}{E'(x-1)}\right) = \sum_{n=0}^{\infty} \frac{1}{n!E'(x-1)^n} \in K_{X-1}.$$
\end{ex}

In this subsection, we will construct embeddings $\varphi_m : K_{X-m} \to K_{X-m-1}$ for $m \in \mathbb{N}$ which formally identify elements of $K_{X-m}$ with elements of $K_{X-m-1}$.
We will show that for all $s \in K_{X-m}$, the image of $s$ under finitely many embeddings has a multiplicative inverse, a logarithm, and an exponential in some $K_{X-m-j}$.

\subsubsection{Defining $\varphi_m$}
For each $n \in \NN$, we will define a map $\varphi_{m,n} : K_{X-m,n} \to K_{X-m-1,n+1}$.
Then we will take $\varphi_m$ to be the common extension of all $\varphi_{m,n}$ for $n \in \NN$.
The base case $\varphi_{m,0}$ requires more set-up than the inductive case $\varphi_{m,n+1}$, just as the construction of $K_{X-m,0}$ was more involved than the construction of $K_{X-m,n+1}$.

We begin by defining $\varphi_{m,0}$ on $\kk[G_{X-m}\Lambda_{X-m}]$ as follows:
\begin{enumerate}
    \item $\displaystyle \varphi_{m,0}\big(\log E'(x)^b\big) = E(x-1)^b \sum_{n = 0}^{\infty}\binom{b}{n}\left(\frac{\log E'(x-1)}{E(x-1)}\right)^n$ \\
    Since $\left(\frac{\log E'(x-1)}{E(x-1)}\right)^{n} < \Gamma_{X-m-1,0}^{small}$ for $n > 0$, this sum is allowed in $K_{X-m-1,0}$.
    
    \item $\displaystyle \varphi_{m,0}\big(E(x)^a\big) = e_{X-m-1}\big(E(x-1)\big)^a$
    
    \item $\displaystyle \varphi_{m,0}\big(E'(x)^a\big) = e_{X-m-1}\big(E(x-1)\big)^aE'(x-1)^a$
    
    \item $\displaystyle \varphi_{m,0}\big(E^{(d)}(x)^a\big) = e_{X-m-1}(E(x-1))^aE'(x-1)^{da} \sum_{n = 0}^{\infty} \binom{a}{n} \delta_d(x)^n$\\
    for $d > 1$, where $\delta_d(x) = \frac{B_d(x-1) - E'(x-1)^d}{E'(x-1)^d} = \frac{\binom{d}{2} E'(x-1)^{d-2}E''(x-1) + \cdots + E^{(d)}(x-1)}{E'(x-1)^d}$.\\
    Since $\mathrm{Supp} (\delta_d(x)^n) < \Gamma_{X-m-1,0}^{small}$ for $n >0$, this sum is allowed in $K_{X-m-1,0}$.
    
    \item Extend $\varphi_{m,0}$ to $\boldsymbol{k}[G_{X-m}\Lambda_{X-m}]$ so that it is a $\boldsymbol{k}$-algebra homomorphism. 
\end{enumerate}

In particular, we have now defined $\varphi_{m,0}$ on $T_{X-m}$.
Before we extend $\varphi_{m,0}$ further, we will show that $\varphi_{m,0}$ is order preserving on $\boldsymbol{k}[\Gamma_{X-m,0}^{small}]$.
Ultimately, this follows from the results of Section \ref{section 3} and Subsection \ref{subsection 4.1}, but the orderings on $\boldsymbol{k}[\Gamma_{X-m,0}^{small}]$ and $K_{X-m-1,0}$ are defined using the maps $\sigma_l$, $\nu_l$, and $\rho_l$ for $l=m,m+1$, which are rather delicate.
For this reason, it takes some maneuvering to use the results proven so far to show that $\varphi_{m,0}$ is order preserving.
We show that $\varphi_{m,0}$ is order preserving on $\boldsymbol{k}[\Gamma_{X-m,0}^{small}]$ in Lemma \ref{phi is order preserving small}.

\begin{lemma}\label{maps commute}
    Assume $\boldsymbol{k}$ and $X$ satisfy the ordering and separation assumptions in Remark \ref{order props}.
    If $t \in T_{X-m} \cap \boldsymbol{k}[G_{X-m}]$, then
        $$\sigma_{m+1}\big(\varphi_{m,0}(t)\big) = \nu_m\big(\sigma_m(t)).$$
\end{lemma}

\begin{proof}
    It suffices to show the result for expressions of the form $\frac{E^{(d)}(x)}{E'(x)} \in T_{X-m} \cap G_{X-m}$ because $\sigma_{m+1}$, $\varphi_{m,0}$, $\nu_m$, and $\sigma_m$ are $\boldsymbol{k}$-algebra homomorphisms.
    First, note that $\sigma_{m+1}$ is defined on $\varphi_{m,0}\left(\frac{E^{(d)}(x)}{E'(x)}\right)$ because
    \begin{align*}
        \varphi_{m,0}\left(\frac{E^{(d)}(x)}{E'(x)}\right) &= \frac{e_{X-m-1}\big(E(x-1)\big) B_d(x-1)}{e_{X-m-1}\big(E(x-1)\big) E'(x-1)} = \frac{B_d(x-1)}{E'(x-1)} \in \boldsymbol{k}[G_{X-m-1}].
    \end{align*}
    If $x$ is the germ of a function and we view all the following expressions as functions, we have
    \begin{align*}
        \frac{E^{(d)}(x)}{E'(x)} 
            &= \sigma_m\left(\frac{E^{(d)}(x)}{E'(x)}\right) 
            &&\text{ by definition of $E_0, E_1,\dots$} \\
            &= \nu_m\left(\sigma_m\left(\frac{E^{(d)}(x)}{E'(x)}\right)\right) 
            &&\text{ by Lemma \ref{log derivative shift}} \\
        \frac{E^{(d)}(x)}{E'(x)} 
            &= \varphi_{m,0}\left(\frac{E^{(d)}(x)}{E'(x)}\right) 
            &&\text{ by the Bell polynomial expression for $E^{(d)}$} \\
            &= \sigma_{m+1}\left(\varphi_{m,0}\left(\frac{E^{(d)}(x)}{E'(x)}\right)\right) 
            &&\text{ by definition of $E_0, E_1,\dots$}
    \end{align*}
    Since $\nu_m\left(\sigma_m\left(\frac{E^{(d)}(x)}{E'(x)}\right)\right)$ and $\sigma_{m+1}\left(\varphi_{m,0}\left(\frac{E^{(d)}(x)}{E'(x)}\right)\right)$ are expressed using the same standard basis $E_0(x-1), E_1(x-1),\dots, E_d(x-1)$ for monomials in $H_{X-m-1}$, they must also be equal in the sense of $\boldsymbol{k}((H_{X-m-1}))$.
    Since the computations of $\nu_m\left(\sigma_m\left(\frac{E^{(d)}(x)}{E'(x)}\right)\right)$ and $\sigma_{m+1}\left(\varphi_{m,0}\left(\frac{E^{(d)}(x)}{E'(x)}\right)\right)$ are the same regardless of whether $x$ is the germ of a function, these two expressions must be equal formally.
\end{proof}

\begin{rmk}\label{extension of maps commute}
    We can extend Lemma \ref{maps commute} beyond $T_{X-m} \cap \boldsymbol{k}[G_{X-m}]$ to include monomials of the form $E(x)^aE'(x)^{-a}$ for any $a \in \boldsymbol{k}$ and products of such monomials with elements of $T_{X-m} \cap \boldsymbol{k}[G_{X-m}]$, since
    \begin{align*}
        \varphi_{m,0}\left(E(x)^aE'(x)^{-a}\right) &= e_{X-m-1}\big(aE(x-1) - aE(x-1)\big)E'(x-1)^{-a} = E'(x-1)^{-a}.
    \end{align*}
    However, we cannot extend Lemma \ref{maps commute} any further because monomials in $\Gamma_{X-m,0}^{small}$ of other forms either involve logarithms, or their image under $\varphi_{m,0}$ is an infinite sum.
\end{rmk}

Define $\displaystyle \mathfrak{s}_{d,N,a}(x) := E(x)^a \frac{E'(x)^{da}}{E(x)^{da}}\sum_{k=0}^N\binom{a}{k} \left(\frac{E(x)^{d-1}E^{(d)}(x)}{E'(x)^d}-1\right)^k$ for $d \ge 2$.

\begin{lemma}\label{first N terms agree sigma}
    The first $N$ terms (at least) of $\sigma_m\big(E^{(d)}(x)^a\big)$ and $\sigma_m(\mathfrak{s}_{d,N,a}(x))$ are equal.
\end{lemma}

\begin{proof}
    Notice that $\mathfrak{s}_{d,N,a}(x)$ is obtained by paralleling the definition of $\sigma_m\big(E^{(d)}(x)^a\big)$:
    we have $\sigma_m(E(x)^a) = E_0(x)^a$, $\sigma_m\left(\frac{E'(x)^{da}}{E(x)^{da}}\right) = E_1(x)^{da}$, and
    \begin{align*}
        \sigma_m\left(\frac{E(x)^{d-1}E^{(d)}(x)}{E'(x)^d}-1\right)
            &= \frac{E_0(x)^d\big(E_1(x)^d + E_1(x)^{d-1}E_2(x) + \cdots + E_1(x)\cdots E_d(x)\big)}{E_0(x)^dE_1(x)^d} - 1
            \\
            &= \frac{E_1(x)^{d-1}E_2(x) + \cdots + E_1(x)\cdots E_d(x)}{E_1(x)^d}.
    \end{align*}
    Let 
        $$\zeta_d(x) = \frac{E_1(x)^{d-1}E_2(x) + \cdots + E_1(x)\cdots E_d(x)}{E_1(x)^d}.$$
    Then we need only check that the first $N$ terms of $\sum_{k=0}^{\infty}\binom{a}{k} \zeta_d(x)^k$ and $\sum_{k=0}^{N}\binom{a}{k} \zeta_d(x)^k$ agree.
    
    If $d=2$, then the sum of the first $N$ terms of $\sum_{k=0}^{\infty}\binom{a}{k} \zeta_d(x)^k$ is exactly $\sum_{k=0}^{N}\binom{a}{k} \zeta_d(x)^k$.
    So assume $d > 2$.
    
    The result holds immediately if $N=1$, so suppose $N >1$.
    Since $\frac{E_2(x)}{E_1(x)}$ is the largest monomial of $\zeta_d(x)$, the largest monomial of any $\zeta_d(x)^k$ is $\left(\frac{E_2(x)}{E_1(x)}\right)^k$, which has $E_1(x)^k$ as its denominator.
    The $N$th monomial $g$ of $\sum_{k=0}^{\infty}\binom{a}{k} \zeta(x)^k$ has $E_1(x)^l$ as its denominator for some $l <N$, so it cannot appear in $\zeta_d(x)^N$.
    Thus the first $N$ terms of $\sigma_m\big(E^{(d)}(x)\big)$ and $\sigma_m(\mathfrak{s}_{d,N,a}(x))$ must agree.
\end{proof}

\begin{cor}\label{first n terms overall}
    Assume $\boldsymbol{k}$ and $X$ satisfy the ordering and separation assumptions in Remark \ref{order props}.
    Let $s=c_1g_1 + \cdots + c_ng_n \in \boldsymbol{k}[G_{X-m}]$, and write
    \begin{align*}
        g_i = \prod_{j=1}^p E(x_j)^{a_{i,j,0}} \cdots E^{(d)}(x_j)^{a_{i,j,d}}
    \end{align*}
    with $x_1 > \cdots > x_p$.
    Suppose $\xi_{x_j}(g_i) = 0$ for all $j=1,\dots,p$.
    (Recall from Definition \ref{small defn} that $\xi_{x_j}(g_i) = a_{i,j,0} + \cdots + a_{i,j,d}$.)
    Then for every $N_0 \in \NN$ there exists some $N_1 \in \NN$ such that the first $N_0$ terms (at least) of $\sigma_m(s)$ and $\sigma_{m}(c_1u_1 + \cdots + c_nu_n)$ are equal, where 
    \begin{align*}
        u_i = \prod_{j=1}^pE(x_j)^{a_{i,j,0}} E'(x_j)^{a_{i,j,1}} \mathfrak{s}_{2,N_1,a_{i,j,2}}(x_j) \cdots \mathfrak{s}_{d,N_1,a_{i,j,d}}(x_j).
    \end{align*}
\end{cor}

\begin{proof}
    Let $N_0 \in \mathbb{N}$. 
    Let $t$ be the finite initial subsum of $\sigma_m(s)$ defined as in Remark \ref{initial subsum of rho(s)}.
    Enumerate $\mathrm{Supp} (\sigma_m(g_i))$ as $(h_{i,k})_{k \in \mathbb{N}}$ with $h_{i,k} > h_{ik+1}$.
    Let $n_i \in \mathbb{N}$ be smallest such that 
        $$\mathrm{Supp} (t) \cap \{h_{i,k} : k > n_i\} = \varnothing.$$
    Let $N_1 = \max(n_1,\dots,n_p)$.
    Let 
        $$u_i = \prod_{j=1}^pE(x_j)^{a_{i,j,0}} E'(x_j)^{a_{i,j,1}} \mathfrak{s}_{2,N_1,a_{i,j,2}}(x_j) \cdots \mathfrak{s}_{d,N_1,a_{i,j,d}}(x_j).$$
    
    So it suffices to show that the first $N_0$ terms of $\nu_m(\sigma_m(s))$ and $\nu_m(\sigma_m(c_1u_1 + \cdots + c_nu_n))$ agree.
    By Lemma \ref{first N terms agree sigma}, the first $N_1$ terms of $\sigma_m(u_i)$ and $\sigma_m(g_i)$ agree.
    By our choice of $N_1$, the first $N_0$ terms of $\sigma_m(c_1u_1 + \cdots + c_nu_n)$ and $\sigma_m(s)$ agree.
    So at least the first $N_0$ terms of $\nu_m(\sigma_m(c_1u_1 + \cdots + c_nu_n))$ and $\nu_m(\sigma_m(s))$ agree, and we are done.
\end{proof}

\begin{rmk}\label{full extension maps commute}
    Observe that
    \begin{enumerate}
        \item $\xi_x(g) = a$ for every monomial $g$ of $\mathfrak{s}_{d,N,a}$,
    
        \item $\varphi_{m,0}(\mathfrak{s}_{d,N,a})$ is a finite sum.
    \end{enumerate}
    Let $g_i$ and $u_i$ be as in Corollary \ref{first n terms overall}.
    By the first observation, no monomial of $\varphi_{m,0}(u_i)$ has any generator of the form $e_{X-m-1}(E(x-1))^a$.
    So $\mathrm{Supp} (\varphi_{m,0}(u_i)) \subset G_{X-m-1}$.
    This along with the second observation means we have $\varphi_{m,0}(u_i) \in \boldsymbol{k}[G_{X-m-1}]$.
    So $\sigma_{m+1}$ is defined on $\varphi_{m,0}(u_i)$.
    In fact, by Lemma \ref{maps commute} and Remark \ref{extension of maps commute}, 
        $$\sigma_{m+1}(\varphi_{m,0}(u_i)) = \nu_m(\sigma_m(u_i)).$$
    So $\sigma_{m+1}(\varphi_{m,0}(c_1u_1 + \cdots + c_nu_n)) = \nu_m(\sigma_m(c_1u_1 + \cdots + c_nu_n))$.
\end{rmk}

\begin{lemma}\label{phi is order preserving small}
    $\varphi_{m,0}$ is order preserving on $\boldsymbol{k}[\Gamma_{X-m,0}^{small}]$.
\end{lemma}

\begin{proof}
    Let $s = s_1\ell_1 + \cdots + s_n\ell_n$ with $s_i \in \boldsymbol{k}[G_{X-m}]$ and $\ell_i \in \Lambda_{X-m}$.
    We will show $s > 0$ if and only if $\varphi_{m,0}(s) > 0$.
    
    Recall from Lemma \ref{rho_0 is injective} that the sign of $\mathrm{Init}(s)$ determines the sign of $s$.
    We will show that $\varphi_{m,0}(s) \vert_{\mathrm{Lv}(\varphi_{m,0}(s))} \in \boldsymbol{k}[G_{X-m-1}]$, which means the sign of $\varphi_{m,0}(s)$ is determined by the sign of $\sigma_{m+1}\left(\varphi_{m,0}(s) \vert_{\mathrm{Lv}(\varphi_{m,0}(s))}\right)$.
    We will then show that $\sigma_{m+1}\left(\varphi_{m,0}(s) \vert_{\mathrm{Lv}(\varphi_{m,0}(s))}\right) = \mathrm{Init}(s)$.
    
    Let $t_i$ be as in Lemma \ref{rho_0 is injective}, and let $N_{i,0} = |\mathrm{Supp} (t_i)|$.
    Let $N_{0} = \max(N_{1,0},\dots,N_{p,0})$.
    For each $i=1,\dots,p$, let $N_{i,1}$ and $c_{i,1}u_{i,1} + \cdots + c_{i,k_i}u_{i,k_i}$ be given by Lemma \ref{first n terms overall}, using $N_0$.
    Define $u_i := c_{i,1}u_{i,1} + \cdots + c_{i,k_i}u_{i,k_i}$.
    
    By Remark \ref{full extension maps commute},
        $$\sigma_{m+1}\big(\varphi_{m,0}(u_i)\big) = \nu_m\big(\sigma_m(u_i)\big) = \rho_0(u_i).$$
    If we write $\ell_i = \prod_{j=1}^p \log E'(x_j)^{b_{i,j}}$, then $\mathrm{Lm}\big(\varphi_{m,0}(\ell_i)\big) = \prod_{j=1}^p E(x_j-1)^{b_{i,j}}$, so 
    \begin{align*}
        \sigma_{m+1}\big(\mathrm{Lm}(\varphi_{m,0}(\ell_i))\big) = \prod_{j=1}^p E_0(x_j-1)^{b_{i,j}}.
    \end{align*}
    Following the proof of Lemma \ref{rho_0 is injective}, $\mathrm{Init}(s) \ne 0$ is the initial subsum of
    \begin{multline*}
        \rho_m(u_1) \prod_{j=1}^p E(x_j-1)^{b_{1,j}} + \cdots + \rho_m(u_n)\prod_{j=1}^p E(x_j-1)^{b_{n,j}} = \\
            \sigma_{m+1}\big(\varphi_{m,0}(u_1)\mathrm{Lm}(\varphi_{m,0}(\ell_1)) + \cdots + \varphi_{m,0}(u_n)\mathrm{Lm}(\varphi_{m,0}(\ell_n))\big) 
    \end{multline*}
    so both these expressions have the same sign as $\mathrm{Init}(s)$.
    Recall from Remark \ref{initial subsum of rho(s)} that $\mathrm{Init}(s)$ is defined by having tuples of exponents of $E_0(x_1-1),\dots,E_0(x_p-1)$ maximized. 
    By the way $\sigma_{m+1}$ is defined, if $g \in \boldsymbol{k}[G_{X-m-1}]$ and $\xi_{x_j}(g) = a$, then the exponent of $E_0(x_j-1)$ in $\sigma_{m+1}(g)$ is $a$.
    So $\mathrm{Init}(s)$ must be the image under $\sigma_{m+1}$ of 
        $$\big(\varphi_{m,0}(u_1)\mathrm{Lm}(\varphi_{m,0}(\ell_1)) + \cdots + \varphi_{m,0}(u_n)\mathrm{Lm}(\varphi_{m,0}(\ell_n)\big)|_w$$ 
    where
        $$w = \mathrm{Lv}\big(\varphi_{m,0}(u_1)\mathrm{Lm}(\varphi_{m,0}(\ell_1)) + \cdots + \varphi_{m,0}(u_n)\mathrm{Lm}(\varphi_{m,0}(\ell_n)\big).$$
    Since for any 
        $$g \in \mathrm{Supp}\Big(\varphi_{m,0}(u_i)\big(\varphi_{m,0}(\ell_i) - \mathrm{Lm}(\varphi_{m,0}(\ell_i))\big)\Big)$$
    $\sigma_{m+1}(g) < \mathrm{Init}(s)$, we must have $\mathrm{Lv}\big(\varphi_{m,0}(s)\big) = w$.
    This completes the proof.
\end{proof}

We would now like to extend $\varphi_{m,0}$ to all of $\Gamma_{X-m,0}$ by defining 
    $$\varphi_{m,0}\big(e_{X-m}(t)\big) = e_{X-m-1}\big(\varphi_{m,0}(t)\big)$$
for $t \in T_{X-m}$, but we must check that $e_{X-m-1}$ is defined on $\varphi_{m,0}(t)$.

\begin{lemma}\label{exp of phi of t}
    Assume $\boldsymbol{k}$ and $X$ satisfy the ordering and separation assumptions in Remark \ref{order props}.
    If $t \in T_{X-m}$, then $e_{X-m-1}\big(\varphi_{m,0}(t)\big)$ is defined in $K_{X-m-1,1}$.
\end{lemma}

\begin{proof}
    Since any element of $T_{X-m}$ is a finite sum of generators, it suffices to show that $e_{X-m-1}$ is defined on the image of generators of $T_{X-m}$ under $\varphi_{m,0}$.
    Let 
        $$t = c\prod_{j=1}^p \log E'(x_j)^{a_{j,-1}} \frac{E''(x_j)^{a_{j,2}} \cdots E^{(d)}(x_j)^{a_{j,d}}}{E'(x_j)^{a_{j,2} + \cdots + a_{j,d}}}$$
    be a generator of $T_{X-m}$, with $x_1 > \cdots > x_p \in X - m$ and $a_{j,-1},a_{j,2},\dots,a_{j,d} \in \mathbb{N}$.
    Notice that for $d \ge 1$,
        $$\varphi_{m,0}(E^{(d)}(x)) = e_{X-m-1}(E(x-1))B_d(x-1)$$
    which matches the difference-differential equation for $E^{(d)}(x)$.
    For each term $g$ of $B_d(x-1)$, $\xi_{x-1}(g) \ge 1$, with $\xi_{x-1}(g) = 1$ only for the smallest term $g = E^{(d)}(x-1)$.
    So for each term $g$ of $\delta_d(x) = \frac{B_d(x-1) - E'(x-1)^d}{E'(x-1)^d}$, $-(d-1) \le \xi_{x-1}(g) \le -1$, with $\xi_{x-1}(g) = -(d-1)$ only for the smallest term $g = \frac{E^{(d)}(x-1)}{E'(x-1)^d}$.
    
    By the definition of $\varphi_{m,0}$, we can write
    \begin{align*}
        \varphi_{m,0} \Bigg(c\prod_{j=1}^p &\log E'(x_j)^{a_{j,-1}} \frac{E''(x_j)^{a_{j,2}} \cdots E^{(d)}(x_j)^{a_{j,d}}}{E'(x_j)^{a_{j,2} + \cdots + a_{j,d}}}\Bigg) = \\
        c\prod_{j=1}^p &E(x_j-1)^{a_{j,-1}} E'(x_j-1)^{a_{j,2} + 2a_{j,3} + \cdots + (d-1)a_{j,d}}  \\
            &\left(\sum_{n = 0}^{a_{j,-1}}\binom{a_{j,-1}}{n}\left(\frac{\log E'(x_j-1)}{E(x_j-1)}\right)^n \right) \left(\sum_{n = 0}^{a_{j,2}} \binom{a_{j,2}}{n} \delta_2(x_j)^n\right) \cdots \left(\sum_{n = 0}^{a_{j,d}} \binom{a_{j,d}}{n} \delta_d(x_j)^n\right).
    \end{align*}
    This is a finite sum.
    For every term $g$ in the sum, $\xi_{x_j-1}(g) \ge 0$ for all $j=1,\dots,p$.
    For every term $g$ except the smallest term, $\xi_{x_1-1}(g) + \cdots + \xi_{x_p-1}(g) > 0$.
    Therefore, every term except the smallest term is in $A_{X-m-1,0}$.
    The smallest term is
    \begin{align*}
        c\prod_{j=1}^p &E(x_j-1)^{a_{j,-1}} E'(x_j-1)^{a_{j,2} + 2a_{j,3} + \cdots + (d-1)a_{j,d}} \\
            &\binom{a_{j,-1}}{a_{j,-1}} \left(\frac{\log E'(x_j-1)}{E(x_j-1)}\right)^{a_{j,-1}} 
            \binom{a_{j,2}}{a_{j,2}} \left(\frac{E''(x_j-1)}{E'(x_j-1)^2}\right)^{a_{j,2}} \cdots
            \binom{a_{j,d}}{a_{j,d}} \left(\frac{E^{(d)}(x_j-1)}{E'(x_j-1)^d}\right)^{a_{j,d}} = \\
        c\prod_{j=1}^p &\log E'(x_j-1)^{a_{j,-1}} \frac{E''(x_j-1)^{a_{j,2}} \cdots E^{(d)}(x_j-1)^{a_{j,d}}}{E'(x_j-1)^{a_{j,2} + \cdots + a_{j,d}}}
    \end{align*}
    and is an element of $T_{X-m-1}$.
    Thus $e_{X-m-1}(t) \in K_{X-m-1,1}$.
\end{proof}

\begin{rmk}\label{phi_m transfers t}
    In the proof of Lemma \ref{exp of phi of t}, observe that the smallest term of $\varphi_{m,0}(t)$ is exactly $t$ with each $x_j$ replaced by $x_j-1$.
\end{rmk}

We can now extend $\varphi_{m,0}$ to all of $\Gamma_{X-m,0}$.
\begin{enumerate}
    \setcounter{enumi}{5}
    \item $\displaystyle \varphi_{m,0}\big(e_{X-m}(t)\big) = e_{X-m-1}(\varphi_{m,0}(t))$
    
    \item Extend $\varphi_{m,0}$ to $\boldsymbol{k}[\Gamma_{X-m,0}]$ so that it is a $\boldsymbol{k}$-algebra homomorphism.
\end{enumerate}

Finally, we would like to define 
    $$\varphi_{m,0} \left(\sum_{M \in \Gamma_{X-m,0}}c_MM\right) = \sum_{M \in \Gamma_{X-m,0}} c_M \varphi_{m,0}(M)$$
but we must check that $\sum_{M \in \Gamma_{X-m,0}} c_M \varphi_{m,0}(M)$ is a valid sum in $K_{X-m-1,1}$.
We will check this using the next two lemmas.

\begin{lemma}\label{M_1 > M_2 in different cosets}
    Let $M_i = g_i \cdot \ell_i \cdot e_{X-m}(t_i) \in \Gamma_{X-m,0}$ for $i=1,2$.
    If $M_1$ and $M_2$ are in different cosets of $\Gamma_{X-m,0}/\Gamma_{X-m,0}^{small}$ and $M_1 > M_2$, then $\mathrm{Supp} (\varphi_{m,0}(M_1)) > \mathrm{Supp} (\varphi_{m,0}(M_2))$.
\end{lemma}

\begin{proof}
    Let $x_1 > \cdots > x_p$ list the elements of $X$ that appear in $M_1, M_2$.
    For $i=1,2$, we can write
    \begin{align*}
        \varphi_{m,0}(M_i) &= e_{X-m-1}\left(\varphi_{m,0}(t_i) + \sum_{j=1}^p \xi_{x_j}(M_i)E(x_j-1)\right)s_i
    \end{align*}
    for some $s_i \in K_{X-m-1,0}$.
    So it suffices to show that
        $$\varphi_{m,0}(t_1) + \sum_{j=1}^p \xi_{x_j}(M_1)E(x_j-1) > \varphi_{m,0}(t_2) + \sum_{j=1}^p \xi_{x_j}(M_2)E(x_j-1).$$
    
    $M_1 > M_2$ means 
    \begin{align*}
        \rho_m(t_1) + \sum_{j=1}^p \xi_{x_j}(M_1)E_0(x_j-1) &> \rho_{m}(t_2) + \sum_{j=1}^p \xi_{x_j}(M_2)E_0(x_j-1) \\
        \text{i.e., } \hspace{.3cm} \rho_m(t_1 - t_2) &> \sum_{j=1}^p \big(\xi_{x_j}(M_2) - \xi_{x_j}(M_1)\big)E_0(x_j-1).
    \end{align*}
    Write $t_1-t_2 = h_{1}\lambda_{1} + \cdots + h_{n}\lambda_{n}$ as a sum with $h_{i} \in G_{X-m}$ and $\lambda_{i} \in \Lambda_{X-m}$.
    Note that the exponents in all $h_{i}$ and $\lambda_{i}$ are integers since $t_{i} \in T_{X-m}$.
    Let $\mathrm{Init}(t_1-t_2)$ be the initial subsum of $\rho_m(t_1-t_2)$ defined in Remark \ref{initial subsum of rho(s)}.
    We know $\nu_m(\sigma_m(h_{i})) = \sigma_{m+1}(\varphi_{m,0}(h_{i}))$ by Lemma \ref{maps commute}.
    We also know the only term of $\rho_m(\lambda_i)$ that can contribute to $\mathrm{Init}(t_1-t_2)$ is $\mathrm{Lm}\big(\rho_m(\lambda_i)\big)$ for all $i=1,\dots,n$.
    Let $\mathrm{Lm}\big(\rho_m(\lambda_i)\big) = \prod_{j=1}^p E_0(x_j-1)^{b_{i,j}}$.
    So we must have 
    \begin{align*}
        \sum_{i=1}^n \sigma_{m+1}(\varphi_{m,0}(h_{i})) \prod_{j=1}^p E_0(x_j-1)^{b_{i,j}} > \sum_{j=1}^p \big(\xi_{x_j}(M_2) - \xi_{x_j}(M_1)\big)E_0(x_j-1)
    \end{align*}
    in $\kk((H_{X-m-1}))$.
    Since $\sigma_{m+1}$ is order preserving, we must have
    \begin{align*}
        \sum_{i=1}^n\varphi_{m,0}(h_i) \prod_{j=1}^p E(x_j-1)^{b_{i,j}} > \sum_{j=1}^p \big(\xi_{x_j}(M_2) - \xi_{x_j}(M_1)\big)E(x_j-1).
    \end{align*}
    in $\kk[G_{X-m-1}\Lambda_{X-m-1}]$.
    Just as in Lemma \ref{phi is order preserving small}, $\varphi_{m,0}(t_1 - t_2)|_{\mathrm{Lv}(\varphi_{m,0}(t_1 - t_2))}$ is a subsum of $\sum_{i=1}^n\varphi_{m,0}(h_i) \prod_{j=1}^p E(x_j-1)^{b_{i,j}}$.
    So 
        $$\varphi_{m,0}(t_1 - t_2) > \sum_{j=1}^p \big(\xi_{x_j}(M_2) - \xi_{x_j}(M_1)\big)E_0(x_j-1)$$
    as desired.
\end{proof}

\begin{lemma}
    Assume $\boldsymbol{k}$ and $X$ satisfy the ordering and separation assumptions in Remark \ref{order props}.
    If $s \in K_{X-m}$, then $\varphi_{m,0}(s)$ is an element of $K_{X-m-1,1}$.
\end{lemma}

\begin{proof}
    Let $s = \sum_{M \in \Gamma_{X-m,0}} c_MM$. 
    We will show $\varphi_{m,0}\left(\sum_{M \in \Gamma_{X-m,0}} c_MM\right)$ is a valid sum in $K_{X-m-1,1}$.
    Since each $M \in \Gamma_{X-m,0}$ is a finite product of generators, $\varphi_{m,0}$ is well defined on monomials.
    So it suffices to check that $(c_M\varphi_{m,0}(M) : M \in \Gamma_{X-m,0})$ is summable in the sense of Definition \ref{modified Hahn series}, i.e.,
    \begin{enumerate}
        
        \item For each coset $N \Gamma_{X-m-1,0}^{small}$ of $\Gamma_{X-m-1,1}/\Gamma_{X-m-1,0}^{small}$, there are only finitely many $M \in \Gamma_{X-m,0}$ such that $c_M \ne 0$ and $N \Gamma_{X-m-1,0}^{small} \cap \mathrm{Supp} (\varphi_{m,0}(M)) \ne \varnothing$.
        
        \item $\displaystyle \bigcup_{M \in \mathrm{Supp} (s)} \mathrm{Supp} (\varphi_{m,0}(M))$ is reverse well-ordered in $\Gamma_{X-m-1,1}$.
    \end{enumerate}
    
    We start with (1).
    Let $N = g \cdot \ell \cdot e_{X-m-1}(t)\cdot e_{X-m-1}(\alpha) \in \Gamma_{X-m-1,1}$ with $g \in G_{X-m-1}$, $\ell \in \Lambda_{X-m-1}$, $t \in T_{X-m-1}$, and $\alpha \in A_{X-m-1,0}$. 
    Note that $t$ and $\alpha$ are fixed across all elements of $N \Gamma_{X-m-1,0}^{small}$. 
    Let $M_i = g_i \cdot \ell_i \cdot e_{X-m}(t_i) \in \Gamma_{X-m,0}$
    for $i=1,2$, and let $x_1 > \cdots > x_p$ list the elements of $X$ that appear in $M_1,M_2$.
    Suppose $N \Gamma_{X-m-1,0}^{small} \cap \mathrm{Supp} (\varphi_{m,0}(M_i)) \ne \varnothing$ for $i=1,2$.
    
    We can write
    \begin{align*}
        \varphi_{m,0}(M_i) &= e_{X-m-1}\left(\varphi_{m,0}(t_i) + \sum_{j=1}^p \xi_{x_j}(M_i)E(x_j-1)\right)s_i
    \end{align*}
    for some $s_i \in K_{X-m-1,0}$.
    $M_1$ and $M_2$ must be in the same coset of $\Gamma_{X-m,0}/\Gamma_{X-m,0}^{small}$, or else by Lemma \ref{M_1 > M_2 in different cosets} we would have
        $$\varphi_{m,0}(t_1) + \sum_{j=1}^p \xi_{x_j}(M_1)E(x_j-1) \ne \varphi_{m,0}(t_2) + \sum_{j=1}^p \xi_{x_j}(M_2)E(x_j-1)$$
    which contradicts that $N \Gamma_{X-m-1,0}^{small} \cap \mathrm{Supp} (\varphi_{m,0}(M_i)) \ne \varnothing$ for both $i=1,2$.
    Since there are only finitely many $M$ in the same coset as $M_1,M_2$ such that $c_{M} \ne 0$, we have proven (1).
    
    Now we prove (2).
    Let $B \subset \displaystyle \bigcup_{M \in \mathrm{Supp} (s)} \mathrm{Supp} (\varphi_{m,0}(M))$.
    We will find a largest element of $B$.
    Since $s|_w$ is finite for any coset $w \in \Gamma_{X-m,0}/\Gamma_{X-m,0}^{small}$, $B \cap \mathrm{Supp}(\varphi_{m,0}\left(s|_w\right))$ has a largest element.
    Now suppose $M_1$ and $M_2$ are two elements of $\Gamma_{X-m,0}$ in different cosets of $\Gamma_{X-m,0}/\Gamma_{X-m,0}^{small}$ with $M_1 > M_2$.
    By Lemma \ref{M_1 > M_2 in different cosets}, we must have $\mathrm{Supp} (\varphi_{m,0}(M_1)) > \mathrm{Supp} (\varphi_{m,0}(M_2))$.
    So the largest element of $B$ is the largest element of 
        $$B \cap \mathrm{Supp} \left(\varphi_{m,0}\left(s|_{\mathrm{Lv}(s)}\right)\right).$$
    So $\displaystyle \bigcup_{M \in \mathrm{Supp} (s)} \mathrm{Supp} (\varphi_{m,0}(M))$ is reverse well-ordered in $\Gamma_{X-m-1,1}$.
\end{proof}

We repeat the full definition of $\varphi_{m,0} : K_{X-m,0} \to K_{X-m-1,1}$:
\begin{enumerate}
    \item $\displaystyle \varphi_{m,0}\big(\log E'(x)^b\big) = E(x-1)^b \sum_{n = 0}^{\infty}\binom{b}{n}\left(\frac{\log E'(x-1)}{E(x-1)}\right)^n$ \\
    Since $\left(\frac{\log E'(x-1)}{E(x-1)}\right)^{n} < \Gamma_{X-m-1,0}^{small}$ for $n > 0$, this sum is allowed in $K_{X-m-1,0}$.
    
    \item $\displaystyle \varphi_{m,0}\big(E(x)^a\big) = e_{X-m-1}\big(E(x-1)\big)^a$
    
    \item $\displaystyle \varphi_{m,0}\big(E'(x)^a\big) = e_{X-m-1}\big(E(x-1)\big)^aE'(x-1)^a$
    
    \item $\displaystyle \varphi_{m,0}\big(E^{(d)}(x)^a\big) = e_{X-m-1}(E(x-1))^aE'(x-1)^{da} \sum_{n = 0}^{\infty} \binom{a}{n} \delta_d(x)^n$\\
    for $d > 1$, where $\delta_d(x) = \frac{B_d(x-1) - E'(x-1)^d}{E'(x-1)^d} = \frac{\binom{d}{2} E'(x-1)^{d-2}E''(x-1) + \cdots + E^{(d)}(x-1)}{E'(x-1)^d}$\\
    Since $\mathrm{Supp} (\delta_d(x)^n) < \Gamma_{X-m-1,0}^{small}$ for $n>0$, this sum is allowed in $K_{X-m-1,0}$.
    
    \item Extend $\varphi_{m,0}$ to $\boldsymbol{k}[G_{X-m}\Lambda_{X-m}]$ so that it is a $\boldsymbol{k}$-algebra homomorphism. 
    
    \item $\displaystyle \varphi_{m,0}\big(e_{X-m}(t)\big) = e_{X-m-1}(\varphi_{m,0}(t))$
    
    \item $\displaystyle \varphi_{m,0}(g_1 \cdots g_n) = \varphi_{m,0}(g_1) \cdots \varphi_{m,0}(g_n)$ for generators $g_1,\dots,g_n \in \Gamma_{X-m,0}$
    
    \item $\displaystyle \varphi_{m,0} \left(\sum_{M \in \Gamma_{X-m,0}}c_MM\right) = \sum_{M \in \Gamma_{X-m,0}} c_M \varphi_{m,0}(M)$. 
\end{enumerate}

\begin{cor}
    $\varphi_{m,0}$ is order preserving on $K_{X-m,0}$.
\end{cor}

\begin{proof}
    Let $s \in K_{X-m,0}$.
    By Lemma \ref{M_1 > M_2 in different cosets}, if $M_1 \in \mathrm{Supp} \left(s|_{\mathrm{Lv}(s)}\right)$ and $M_2 \not \in \mathrm{Supp} \left(s|_{\mathrm{Lv}(s)}\right)$, then $\mathrm{Supp} (\varphi_{m,0}(M_1)) > \mathrm{Supp} (\varphi_{m,0}(M_1))$.
    So the sign of $\varphi_{m,0}(s)$ is determined by the sign of $\varphi_{m,0}\left(s|_{\mathrm{Lv}(s)}\right)$.
    The sign of $\varphi_{m,0}\left(s|_{\mathrm{Lv}(s)}\right)$ is the same as the sign of $\varphi_{m,0}\left(\frac{s|_{\mathrm{Lv}(s)}}{\mathrm{Lm}(s)}\right)$, and $\frac{s|_{\mathrm{Lv}(s)}}{\mathrm{Lm}(s)} \in \boldsymbol{k}[\Gamma_{X-m,0}^{small}]$.
    Apply Lemma \ref{phi is order preserving small} to finish the proof.
\end{proof}

Now given $\varphi_{m,n} : K_{X-m,n} \to K_{X-m-1,n+1}$, define $\varphi_{m,n+1} : K_{X-m,n+1} \to K_{X-m-1,n+2}$ as follows:
    $$\displaystyle \varphi_{m,n+1} \left(\sum f_ae_{X-m}(a)\right) = \sum \varphi_{m,n}(f_a) e_{X-m}(\varphi_{m,n}(a))$$
for $f_a \in K_{X-m,n}$ and $a \in A_{X-m,n}$.
This is a valid sum in $K_{X-m-1,n+2}$, and since $\varphi_{m,n}$ is order preserving, so is $\varphi_{m,n+1}$.
Let $\varphi_m : K_{X-m} \to K_{X-m-1}$ be the common extension of all $\varphi_{m,n}$, which is also order preserving.

\subsection{Finding logarithms, exponentials, and multiplicative inverses}
In the remainder of this section, we will show that for every $s \in K_{X-m}$, there is some $l \in \mathbb{N}$ such that
\begin{enumerate}
    \item $e_{X-m-l-1}\big((\varphi_{m+l} \circ \cdots \circ \varphi_m)(s)\big)$ is defined in $K_{X-m-l-1}$ (Lemma \ref{exp closure}), 
    
    \item if $s>0$ then there is some $s' \in K_{X-m-l-1}$ such that 
        $$(\varphi_{m+l} \circ \cdots \circ \varphi_m)(s) = e_{X-m-l-1}(s')$$
    i.e., this element has a logarithm (Lemma \ref{log closure}), and
    
    \item $(\varphi_{m+l} \circ \cdots \circ \varphi_m)(s)$ has a multiplicative inverse in $K_{X-m-l-1}$ (Lemma \ref{field}).
\end{enumerate}
We find multiplicative inverses with the more general result that if $s_1,\dots,s_n \in K_{X-m}$ with $\mathrm{Supp} (s_i) < 1$, then there is $l \in \mathbb{N}$ such that
    $$\sum_{j_1,\dots,j_n = 0}^{\infty} a_{j_1,\dots,j_n}\prod_{i=1}^n(\varphi_{m+l} \circ \cdots \circ \varphi_m)(s_i)^{j_i} \in K_{X-m-l-1}$$
for all $(a_{j_1,\dots,j_n})_{\overline{j} \in \mathbb{N}^n} \in \boldsymbol{k}^{\omega}$ (Lemma \ref{R_an closure}).
This result will be necessary in order to show our construction is closed under restricted analytic functions.

\begin{lemma}\label{phi transfers exp}
    If $s = e_{X-m}(s')$ for some $s' \in B_{X-m}^*$, then $\varphi_m(s) = e_{X-m-1}\big(\varphi_m(s')\big)$.
\end{lemma}

\begin{proof}
    Since $s' \in B_{X-m,n}^* $ for some $n \in \mathbb{N}$, we can write
        $$s' = a + t + r + \epsilon$$
    with $a \in A_{X-m,n-1} \oplus \cdots \oplus A_{X-m,0}$, $t \in T_{X-m}$, $r \in \boldsymbol{k}$, and $\epsilon \in \mathfrak{m}(B_{X-m,1}) \oplus \cdots \oplus \mathfrak{m}(B_{X-m,n})$.
    So
    \begin{align*}
        \varphi_m(s) &= \varphi_m\left(e_{X-m}(a) e_{X-m}(t) e_{X-m}(r) \cdot \sum_{k=0}^{\infty}\frac{\epsilon^{k}}{k!} \right) \\
            &= e_{X-m-1}\big(\varphi_m(a)\big) e_{X-m-1}\big(\varphi_m(t)\big) e_{X-m-1}\big(\varphi_m(r)\big) \sum_{k=0}^{\infty}\frac{\varphi_m(\epsilon)^{k}}{k!}  \\
            &= e_{X-m-1}\big(\varphi_m(a + t + r + \epsilon) \big) \\
            &= e_{X-m-1}\big(\varphi_m(s')\big)
    \end{align*} 
    as claimed.
\end{proof}

\begin{lemma}\label{dominant leading monomial for monomials}
    If $M \in \Gamma_{X-m,0}$, then there is $s' \in K_{X-m-2}$ such that $\varphi_{m+1}\big(\varphi_m(M)\big) = e_{X-m-2}(s')$.
\end{lemma}

\begin{proof}
    Let $M = g \cdot \ell \cdot e_{X-m}(t)$ with $g \in G_{X-m}$, $\ell \in \Lambda_{X-m}$, and $t \in T_{X-m}$.
    Since 
        $$\varphi_{m+1}\big(\varphi_m(e_{X-m}(t))\big) = e_{X-m-2}\big(\varphi_{m+1}(\varphi_m(t))\big)$$
    and $\varphi_m$ and $\varphi_{m+1}$ are ring homomorphisms, it suffices to show that for any $x \in X-m$, $d \in \mathbb{N}$, $a,b \in \boldsymbol{k}$, we can find $s_1,s_2 \in K_{X-m-2}$ such that
    \begin{align*}
        \varphi_{m+1}\left(\varphi_m\left(E^{(d)}(x)^{a}\right)\right) 
            &= e_{X-m-2}(s_1) \\
        \varphi_{m+1}\left(\varphi_m\left(\log E'(x)^{b}\right)\right)
            &= e_{X-m-2}(s_2).
    \end{align*}
    
    First,
        $$\varphi_m\left(E^{(d)}(x)^a\right)= e_{X-m-1}(E(x-1))^aE'(x-1)^{da} \sum_{k = 0}^{\infty} \binom{a}{k} \delta_d(x)^k.$$
    Recall that $\mathrm{Supp}( \delta_d(x)) < \Gamma_{X-m-1,0}^{small}$, so 
        $$\mathrm{Supp} \left(\sum_{k = 1}^{\infty} \binom{a}{k} \delta_d(x)^k\right) < \Gamma_{X-m-1,0}^{small}.$$
    Thus $\mathrm{Supp} \left(\left(\sum_{k = 1}^{\infty} \binom{a}{k} \delta_d(x)^k\right)^l\right) < \Gamma_{X-m-1,0}^{small}$ for all $0 < l \in \mathbb{N}$.
    So we can express
    \begin{align*}
        E'(x-1)^{da} \sum_{k=0}^{\infty} \binom{a}{k} \epsilon^k&=e_{X-m-1}\left(da\log E'(x-1)+\sum_{l=1}^{\infty}\frac{(-1)^{l+1}}{l}\left(\sum_{k = 1}^{\infty} \binom{a}{k} \delta_d(x)^k\right)^l\right).
    \end{align*}
    By Lemma \ref{phi transfers exp}, we can take 
        $$s_1 = \varphi_{m+1}\left(aE(x-1) + da\log E'(x-1)+\sum_{l=1}^{\infty}\frac{(-1)^{l+1}}{l}\left(\sum_{k = 1}^{\infty} \binom{a}{k} \delta_d(x)^k\right)^l\right).$$
    
    Next we will find $s_2$. 
    We can calculate
    \begin{align*}
        \varphi_{m+1}\left(\varphi_m\left(\log E'(x)^b\right)\right) &= \varphi_{m+1}\left(E(x-1)^b \sum_{k = 0}^{\infty}\binom{b}{k}\left(\frac{\log E'(x-1)}{E(x-1)} \right)^k\right) \\
            &= e_{X-m-2}\big(E(x-2)\big)^b \sum_{k=0}^{\infty}\binom{b}{k} \frac{\big(E(x-2) + \log E'(x-2)\big)^k}{e_{X-m-2}\big(E(x-2)\big)^k}.
    \end{align*}
    Since $\mathrm{Supp} \left(\sum_{k=1}^{\infty} \binom{b}{k} \frac{\big(E(x-2) + \log E'(x-2)\big)^k}{e_{X-m-2}\big(E(x-2)\big)^k}\right) < \Gamma_{X-m-2,0}^{small}$, we can now bring this sum inside $e_{X-m-2}$ to get $\varphi_{m+1}\left(\varphi_m\left(\log E'(x)^b\right)\right)= e_{X-m-2}(s_2)$ where
    \begin{align*}
        s_2 &= bE(x-2) + \sum_{l=1}^{\infty}\frac{(-1)^{l+1}}{l}\left(\sum_{k=1}^{\infty}\binom{b}{k}\frac{\big(E(x-2) + \log E'(x-2)\big)^k}{e_{X-m-2}\big(E(x-2)\big)^k}\right)^l. \qedhere
    \end{align*}
\end{proof}

\begin{cor}\label{log hat = Lm(log)}
    The order defined in Definition \ref{monomial comp advanced} on $\Gamma_{X,0}$ using $\widehat{\log}$ matches the order given by defining 
        $$g \cdot \ell > e_{T_{X}}(t) \text{ if and only if } s' > \rho_{2}\big((\varphi_1 \circ \varphi_0)(t)\big)$$
    where $s'$ is such that $(\varphi_{1} \circ \varphi_0)(g \cdot \ell) = e_{X-2}(s')$.
\end{cor}

\begin{proof}
    Write $g \cdot \ell = \prod_{j=1}^p \log E'(x_j)^{b_j}E(x_j)^{a_{j,0}} \cdots E^{(d)}(x_j)^{a_{j,d}}$.
    Following Lemma \ref{dominant leading monomial for monomials}, for each $j=1,\dots,p$ find $s_{j,l}$ such that $\varphi_{1}(\varphi_0(E^{(l)}(x_j)^{a_{j,l}})) = e_{X-2}(s_{j,l})$, and find $s'$ such that $\varphi_{1}(\varphi_0(g \cdot \ell)) = e_{X-2}(s')$.
    Then $\Lt(s_{j,l}) = a_{j,l}e_{X-2}\big(E(x_j-2)\big)$ which matches $\widehat{\log}\big(E^{(d)}(x_j)\big) = a_{j,l}E_0(x-1)$.
    Also, since $s' - \sum_{j=1}^p\sum_{l=1}^d s_{j,l} \in K_{X-2,0}$, its terms do not come into play when comparing $s'$ and $\rho_{2}\big((\varphi_1 \circ \varphi_0)(t)\big)$.
    So 
    \begin{align*}
        s' > \rho_{2}\big((\varphi_1 \circ \varphi_0)(t)\big) \text{ if and only if } \widehat{\log}(\rho_0(g \cdot \ell)) > \rho_0(t) 
    \end{align*}
    since $\varphi_0$, $\varphi_1$, and $\rho_0$ are order preserving.
\end{proof}

\begin{lemma}\label{dominant leading monomial}
    Let $s \in \boldsymbol{k}[\Gamma_{X-m,0}^{small}]$ and let $d$ be the value of the largest derivative appearing in $s$.
    Let $\tilde{s} = (\varphi_{m+d} \circ \cdots \circ \varphi_m)(s)$.
    Then $\tilde{s}|_{\mathrm{Lv}(\tilde{s})}$ is a single term of the form 
        $$ce_{X-m-d-1}(\alpha)\prod_{j=1}^pE'(x_j-d-1)^{a_j}$$
    with $c \in \boldsymbol{k}^{\times}$ $\alpha \in A_{X-m-d-1,0} \oplus \cdots \oplus A_{X-m-d-1,d}$ and $a_1,\dots,a_p \in \boldsymbol{k}$.
\end{lemma}

\begin{proof}
    In Lemma \ref{phi is order preserving small}, we showed that $\varphi_{m,0}(s)|_{\mathrm{Lv}(\varphi_{m,0}(s))} \in \boldsymbol{k}[G_{X-m-1}]$, and 
    \begin{align*}
        \mathrm{Lm}\big(\rho_m(s)\big) = \mathrm{Lm}\Big(\sigma_{m+1}\left(\varphi_{m,0}(s)|_{\mathrm{Lv}(\varphi_{m,0}(s))}\right)\Big).
    \end{align*}
    Write
    \begin{align*}
        \mathrm{Lm}\big(\rho_m(s)\big) &= \prod_{j=1}^p E_0(x_j-1)^{a_{j,0}} E_1(x_j-1)^{a_{j,1}} \cdots E_d(x_j-1)^{a_{j,d}}.
    \end{align*}
    
    We now inductively define a sequence $s=s_0, s_1, \dots, s_{d}$ with $s_n \in \boldsymbol{k}[\Gamma_{X-m-n,0}^{small}]$ for $n=0,\dots,d$.
    Let $s_0 = s$.
    Let $s_1 = \frac{\varphi_{m,0}(s_0)|_{\mathrm{Lv}(\varphi_{m,0}(s_0))}}{\prod_{j=1}^p E(x_j-1)^{a_{j,0}}} \in \boldsymbol{k}[\Gamma_{X-m-1,0}^{small}]$, so that
    \begin{align*}
        \mathrm{Lm}\big(\sigma_{m+1}(s_1)\big)
            &= \prod_{j=1}^p E_1(x_j-1)^{a_{j,1}} \cdots E_d(x_j-1)^{a_{j,d}} \\
        \mathrm{Lm}\big(\rho_{m+1}(s_1)\big)
            &=\prod_{j=1}^p E_0(x_j-2)^{a_{j,1}}E_1(x_j-2)^{a_{j,1} + a_{j,2}} E_2(x_j-2)^{a_{j,3}} \cdots E_{d-1}(x_j-2)^{a_{j,d}}.
    \end{align*}
    
    Now given $s_n$ for $n<d$ such that
    \begin{align*}
        \mathrm{Lm}\big(\rho_{m+n}(s_{n})\big) 
            = \prod_{j=1}^p &E_0(x_j-n-1)^{a_{j,1} + \cdots + a_{j,n}} E_1(x_j-n-1)^{a_{j,1} + \cdots + a_{j,n+1}} \\
            &E_2(x-n-1)^{a_{j,n+2}} \cdots E_{d-n}(x_j-n-1)^{a_{j,d}}
    \end{align*}
    let 
        $$s_{n+1} = \frac{\varphi_{m+n,0}(s_n)|_{\mathrm{Lv}(\varphi_{m+n,0}(s_n))}}{\prod_{j=1}^p E(x_j-n-1)^{a_{j,1} + \cdots + a_{j,n}}}.$$
    Then $s_{n+1} \in \boldsymbol{k}[\Gamma_{X-m-n-1,0}^{small}]$ and 
    \begin{align*}
        \mathrm{Lm}\big(\sigma_{m+n+1}(s_{n+1})\big) &= \sigma_{m+n+1}\left(\frac{\varphi_{m+n,0}(s_n)|_{\mathrm{Lv}(\varphi_{m+n,0}(s_n))}}{\prod_{j=1}^p E_0(x_j-n-1)^{a_{j,1} + \cdots + a_{j,n}}}\right) \\
            &= \frac{\sigma_{m+n+1}\big(\varphi_{m+n,0}(s_n)|_{\mathrm{Lv}(\varphi_{m+n,0}(s_n))}\big)}{\prod_{j=1}^p E_0(x_j-n-1)^{a_{j,1} + \cdots + a_{j,n}}} \\
            &= \prod_{j=1}^p E_1(x_j-n-1)^{a_{j,1} + \cdots + a_{j,n+1}} \cdots E_{d-n}(x_j-n-1)^{a_{j,d}}
    \end{align*}
    by Lemma \ref{phi is order preserving small}.
    So
    \begin{align*}
        \mathrm{Lm}\big(\rho_{m+n+1}(s_{n+1})\big) 
            = \mathrm{Lm}&\big(\sigma_{m+n+1}(\mathrm{Lv}(\varphi_{m+n+1,0}(s_{n+1})))\big) \\
            = \prod_{j=1}^p &E_0(x_j-n-2)^{a_{j,1} + \cdots + a_{j,n+1}} E_1(x_j-n-2)^{a_{j,1} + \cdots + a_{j,n+1} + a_{j,n+2}} \\
                & E_2(x_j-n-2)^{a_{j,n+3}}\cdots E_{d-n-1}(x_j-n-2)^{a_{j,d}} 
    \end{align*}
    and we can continue the induction.
    
    Note that 
        $$\mathrm{Lm}\big(\rho_{m+d}(s_{d})\big) = \prod_{j=1}^p E_0(x_j-d-1)^{a_{j,1} + \cdots + a_{j,d}} E_1(x_j-d-1)^{a_{j,1} + \cdots + a_{j,d}}.$$
    Let $\mathrm{Init}(s_d)$ be as in Remark \ref{initial subsum of rho(s)}.
    Since $s_d \in \boldsymbol{k}[G_{m+d}]$ involves no log generators, and the only possible exponents of $E_l(x_j-d-1)$ for $l \ge 2$ in $\mathrm{Init}(s_d)$ are natural numbers, the preimage of $\mathrm{Init}(s_d)$ under $\sigma_{m+d+1}$ must be the single term
        $$c\prod_{j=1}^pE'(x_j-d-1)^{a_{j,1} + \cdots + a_{j,d}} = \varphi_{m+d}(s_d)|_{\mathrm{Lv}(\varphi_{m+d}(s_d))}$$
    for some $0 \ne c \in \boldsymbol{k}$.
    Let $t_0 = (\varphi_{m+d} \circ \cdots \circ \varphi_{m+1})\left(\prod_{j=1}^p E(x_j-1)^{a_{j,0}}\right) \in A_{X-m-d-1,d}$, and for $l=1,\dots, d-1$, let $t_l = (\varphi_{m+d} \circ \cdots \circ \varphi_{m+l+1})\left(\prod_{j=1}^p E(x_j-l-1)^{a_{j,1} + \cdots + a_{j,l}}\right) \in A_{X-m-d-1,d-l}$.
    Then
    \begin{align*}
        \tilde{s}|_{\mathrm{Lv}(\tilde{s})} &= (\varphi_{m+d} \circ \cdots \circ \varphi_m)(s|_{\mathrm{Lv}(s)})|_{\mathrm{Lv}(\tilde{s})} \\
            &= (\varphi_{m+d} \circ \cdots \circ \varphi_{m+1})\left(s_1\prod_{j=1}^p E(x_j-1)^{a_{j,0}}\right)\Bigg|_{\mathrm{Lv}(\tilde{s})} \\
            &= (\varphi_{m+d} \circ \cdots \circ \varphi_{m+1})(s_1)|_{\mathrm{Lv}(\tilde{s})} \cdot t_0 \\
            &= (\varphi_{m+d} \circ \cdots \circ \varphi_{m+2})\left(s_2\prod_{j=1}^p E(x_j-2)^{a_{j,1}}\right)\Bigg|_{\mathrm{Lv}(\tilde{s})} \cdot t_0 \\
            &= (\varphi_{m+d} \circ \cdots \circ \varphi_{m+2})(s_2)|_{\mathrm{Lv}(\tilde{s})} \cdot t_0 t_1 \\
            & \vdots \\
            &= \varphi_{m+d}\left(s_d\prod_{j=1}^p E(x_j-d)^{a_{j,1} + \cdots + a_{j,d}}\right)\Bigg|_{\mathrm{Lv}(\tilde{s})} t_0 \cdots t_{d-2} \\
            &=c\prod_{j=1}^pE'(x_j-d-1)^{a_{j,1} + \cdots + a_{j,d}} \cdot t_0 \cdots t_{d-1}.
    \end{align*}
    So $\tilde{s}|_{\mathrm{Lv}(\tilde{s})}$ is a single term of the specified form.
\end{proof}

The following lemma shows that we can eventually find a logarithm of any positive element of $K_{X-m}$, after applying enough embeddings.

\begin{lemma}\label{log closure}
    For all $s \in (K_{X-m})^{>0}$, there is $l \in \mathbb{N}$ and $s' \in K_{X-m-l-1}$ such that
        $$(\varphi_{m+l} \circ \cdots \circ \varphi_m)(s) = e_{X-m-l-1}(s').$$
\end{lemma}

\begin{proof}
    Let $\mathrm{Lm}(s) = e_{X-m}(\alpha + t)\cdot g \cdot \ell$ with $\alpha \in A_{X-m,n} \oplus \cdots \oplus A_{X-m,0}$, $t \in T_{X-m}$, $g \in G_{X-m}$, and $\ell \in \Lambda_{X-m}$.
    Let $x_1 > \cdots > x_p$ list the elements of $X$ that appear in $g$ or $\ell$ and let $a_j = \xi_{x_j}(g)$ for $j=1,\dots,p$.
    Then we can write
    \begin{align*}
        s &= s_0 \cdot e_{X-m}(\alpha + t)\prod_{j=1}^pE(x_j)^{a_j} 
    \end{align*}
    with $\mathrm{Lv}(s_0) = \Gamma_{X-m,0}^{small}$.
    Let $d$ be the largest derivative that appears in $s_0|_{\mathrm{Lv}(s_0)}$.
    Then 
    \begin{multline*}
        (\varphi_{m+d} \circ \cdots \circ \varphi_m)\left(e_{X-m}(\alpha + t)\prod_{j=1}^pE(x_j)^{a_j} \right) = \\
            e_{X-m-d-1}\left((\varphi_{m+d} \circ \cdots \circ \varphi_m)(\alpha + t) + \sum_{j=1}^pa_j(e_{X-m-d-1})^{\circ (d+1)}\big(E(x_j-d-1)\big)\right).
    \end{multline*}
    By Lemma \ref{dominant leading monomial}, we can write 
    \begin{align*}
        (\varphi_{m+d} \circ \cdots \circ \varphi_m)(s_0) = \left(ce_{X-m-d-1}(\alpha')\prod_{j=1}^pE'(x_j-d-1)^{a_{j,1} + \cdots + a_{j,d}}\right)(1 + \epsilon)
    \end{align*}
    with $\mathrm{Supp}(\epsilon) < \Gamma_{X-m-d-1,0}^{small}$.
    Note that $c > 0$ since $s > 0$ and all $\varphi_{k}$ are order preserving for $k \in \mathbb{N}$.
    So $\log c \in \boldsymbol{k}$ is defined, and we can express $(\varphi_{m+d} \circ \cdots \circ \varphi_m)(s) = e_{X-m-d}(s')$ where
    \begin{align*}
        s' = \log c &+ \alpha' + \sum_{j=1}^p (a_{j,1} + \cdots + a_{j,d})\log E'(x_j-d-1) + \sum_{k=1}^{\infty} \frac{(-1)^{k+1}}{k}\epsilon^k \\
            &+ (\varphi_{m+d} \circ \cdots \circ \varphi_m)(\alpha + t) + \sum_{j=1}^pa_j(e_{X-m-d-1})^{\circ (d+1)}\big(E(x_j-d-1)\big).
    \end{align*}
    Since $\mathrm{Supp}(\epsilon) < \Gamma_{X-m-d-1,0}^{small}$, we have $\sum_{k=1}^{\infty} \frac{(-1)^{k+1}}{k}\epsilon^k \in K_{X-m-d-1}$. 
\end{proof}

\begin{lemma}\label{exp(exp(s)) is defined}
    If $s \in K_{X-m}$ and $e_{X-m}(s)$ is defined, so are $e_{X-m}\big(\pm e_{X-m}(s)\big)$.
\end{lemma}

\begin{proof}
    Suppose $s \in B_{X-m,n}^*$, so we can write $s = a + t + c + \epsilon$ with $\mathrm{Supp} (a) > \Gamma_{X-m,0}^{small}$ (or $a=0$), $t \in T_{X-m}$, $c \in \boldsymbol{k}$, and $\mathrm{Supp} (\epsilon) < \Gamma_{X-m,0}^{small}$ (or $\epsilon = 0)$.
    Note that $e_{X-m}(a + t)\epsilon^k \not \in \Gamma_{X-m,0}^{small}$ for all $k \in \mathbb{N}$.
    So we can separate
    \begin{align*}
        e_{X-m}(a + t + c + \epsilon) = e_{X-m}(c)\sum_{k=0}^{\infty} \frac{e_{X-m}(a + t)\epsilon^k}{k!} 
            = a' + \epsilon' 
    \end{align*}
    with $\mathrm{Supp} (a') > \Gamma_{X-m,0}^{small}$ and $\mathrm{Supp} (\epsilon') < \Gamma_{X-m,0}^{small}$.
    So $e_{X-m}(a + t + c + \epsilon) \in B_{X-m,n+1}^*$.
    The same argument shows $e_{X-m}(-a - t - c - \epsilon) \in B_{X-m,n+1}^*$.
\end{proof}

\begin{cor}\label{exp closure}
    For all $s \in K_{X-m}$, there is $l \in \mathbb{N}$ such that 
        $$e_{X-m-l-1}\big((\varphi_{m+l} \circ \cdots \circ \varphi_m)(s)\big)$$
    is defined.
\end{cor}

\begin{proof}
    Let $s \in K_{X-m}$, and write $s = u|s|$ with $u = \pm 1$.
    By Lemma \ref{log closure}, we know there is some $l \in \mathbb{N}$ and $s' \in K_{X-m-l-1}$ such that
        $$(\varphi_{m+l} \circ \cdots \circ \varphi_m)(|s|) = e_{X-m-l-1}(s').$$
    By Lemma \ref{exp(exp(s)) is defined}, $e_{X-m-l-1}\big(u \cdot e_{X-m-l-1}(s')\big)$ is defined.
    So
    \begin{align*}
        e_{X-m-l-1}\big(u \cdot e_{X-m-l-1}(s')\big) 
            &= e_{X-m-l-1}\big(u(\varphi_{m+l} \circ \cdots \circ \varphi_m)(|s|)\big) \\
            &= e_{X-m-l-1}\big((\varphi_{m+l} \circ \cdots \circ \varphi_m)(s)\big)
    \end{align*}
    is defined.
\end{proof}

\begin{lemma}\label{R_an closure lemma}
    Let $s_1,\dots,s_n \in K_{X-m}$ with $\mathrm{Supp} (s_i) < 1$.
    Then there is $l \in \mathbb{N}$ such that 
        $$\sum_{j_1,\dots,j_n = 0}^{\infty} a_{j_1,\dots,j_n}\prod_{i=1}^n(\varphi_{m+l} \circ \cdots \circ \varphi_{m})(s_i)^{j_i} \in K_{X-m-l-1}$$
    for any $(a_{j_1,\dots,j_n})_{j \in \mathbb{N}^n} \in (\boldsymbol{k}_{\alpha})^{\omega}$, for any $\alpha \in \mathcal{A}$.
\end{lemma}

\begin{proof}
    It suffices to show that for some $l \in \mathbb{N}$, 
        $$\mathrm{Supp}\big((\varphi_{m+l} \circ \cdots \circ \varphi_{m})(s_i)\big) < \Gamma_{X-m-l-1,0}^{small}$$
    for all $i=1,\dots,n$.
    If $\mathrm{Supp}\big((\varphi_{m+l} \circ \cdots \circ \varphi_{m})(s_i)\big) < \Gamma_{X-m-l-1,0}^{small}$ for some $l$, then $\mathrm{Supp}\big((\varphi_{m+l'} \circ \cdots \circ \varphi_{m})(s_i)\big) < \Gamma_{X-m-l'-1,0}^{small}$ for all $l' > l$.
    So we will find $l_i$ that works for $i$ and take $l = \max(l_1,\dots,l_n)$.
    
    If $\mathrm{Supp} (s_i) < \Gamma_{X-m,0}^{small}$, then we can use $l=0$, so we may assume $\mathrm{Lv}(s_i) = \Gamma_{X-m,0}^{small}$.
    Let $d_i$ be the largest derivative that appears in $s_i|_{\mathrm{Lv}(s_i)}$.
    Let $x_1 > \cdots > x_p$ list the elements of $X$ that appear in $s_1|_{\mathrm{Lv}(s_1)}, \dots, s_n|_{\mathrm{Lv}(s_n)}$.
    Let $\widetilde{s_i} = (\varphi_{m+d} \circ \cdots \circ \varphi_m)(s_i)$.
    Since $\mathrm{Supp} (s_i) < 1$, we must have $\mathrm{Supp}(\widetilde{s_i}) < 1$.
    By Lemma \ref{dominant leading monomial}, we can write
        $$\widetilde{s_i}|_{\mathrm{Lv}(\widetilde{s_i})} = c_i e_{X-m-d-1}(\alpha_i)\prod_{j=1}^pE'(x_j-d-1)^{a_{i,j}}$$
    with $\alpha_i \in A_{X-m-d-1,0} \oplus \cdots \oplus A_{X-m-d-1,d}$.
    No monomial of this form can be in $\Gamma_{X-m-d-1,0}^{small}$, so we must have $\mathrm{Supp} (\widetilde{s_i}) < \Gamma_{X-m-d-1,0}^{small}$.
\end{proof}

\begin{cor}\label{R_an closure}
    Let $s_1,\dots,s_n \in K_{X-m}$ with $\mathrm{Supp} (s_i) < 1$.
    Let $f \in \mathbb{R}\{X_1,\dots,X_n\}$, and let $(r_1,\dots,r_n) \in \boldsymbol{k}_{\alpha} \cap [-1,1]^n$ for some $\alpha \in \mathcal{A}$.
    Then there is $l \in \mathbb{N}$ such that 
        $$f\big(r_1 + (\varphi_{m+l} \circ \cdots \circ \varphi_{m})(s_1), \dots, r_n + (\varphi_{m+l} \circ \cdots \circ \varphi_{m})(s_n)\big) \in K_{X-m-l-1}.$$
\end{cor}

\begin{proof}
    The coefficient of any monomial 
        $$\prod_{i=1}^n(\varphi_{m+l} \circ \cdots \circ \varphi_{m})(s_i)^{j_i}$$
    can be shown to converge because $f$ converges in a neighborhood of $[-1,1]^n$.
\end{proof}

\begin{cor}\label{field}
    For all $s \in K_{X-m}$, there is $l \in \mathbb{N}$ such that $(\varphi_{m+l} \circ \cdots \circ \varphi_{m})(s)$ has a multiplicative inverse in $K_{X-m-l-1}$.
\end{cor}

\begin{proof}
    We can write $s = \mathrm{Lm}(s)(1 + \epsilon)$, with $\mathrm{Supp} (\epsilon) < 1$.
    By Lemma \ref{R_an closure}, there is $l \in \mathbb{N}$ such that
        $$\sum_{k=0}^{\infty} \big(-(\varphi_{m+l} \circ \cdots \circ \varphi_{m})(\epsilon)\big)^k$$
    is defined in $K_{X-m-l-1}$.
    Then
        $$(\varphi_{m+l} \circ \cdots \circ \varphi_{m})\left(\frac{1}{\mathrm{Lm}(s)}\right) \cdot \sum_{k=0}^{\infty} \big(-(\varphi_{m+l} \circ \cdots \circ \varphi_{m})(\epsilon)\big)^k$$
    the multiplicative inverse of $(\varphi_{m+l} \circ \cdots \circ \varphi_{m})(s)$ in $K_{X-m-l-1}$.
\end{proof}

\begin{defn}
    Let $D_X(\boldsymbol{k})$ be the direct limit of the directed system 
        $$(K_{X-m}, \varphi_{m+l} \circ \cdots \circ \varphi_m)_{l,m \in \mathbb{N}}.$$
\end{defn}

By Lemmas \ref{log closure} and \ref{exp closure} and Corollaries \ref{R_an closure} and \ref{field}, $D_X(\kk)$ can be made into a model of $T_{\mathrm{an}}(\exp,\log)$.

\begin{lemma}\label{log-exp bounded}
    For all $s \in D_X(\boldsymbol{k})_{+\infty}$ and $x \in X$, there are $n_1,n_2 \in \mathbb{N}$ such that
        $$\log^{\circ n_1}E(x) < s < \exp^{\circ n_2}E(x).$$
\end{lemma}

\begin{proof}
    Identify $s$ with a corresponding element $s'$ of some $K_{X-m,n}$, and let $x \in X$.
    Then $(e_{X-m})^{\circ k}\big(E(x-m)\big) \in A_{X-m,k}$ for all $k \in \mathbb{N}$.
    So
        $$s' < (e_{X-m})^{\circ (n+1)}\big(E(x)\big) \in e_{X-m}(A_{X-m,n})$$
    and thus $s' < \exp^{\circ (n+1-m)}E(x)$ in $D_X(\boldsymbol{k})$.
    

    Now we must find $n_1$ such that
        $$\log^{\circ n_1}E(x) < s.$$
    As in the proof of Lemma \ref{log closure}, write
        $$s' = s_0 \cdot e_{X-m}(\alpha + t)\prod_{j=1}^p E(x_j-m)^{a_j}$$
    with $\alpha \in A_{X-m,n} \oplus \cdots \oplus A_{X-m,0}$, $t \in T_{X-m}$, and $\mathrm{Lv}(s_0) = \Gamma_{X-m,0}^{small}$.
    Since $s$ is positive and infinite, so is $s'$.
    
    We now split into cases based on the form of $\mathrm{Lm}(s')$.
    First, if $\alpha \ne 0$, then we must have $\alpha > 0$ since $s'$ is positive and infinite.
    Then
        $$s' > E(x-m) \in K_{X-m,0}$$
    so $s > \log^{\circ m}E(x)$.
    
    Second, if $\alpha = 0$ and $e_{X-m}(t)\prod_{j=1}^p E(x_j-m)^{a_j} \ne 1$, then we must have 
        $$e_{X-m}(\alpha + t)\prod_{j=1}^p E(x_j-m)^{a_j} > 1$$
    again since $s'$ is infinite.
    Let $y$ be the smallest element of $X$ appearing in $\mathrm{Lm}(t)$ (or $y=+\infty$ if $t \ne 0)$, and let $x_* = \min(y,x_p)$.
    We can find $0 < a \in \boldsymbol{k}$ such that
        $$E(x_*-m)^a < e_{X-m}(\alpha + t)\prod_{j=1}^p E(x_j-m)^{a_j}.$$
    If $t \ne 0$ or if $p > 1$, then $a = 1$ works.
    If $t = 0$ and $x_* = x_1$, then let $a = \frac{|a_1|}{2}$.
    Now since $\mathrm{Lv}(s_0) = \Gamma_{X-m,0}^{small}$, we must have $s_0 > E(x_*-m)^{-a/2}$.
    So
    \begin{align*}
        s' > E(x_*-m)^{a} E(x_*-m)^{-a/2} = E(x_*-m)^{a/2}.
    \end{align*}
    And thus $s > E(x_*-m)^{a/2} > E(x-m-1) = \log^{\circ (m+1)}E(x)$.
    
    Third, suppose $s' = s_0$.
    Let $d$ be the largest derivative appearing in $s_0|_{\mathrm{Lv}(s_0)}$.
    Let $\widetilde{s_0} = (\varphi_{m+d} \circ \cdots \circ \varphi_m)(s_0)$, which must be positive and infinite since $s'$ is.
    By Lemma \ref{dominant leading monomial}, we can write
        $$\widetilde{s_0}|_{\mathrm{Lv}(\widetilde{s_0})} = ce_{X-m-d-1}(\alpha)\prod_{j=1}^pE'(x_j-d-1)^{a_j}$$
    with $c > 0$, $\alpha \in A_{X-m-d-1,0} \oplus \cdots \oplus A_{X-m-d-1,d}$, and $a_j \in \boldsymbol{k}$.
    Since $\widetilde{s_0} \in K_{X-m-d-1}$ is of a form handled by either the first or second case, we get either $\widetilde{s_0} > E(x-m-d-1)$ or $\widetilde{s_0} > E(x_* -m-d-1)^{a/2}$.
    We have
        $$s > E(X-m-d-2) = \log^{\circ (m+d+2)}E(x)$$
    in either case.
\end{proof}

\section{An \texorpdfstring{$\mathcal{L}_{\mathrm{transexp}}$}{Ltransexp} differential series field}

Let $F \vDash T_{\mathrm{transexp}}$.
We will build an increasing sequence $(H_i : i \in \mathbb{N})$ of ordered log-exp differential fields, starting with $H_0 = F((\tau^{-1}))^{le}$ where we take $\tau > F$ for the ordering.
We will build this sequence so that 
    $$M_F := \bigcup_{i \in \mathbb{N}} H_i$$
is a $\mathcal{L}_{\mathrm{an}}(\exp,\log)$-structure closed under $E$, its derivatives $E^{(d)}$, and their functional inverses.
We will then define a derivation on $M_F$ that works like differentiation with respect to $\tau$.
Each $H_{i+1}$ will be constructed from $H_i$ using the constructions in Subsections \ref{subsection 4.1} and \ref{Subsection 4.2} to add new monomials for $E$, its derivatives, and $L$ applied to certain elements of $H_i$.

We do not need to create new monomials for $E$, its derivatives, and $L$ applied to \textit{all} elements of $H_i$.
For example, there is no need to add new monomials for $E$ composed with elements in the same $\mathbb{Z}$-orbit---if we include $E(x)$ when building $H_1$, then we automatically have expressions that represent $E(x + k)$, $k \in \mathbb{Z}$.
Also, if $x,y \in H_i$ are ``too close" in the sense that $x > y$ but $E(x) < E(y)^a$ for some $a$, then we cannot add both $E(x)$ and $E(y)$ as new elements over a field of coefficients containing $a$ because the separation assumptions of Remark \ref{order props} would not be satisfied.
This will not be an issue because we will be able to express $E(x)$ and $E(y)$ in terms of each other in that case.

\subsection{Constructing \texorpdfstring{$H_{i+1}$}{Hi+1} from \texorpdfstring{$H_i$}{Hi}}\label{subsection 5.1}

Suppose $H_0, \dots, H_i$ have been constructed.
If $i > 0$, assume we have also constructed order-preserving embeddings $\iota_j : H_j \to H_{j+1}$ for $j = 0,\dots, i-1$ such that if $j > 0$, $z \in H_{j-1}$, and $E(z)$ is defined in $H_{j}$, then 
    $$\iota_{j}\big(E(z)\big) = E\big(\iota_{j-1}(z)\big).$$
For each $j=0,\dots,i$, let
    $$\mathrm{Fin}(H_j) = \left\{f \in H_j : \exists n \in\mathbb{N} \big(|f| \le n\big)\right\}.$$

Since $\mathrm{Fin}(H_j)$ is a convex subgroup of $H_j$, the quotient $H_j/\mathrm{Fin}(H_j)$ inherits an order from $H_j$.
For each $j=1,\dots,i-1$, $\iota_j(H_j)/\mathrm{Fin}\left(\iota_j\left(H_j\right)\right)$ is a (non-convex) subgroup of $H_{j+1}/\mathrm{Fin}\left(H_{j+1}\right)$.

\begin{defn}
    If $s \in H_{j+1}$ and $s > F$, define $\lambda(s)$ to be the unique (if it exists) one of
    \begin{enumerate}
        \item $\beta = \big\{E\big(L^{\circ (j+2)}(\tau)\big)\big\}$ or
        
        \item $\beta \in H_j/\mathrm{Fin} (H_j)$
    \end{enumerate}
    such that for some $z \in \beta$ and $n \in \mathbb{N}$, $E(z\pm n)$ are defined in $H_{j+1}$ and 
        $$E(z-n) < s < E(z+n).$$
\end{defn}
Suppose also that for $j = 0,\dots,i-1$ and each $s \in (H_{j+1})_{>F}$, $\lambda(s)$ exists. 
After defining $H_{i+1}$, we will prove by induction in Lemma \ref{lambda exists} that $\lambda(s)$ exists for all $s \in (H_{i+1})_{>F}$.

\begin{rmk}
    The goal of this subsection is to construct $H_{i+1}$ from $H_i$ as follows:
    Any coset $\alpha \in H_i/\mathrm{Fin}(H_i)$ with $\alpha > F$, viewed as a subset of $H_i$, consists only of positive elements of $H_i$ that are infinite relative to $F$.
    For each such $\alpha$, we will define $X_{\alpha} \subseteq \alpha$ and create new monomials $E^{(d)}(x)$ and $\log E'(x)$ for all $x \in X_{\alpha}$ and $d \in \mathbb{N}$.
    We will do this by building a field $C_{\bar{\alpha},i} \vDash T_{\mathrm{an}}(\exp,\log)$ for each finite increasing sequence 
        $$\bar{\alpha} = \{\alpha_0<\alpha_1<\dots<a_k\} \subset \big(H_i/\mathrm{Fin}(H_i)\big)_{>F}$$
    using the constructions in Subsections \ref{subsection 4.1} and \ref{Subsection 4.2}.
    We will then define $H_{i+1}$ to be the direct limit of $\big(C_{\bar{\alpha},i} : \bar{\alpha} \subset (H_i/\mathrm{Fin}(H_i))_{>F}\big)$, and show that for each infinite (relative to $F$) $s \in H_{i+1}$, $\lambda(s)$ exists. 
    Finally, we will define an order preserving embedding $\iota_i : H_i \to H_{i+1}$.
\end{rmk}

\subsubsection{Constructing $C_{\bar{\alpha},i}$} 

First suppose $\bar{\alpha}$ is the empty sequence.
We will build $C_{\varnothing,i}$ in two steps.
First, let $X_i' = \{L^{\circ(i+2)}(\tau) \}$, a single element set.
Note that $X_i'$ and $F$ trivially satisfy the ordering and separation assumptions of Remark \ref{order props}.
Since $F \vDash T_{\mathrm{an}}(\exp,\log)$, we can build $D_{X_i'}(F) \vDash T_{\mathrm{an}}(\exp,\log)$, using the one-element family $\{F\}$ of subfields of $F$.

Second, let $X_i= \{L^{\circ(i+1)}(\tau) \}$, another single element set.
We would like to define $C_{\varnothing,i} := D_{X_i}(D_{X_i'}(F))$, so we must check that $X_i$ and $D_{X_i'}(F)$ satisfy the ordering and separation assumptions of Remark \ref{order props}.
Most of the assumptions are satisfied trivially, but we must show that 
    $$E\big(L^{\circ(i+1)}(\tau)-m\big) > D_{X_i'}(F)$$
for all $m \in \mathbb{N}$.
We have not yet defined how $E\big(L^{\circ(i+1)}(\tau)-m\big)$ compares to elements of $D_{X_i'}(F)$, but by Lemma \ref{log-exp bounded}, any element $s \in D_{X_i'}(F)$ is bounded above by 
    $$\exp^{\circ n}E\big(L^{\circ(i+2)}(\tau)\big) = \exp^{\circ n}L^{\circ(i+1)}(\tau)$$
for some $n \in \mathbb{N}$.
Since we want $E > \exp^{\circ l}$ for any $l \in \mathbb{N}$, we extend the ordering so that 
\begin{align*}
    E\big(L^{\circ(i+1)}(\tau)-m\big) 
        &= \log^{\circ m} E\big(L^{\circ(i+1)}(\tau)\big) \\
        &> \log^{\circ m}\big(\exp^{\circ (n+m)}(L^{\circ (i+1)}(\tau))\big) \\
        &= \exp^{\circ n}\big(L^{\circ (i+1)}(\tau)\big) \\
        &> s
\end{align*}
for all $m \in \mathbb{N}$. 
Since the assumptions of Remark \ref{order props} are now satisfied, we may take 
    $$C_{\varnothing,i} := D_{X_i}(D_{X_i'}(F)) \vDash T_{\mathrm{an}}(\exp,\log).$$
again using the one-element family $\{D_{X_i'}(F)\}$ of subfields of $D_{X_i'}(F)$.

Next, for $\bar{\alpha} = \{\alpha_0 < \cdots < \alpha_k\}$ an increasing sequence in $\big(H_i/\mathrm{Fin}(H_i)\big)_{>F}$ and $\beta > \bar{\alpha}$, we will define the sets $C_{\bar{\alpha}\frown \beta,i}$ to be $D_{X_{\beta}}(C_{\bar{\alpha},i})$ where $X_{\beta}$ depends on $\beta$ as described in Definition \ref{define X_alpha}.
We use the following family of subfields of $C_{\bar{\alpha},i}$ in the construction of $D_{X_{\beta}}(C_{\bar{\alpha},i})$ from Subsection \ref{Subsection 4.2}:
\begin{align*}
    \kappa_{\varnothing,i} 
        &:= \{C_{\varnothing,i}\} 
        & \text{if } \bar{\alpha} = \varnothing\\
    \kappa_{\alpha_0,i} 
        &:= \{D_X(C_{\varnothing,i}) : X \subset X_{\alpha_0} \text{ finite}\} 
        &\text{if } k=0\\
    \kappa_{\{\alpha_0 < \cdots < \alpha_{k}\},i} 
        &:= \{D_X(\boldsymbol{k}) : X \subset X_{\alpha_l} \text{ finite}, \boldsymbol{k} \in \kappa_{\{\alpha_0 < \cdots < \alpha_k\},i}\}
        & \text{if } k>0.
\end{align*}

Next, we define $X_{\beta}$ for $\beta \in \big(H_i/\mathrm{Fin}(H_i)\big)_{>F}$ as follows: 

\begin{defn}\label{define X_alpha}
    If $i > 0$ and $\beta \in \big(\iota_{i-1}(H_{i-1})/\mathrm{Fin}(\iota_{i-1}(H_{i-1}))\big)_{>F}$, then let 
        $$X_{\beta}' = \iota_{i-1}\big(X_{(\iota_{i-1})^{-1}(\beta)}\big).$$
    If $i > 0$ and $E\big(L^{\circ(i+1)}\big) \in \beta$, let $X_{\beta}' = \big\{E\big(L^{\circ(i+1)}\big)\big\}$.
    Otherwise, let $X_{\beta}' = \varnothing$.
    
    Extend $X_{\beta}'$ to a maximal set $X_{\beta}$ of representatives from $\beta$ satisfying the following conditions:
    \begin{enumerate}
        \item The elements of $X_{\beta}$ are not too close:
        Suppose $i>0$, $x > y$, and $\frac{1}{x-y}$ is infinite.
        If 
            $$\beta \le \lambda\left(\frac{1}{x-y}\right)$$
        then at most one of $x,y$ is in $X_{\beta}$.
    
        \item The elements of $X_{\beta}$ not too far apart: there exists $r : X_{\beta} \times X_{\beta} \to \mathbb{Q} \cap (0,1)$ such that for all $x,y \in X_{\beta}$ with $x > y$, we have $x-y < r(x,y)$
    \end{enumerate}
\end{defn}

\begin{rmk}
    The intuition for why the first condition means $x$ and $y$ are ``not too close" comes from its equivalence via Lemma \ref{condition 1} to the ordering and separation assumption from Remark \ref{order props} requiring that $E(x-m) > E(y-m)^a$ whenever $x>y$, $m \in \mathbb{N}$, and $a$ is in the field $\boldsymbol{k}$ of coefficients and exponents.
    We use the equivalent condition here to make it clear that the same set $X_{\beta}$ will satisfy the ordering and separation assumptions regardless of which field of coefficients we pair it with.
\end{rmk}

\begin{rmk}
    The cases used to define $X_{\beta}'$ in Definition \ref{define X_alpha} are disjoint.
    Every infinite monomial in $H_{i-1}$ is bounded below by $\log^{\circ l}(\tau)$ for some $l \in \mathbb{N}$ if $i=1$, and by 
        $$\log ^{\circ l} E\big(L^{\circ i}(\tau)\big) \in D_{X_{i-2}'}(F)$$
    for some $l \in \mathbb{N}$ if $i > 1$ by Lemma \ref{log-exp bounded}. 
    If $E\big(L^{\circ(i+1)}\big) \in \beta$, then 
        $$\beta < \big(\iota_{i-1}(H_{i-1})/\mathrm{Fin}(\iota_{i-1}(H_{i-1}))\big)_{>F}.$$
    So the different definitions of $X_{\beta}'$ do not conflict.
\end{rmk}

If $i > 0$ and $\beta \in \iota_{i-1}(H_{i-1})/\mathrm{Fin}(\iota_{i-1}(H_{i-1}))$, we must check that $X_{\beta}'$ satisfies the two conditions in Definition \ref{define X_alpha}, so that it is possible to extend it to a maximal set $X_{\beta}$ that satisfies these conditions.
\begin{enumerate}
    \item To see that $X_{\beta}'$ satisfies Condition (1), let $x > y \in X_{\beta}'$ and suppose $\frac{1}{x-y}$ is infinite.
    
    If $i=1$, then we must have $\iota_{0}^{-1}\left(\frac{1}{x-y}\right) \in H_0$.
    Then 
        $$\lambda\left(\frac{1}{x-y}\right) = \big\{E\big(L(L(\tau))\big)\big\}$$
    since $H_0$ is exponentially bounded and we identify $E\big(E(L(L(\tau)))+n\big)$ with $\exp^{\circ n}(\tau)$ for all $n \in \mathbb{Z}$. 
    So we have $\beta > \lambda\left(\frac{1}{x-y}\right)$.
    
    If $i > 1$, then $(\iota_{i-1})^{-1}(x),(\iota_{i-1})^{-1}(y) \in X_{(\iota_{i-1})^{-1}(\beta)}$ means that
        $$(\iota_{i-2})^{-1}(\beta) > \lambda\left(\frac{1}{(\iota_{i-1})^{-1}(x-y)}\right).$$
    Since $\iota_{i-1}\big(E(z)\big) = E\big(\iota_{i-2}(z)\big)$, we must have $\iota_{i-2}\big(\lambda(s)\big) = \lambda\big(\iota_{i-1}(s)\big)$ because 
    \begin{align*}
        E(z-n) < s < E(z+n) \text{ if and only if }& \iota_{i-1}\big(E(z-n)\big) < \iota_{i-1}(s) < \iota_{i-1}\big(E(z+n)\big) \\
            \text{ i.e., }& E\big(\iota_{i-2}(z)-n\big) < \iota_{i-1}(s) < E\big(\iota_{i-2}(z)-n\big).
    \end{align*}
    So $\beta > \iota_{i-2}\left(\lambda\left(\frac{1}{(\iota_{i-1})^{-1}(x-y)}\right)\right) = \lambda\left(\frac{1}{x-y}\right)$.
    
    \item Since $X_{(\iota_{i-1})^{-1}(\beta)}$ satisfies Condition (2) and $\iota_{i-1}$ is order preserving, $X_{\beta}'$ satisfies Condition (2).
\end{enumerate}
Since $X_{\beta}'$ satisfies Conditions (1) and (2) in Definition \ref{define X_alpha}, it can be extended to a maximal set $X_{\beta}$ that satisfies the two conditions.

In order to define $C_{\bar{\alpha} \frown \beta,i} := D_{X_{\beta}}(C_{\bar{\alpha},i})$
we must check that the ordering and separation assumptions of Remark \ref{order props} are satisfied in each case.
Any $X_{\beta}$ is a subset of an ordered field, and each $C_{\bar{\alpha},i}$ is constructed to be an ordered exponential field.
Note that Condition (2) in Definition \ref{define X_alpha} of the sets $X_{\beta}$ matches the final ordering and separation assumption of Remark \ref{order props}. 
So all that remains to show are assumptions (3a) and (3b).

First, we will check that for all $x \in X_{\beta}$ and $m \in \mathbb{N}$, we have $E(x-m) > C_{\bar{\alpha},i}$.

\begin{lemma}
    For all $x \in X_{\beta}$ and $m \in \mathbb{N}$, we have $E(x-m) > C_{\bar{\alpha},i}$.
\end{lemma}

\begin{proof}
    Let $s \in C_{\bar{\alpha},i}$.
    
    First suppose $\bar{\alpha} = \varnothing$.
    Since $x \in X_{\beta} \subset \beta \in \big(H_i/\mathrm{Fin}(H_i)\big)_{>F}$, $x$ is positive and infinite.
    The smallest positive infinite elements of $H_i$ come from $D_{X_{i-1}'}(F)$, where $X_{i-1}' = \{L^{\circ (i+1)}(\tau)\}$.
    Let $z = L^{\circ (i+1)}(\tau)$.
    By Lemma \ref{log-exp bounded}, there is some $l \in \mathbb{N}$ such that 
        $$x > \log^{\circ l}E(z).$$
    Recall that $X_i = \{L^{\circ (i+1)}(\tau)\} = \{z\}$ as well.
    So by Lemma \ref{log-exp bounded} there is some $n \in \mathbb{N}$ such that 
        $$s < \exp^{\circ n}E(z).$$
    We have not yet defined how $E(x-m)$ compares to elements of $C_{\varnothing,i}$, but since we want $E > \exp^{\circ j}$ for all $j \in \mathbb{N}$, we extend the ordering so that
    \begin{align*}
        E(x-m) &= \log^{\circ m}E(x) \\
            &> \log^{\circ m}E\big(\log^{\circ l}E(z)\big) \\
            &> \log^{\circ m}\exp^{\circ (n + m + l)}\big(\log^{\circ l}E(z)\big) \\
            &= \exp^{\circ n}E(z) \\
            &>s
    \end{align*}
    for all $m \in \mathbb{N}$.
    
    If $\bar{\alpha} \ne \varnothing$, let $s \in C_{\bar{\alpha},i}$.
    By Lemma \ref{log-exp bounded}, there is some $z \in X_{\alpha_{k}}$ and $n \in \mathbb{N}$ such that $s < \exp^{\circ n}E(z)$.
    We have not yet defined how $E(x-m)$ compares to elements of $C_{\bar{\alpha},i}$, but since $\beta > \alpha_{k}$ in $H_i/\mathrm{Fin}(H_i)$, we know $x - m > z + n$ for any $m \in \mathbb{N}$.
    So we extend the ordering to have
        $$E(x - m) > E(z + n) = \exp^{\circ n}E(z) > s$$
    for all $m \in \mathbb{N}$.
\end{proof}

To finish checking the assumptions of Remark \ref{order props}, we must show that 
    $$E(x-m) > E(y-m)^a$$
for all $m \in \mathbb{N}$, $x,y \in X_{\beta}$ with $x > y$, and $a \in C_{\bar{\alpha},i}$.
We prove this from the first condition on $X_{\beta}$ from Definition \ref{define X_alpha}.

\begin{lemma} \label{condition 1}
    Let $i \in \mathbb{N}$ and let $\beta$ be larger than all elements of $\bar{\alpha} \subset (H_i / \mathrm{Fin}(H_i))_{>F}$, so that $E(x-m) > C_{\bar{\alpha},i}$ for all $x \in X_{\beta}$ and $m \in \mathbb{N}$.
    Then the following are equivalent:
    \begin{enumerate}
        \item $X_{\beta}$ satisfies Condition (1) of Definition \ref{define X_alpha}.
        
        \item For $x,y \in X_{\beta}$ with $x > y$, $m \in \mathbb{N}$, and $a \in C_{\bar{\alpha},i}$, we have $E(x-m) > E(y-m)^a$.
    \end{enumerate}
\end{lemma}

\begin{proof}
    Let $x,y \in X_{\beta}$ with $x > y$.
    Before proving the equivalence in the statement of the lemma, we will need to extend the partial order on expressions involving $E$, $L$, and log to make sense of ``$E(x-m) > E(y-m)^a$": for any $m \in \mathbb{N}$ and any $a \in (C_{\bar{\alpha},i})_{>1}$, we can have
    \begin{align*}
        E(x-m) > E(y-m)^a &\text{ if and only if } E(x-m-2) > E(y-m-2)+\log a\\
            &\text{ if and only if } x-m-2 > L\big(E(y-m-2)+\log a\big) \\
            &\text{ if and only if } x-y > L\big(E(y-m-2)+\log a\big)-L\big(E(y-m-2)\big).
    \end{align*}
    Expanding $L\big(E(y-m-2)+\log a\big)$ using the Taylor series for $L$ at $E(y-m-2)$ gives
    \begin{align*}
        L(E(y-m-2)+\log a)-L(E(y-m-2)) &= \sum_{j=1}^{\infty} \frac{L^{(j)}(E(y-m-2)) (\log a)^j}{j!} \\
            &= \frac{\log a}{E'(y-m-2)} + \frac{(\log a)^2 E''(y-m-2)}{2E'(y-m-2)^3} + \cdots
    \end{align*}
    which is a valid sum in the structure $D_{X'}(C_{\bar{\alpha},i})$ for any $X' \ni y$ such that $X'$ and $C_{\bar{\alpha},i}$ satisfy the ordering and separation assumptions of Remark \ref{order props}.
    
    To prove the lemma, we split into two cases based on whether $x-y$ is infinitesimal with respect to $H_i$.
    First suppose $x-y$ is not infinitesimal, so $\frac{1}{x-y} \in \mathrm{Fin}(H_i)$.
    Then Condition (1) is trivially satisfied.
    Since there is some $n \in \mathbb{N}$ such that $\frac{1}{x-y} < n$, and since $\frac{E(y-m)}{\log a}$ is infinite, we have $\frac{1}{x-y} < \frac{E(y-m)}{\log a}$.
    By the computations above, this shows that $E(x-m) > E(y-m)^a$.
    
    Now suppose $x-y$ is infinitesimal with respect to $H_i$, so $\frac{1}{x-y} \not \in \mathrm{Fin}(H_i)$.
    Fix $z \in \lambda\left(\frac{1}{x-y}\right)$, and $n \in \mathbb{N}$ such that 
        $$E(z-n) < \frac{1}{x-y} < E(z+n).$$
    If the second condition of the lemma fails, then $E(x-m) \le E(y-m)^a$ for some $m \in \mathbb{N}$ and $a \in (C_{\bar{\alpha},i})_{>1}$. 
    So 
    \begin{align*}
        E(z + n) > \frac{1}{x-y} > \frac{E'(y-m-2)}{2\log a} > E(y-m-2).
    \end{align*}
    So $z + n > y-m-2$, and thus $\lambda\left(\frac{1}{x-y}\right) \ge \beta$.
    This contradicts Condition (1) of Definition \ref{define X_alpha} because we assumed that both $x$ and $y$ are in $X_{\beta}$.
    Conversely, suppose $X_{\beta}$ fails Condition (1) due to $x$ and $y$, i.e., $\lambda\left(\frac{1}{x-y}\right) \ge \beta$.
    Then there is some $m \in \mathbb{N}$ such that $z - n > y-m$.
    So
    \begin{align*}
        \frac{1}{x-y} > E(z-n) > E(y-m) > \frac{E'(y-m-2)}{2\log a}.
    \end{align*}
    By the computations above, this means $E(x-m) \le E(y-m)^a$.
\end{proof}

So the ordering and separation assumptions in Remark \ref{order props} are satisfied by $X_{\beta}$ and $C_{\bar{\alpha},i}$.
So we can define 
\begin{align*}
    C_{\bar{\alpha}\frown \beta,i} := D_{X_{\beta}}(C_{\bar{\alpha},i}).
\end{align*}

\subsubsection{Defining $H_{i+1}$ and $\iota_i : H_i \to H_{i+1}$}

We would like to define $H_{i+1}$ as the direct limit of the sets $C_{\bar{\alpha},i}$, so we must show that these sets form a directed system ordered by inclusion.

\begin{lemma}\label{C's form directed system}
    Let $\bar{\alpha} = \{\alpha_0 < \cdots < \alpha_k\} \subset (H_i/\mathrm{Fin}(H_i))_{>F}$, and let $\beta \in (H_i/\mathrm{Fin}(H_i))_{>F}$ such that $\beta \ne \alpha_{j}$ for $j=1,\dots,k$.
    Then $C_{\bar{\alpha},i}$ is an $\mathcal{L}_{\an,\exp}(\log)$-substructure of $C_{\bar{\alpha}\cup \beta,i}$.
\end{lemma}

\begin{proof}
    If $\beta > \alpha_k$, then $C_{\bar{\alpha},i} \subset D_{C_{\bar{\alpha},i}}(X_{\beta}) = C_{\bar{\alpha}\frown \beta,i} = C_{\bar{\alpha}\cup \beta,i}$.
    If $\beta < \alpha_0$, then
    \begin{align*}
        C_{\varnothing,i} \subset C_{\{\beta\},i} 
            & \Rightarrow C_{\{\alpha_0\},i} \subset C_{\{\beta < \alpha_0\},i} \\
            & \Rightarrow C_{\{\alpha_0 < \alpha_1\},i} \subset C_{\{\beta < \alpha_0 < \alpha_1\},i} \\
            &\vdots \\
            & \Rightarrow C_{\{\alpha_0 < \cdots < \alpha_k\},i} \subset C_{\{\beta < \alpha_0 < \cdots < \alpha_k\},i}.
    \end{align*}
    Similarly, if $\alpha_{l} < \beta < \alpha_{l+1}$, then
    \begin{align*}
        C_{\{\alpha_0 < \cdots < \alpha_{l}\},i} \subset C_{\{\alpha_0 < \cdots < \alpha_{l} < \beta\},i} 
            & \Rightarrow C_{\{\alpha_0 < \cdots < \alpha_{l+1}\},i} \subset C_{\{\alpha_0 < \cdots < \alpha_{l} < \beta < \alpha_{l+1}\},i} \\
            &\vdots \\
            & \Rightarrow C_{\{\alpha_0 < \cdots < \alpha_k\},i} \subset C_{\{\alpha_0 < \cdots < \alpha_{l} < \beta < \alpha_{l+1} < \cdots < \alpha_k\},i}.
    \end{align*}
    So we can view every element of $C_{\bar{\alpha},i}$ as an element of $C_{\bar{\alpha}\cup\beta,i}$ that does not mention any elements of $X_{\beta}$.
    
    Let $s \in C_{\bar{\alpha},i}$.
    Since the definition of $\exp s$ mentions the same (finitely many) elements of each $X_{\alpha_0}, \dots,X_{\alpha_k}$ that $s$ does, the definition of $\exp s$ in $C_{\bar{\alpha},i}$ is the same as its definition in $C_{\bar{\alpha}\cup \beta,i}$.
    The same holds for the other function symbols of $\mathcal{L}_{\an,\exp}(\log)$, so  $C_{\bar{\alpha},i}$ is an $\mathcal{L}_{\an,\exp}(\log)$-substructure of $C_{\bar{\alpha}\cup \beta,i}$.
\end{proof}

\begin{cor}
    $\left(C_{\bar{\alpha},i} : \bar{\alpha} \in (H_i/\big(\mathrm{Fin}(H_i))_{>F}\big)^{<\omega}\right)$ forms a directed system.
\end{cor}

Build $H_{i+1}$ from $H_i$ as the direct limit of the directed system 
    $$\left(C_{\bar{\alpha},i} : \bar{\alpha} \in \big((H_i/\mathrm{Fin}(H_i))_{>F}\big)^{<\omega}\right).$$
$H_{i+1}$ can be made into a model of $T_{\mathrm{an}}(\exp,\log)$ because each $C_{\bar{\alpha},i} \vDash T_{\mathrm{an}}(\exp,\log)$.

\begin{lemma}\label{lambda exists}
    For all $s \in (H_{i+1})_{+\infty}$, $\lambda(s)$ exists.
\end{lemma}

\begin{proof}
    Let $s \in C_{\bar{\alpha},i}$ with $s > F$.
    We divide into several cases:
    \begin{enumerate}
        \item If $s$ is bounded in $D_{X_i'}(F)$, then 
            $$\lambda(s) := \left\{E\big(L^{\circ(i+2)}(\tau)\big)\right\}$$
        works, by Lemma \ref{log-exp bounded}.
        
        \item If $s > D_{X_i'}(f)$ and is bounded in $D_{X_i}\big(D_{X_i'}(F)\big)$, then 
            $$\lambda(s) := \beta$$ 
        for $E\big(L^{\circ(i+1)}(\tau)\big) \in X_{\beta} \subset \beta \subset H_i$ works, again by Lemma \ref{log-exp bounded}.
        
        \item If $s > C_{\varnothing,i}$, then we may write $\bar{\alpha} = \{\alpha_0 < \cdots < \alpha_k\}$ for some $k \ge 0$. 
        If $s$ is bounded in $C_{\{\alpha_0\},i}$, then 
            $$\lambda(s) := \alpha_0$$
        works, by Lemma \ref{log-exp bounded}.
        
        \item Finally, if $s > C_{\{\alpha_0 < \cdots < \alpha_l\},i}$ and is bounded in $C_{\{\alpha_0 < \cdots < \alpha_{l+1}\},i}$ for some $l < k$, then use
            $$\lambda(s) := \alpha_{l+1}$$
        by Lemma \ref{log-exp bounded}.
    \end{enumerate}
    So $\lambda(s)$ is defined for all $s > F$ in $H_{i+1}$.
\end{proof}

We will now define an order preserving embedding $\iota_i : H_i \to H_{i+1}$.
If $i=0$, we think of $\iota_0 : H_0 \to H_1$ as substituting $E(L(\tau))$ for $\tau$. 
Let $s \in H_0 = F((\tau^{-1}))^{le}$.
We follow the notation of \cite{LEseries}, in which $F((\tau^{-1}))^e$ is the direct limit of $(K_n : n \in \mathbb{N})$, and $F((\tau^{-1}))^{le}$ is the direct limit of $(L_n : n \in \mathbb{N})$.
We may identify $s$ with a unique element of $C_{\varnothing,0} \subset H_1$ as follows:
\begin{enumerate}
    \item If $s = \sum a_r \tau^r \in K_0$, then define $\iota_0(s) = \sum a_r E(L(\tau))^r$
    
    \item If $s = \sum f_a e(a) \in K_{n+1}$, then define $\iota_0(s) = \sum \iota_0(f_a) e(\iota_0(a))$
    
    \item For any $s \in L_n$, let $\hat{s} = \eta_n(s) \in F((\tau^{-1}))^e$, and express 
        $$s = \hat{s}\big(\log^{\circ n}(\tau)\big).$$
    \begin{enumerate}
        \item If $\hat{s} = \sum a_r \tau^r \in K_0$, then define 
            $$\iota_{0,n}(\hat{s}) = \sum a_r \log^{\circ n}E(L(\tau))^r.$$
    
        \item If $\hat{s} = \sum f_a e(a) \in K_{m+1}$, then define 
            $$\iota_{0,n}(\hat{s}) = \sum \iota_{0,n}(f_a) e(\iota_{0,n}(a)).$$
    \end{enumerate}
    Define $\iota_0(s) = \iota_{0,n}(\hat{s})$.
\end{enumerate}

Now suppose $i > 0$.
We inductively define $\iota_i$ on the generators built from $X_{i-1}'$, $X_{i-1}$, and $X_{\alpha}$ for $\alpha \in \big(H_{i-1}/\mathrm{Fin}(H_{i-1})\big)_{>F}$.
For all $d \in \mathbb{N}$ and $a \in F$, define
\begin{align*}
    \iota_i\big(E^{(d)}(L^{\circ (i+1)}(\tau))^a\big) = E^{(d)}\big(L^{\circ (i+1)}(\tau)\big)^a \in C_{\varnothing,i}.
\end{align*}
Since the image of each generator is a single generator, we can extend $\iota_i$ so that $\iota_i\big(D_{X_{i-1}'}(F)\big) = D_{X_i}(F) \subset C_{\varnothing,i}$.

For all $d \in \mathbb{N}$ and $a \in D_{X_{i-1}'}(F)$, define
\begin{align*}
    \iota_i\big(E^{(d)}(L^{\circ i}(\tau))^a\big) = E^{(d)}\big(E(L^{\circ (i+1)}(\tau))\big)^{\iota_i(a)} \in C_{\{\beta\},i}
\end{align*}
where $\beta$ is such that $E(L^{\circ (i+1)}(\tau)) \in X_{\beta}$.
Again, since the image of each generator is a single generator, we can extend $\iota_i$ to identify $C_{\varnothing,i-1}$ with an isomorphic copy of itself in $C_{\{\beta\},i}$. 

Now assume we have defined $\iota_i$ on $C_{\bar{\alpha},i}$ for some $\bar{\alpha} \subset \big(H_{i-1}/\mathrm{Fin}(H_{i-1})\big)_{>F}$.
Let $\beta > \bar{\alpha}$.
For all $d \in \mathbb{N}$ and $a \in C_{\bar{\alpha},i}$, define
\begin{align*}
    \iota_i\big(E^{(d)}(x)^a\big) = E^{(d)}\big(\iota_{i-1}(x)\big)^{\iota_i(a)}.
\end{align*}
Again, we can extend $\iota_i$ to identify $D_{X_{\beta}}(C_{\bar{\alpha},i-1})$ with an isomorphic copy of itself in $C_{\iota_{i-1}(\bar{\alpha} \cup \{\beta\}),i}$.

Thus, we have identified $H_i$ with an isomorphic copy of itself $\iota_i(H_i) \subset H_{i+1}$.
We will often suppress the $\iota_0, \iota_1,\dots$ maps and identify each $H_i$ with its isomorphic copy in $H_{i+n}$.

\subsection{Building a structure closed under the symbols of \texorpdfstring{$\mathcal{L}_{\mathrm{transexp}}$}{Ltransexp}}\label{subsection 5.2}

Let $M_F$ be the direct limit of the directed system given by $(H_i, \iota_i)_{i \in \NN}$.
Next we will show that $M_F$ is closed under $E$, its derivatives, $L$, and the inverses of the derivatives of $E$. 

\begin{lemma}\label{E closure}
    For all $s \in (H_i)_{>F}$ and $d \in \mathbb{N}$, $E^{(d)}(s) \in H_{i+1}$.
\end{lemma}

\begin{proof}
    Let $s \in \beta \in (H_i/\mathrm{Fin} (H_i))_{>F}$.
    If $s-l \in X_{\beta}$ for some $l \in \mathbb{Z}$, then $E^{(d)}(s) \in C_{\{\beta\},i}$ by construction.
    
    If $s-l \not\in X_{\beta}$ for any $l \in \mathbb{Z}$, then we must have $i > 0$, and each $s-l$ must have been excluded from $X_{\beta}$ for a reason.
    Let $m \in \mathbb{N}$ and $x \in X_{\beta}$ be the unique elements such that $s-m$ and $x$ violate Condition (1) and satisfy Condition (2) of Definition \ref{define X_alpha}.
    For ease of notation we may assume $m = 0$.
    Let $\beta > \alpha \in (H_i/\mathrm{Fin}(H_i))_{>F}$.
    By Lemma \ref{condition 1}, either $s > x$ and there is some $1 < a \in C_{\{\alpha\},i}$ and $n \in \mathbb{N}$ such that
        $$E(s-n) \le E(x-n)^a$$
    or $x > s$ and there are $a$ and $n$ as above with
        $$E(x-n) \le E(s-n)^a.$$
    Following the first computation in Lemma \ref{condition 1}, this means that
    \begin{align*}
        |s-x| < \frac{1}{E(x-n-3)}.
    \end{align*}
    Since $s-x \in H_i$, we have $\iota_i(s-x) \in C_{\bar{\alpha},i} \subset H_{i+1}$ for some $\bar{\alpha} \in H_i/\mathrm{Fin} (H_i)$.
    So we can represent $E^{(l)}(s-n-3)$ by
    \begin{align*}
        \sum_{k=0}^{\infty} \frac{E^{(l+k)}(x-n-3)}{k!}(s-x)^k \in C_{\bar{\alpha} \cup \{\beta\},i}
    \end{align*}
    for $l \in \mathbb{N}$.
    Then for $m=0,\dots,n-4$, we can represent $E(s-m)$ by 
    \begin{align*}
        \exp^{\circ(n-3-m)}\left(\sum_{k=0}^{\infty} \frac{E^{(k)}(x-n-3)}{k!}(s-x)^k\right).
    \end{align*}
    Finally, we express $E^{(d)}(s)$ in terms of the expressions we have found for $E^{(l)}(s-n-3)$ and $E(s-m)$, $l,m \in \mathbb{N}$, using the Bell polynomial difference-differential equations for $E$.
\end{proof}

We will use the following technical lemma to show that for all positive infinite (relative to $F$) elements $s \in H_i$, we can identify $L(s)$ with an element of $H_{i+1}$.

\begin{lemma}\label{dominant leading monomial is E(x-n)}
    Suppose $X$ and $\boldsymbol{k}$ satisfy the ordering and separation axioms of Remark \ref{order props}.
    Then for any positive infinite (relative to $\boldsymbol{k}$) element $s \in D_X(\boldsymbol{k})$, there are $n_0,n_1,n_2 \in \mathbb{N}$ for which we can express $\log^{\circ n_0}(s)$ in the form 
        $$\log^{\circ n_0}(s) = \exp^{\circ n_1}E(x-n_2) + s_0$$
    with $\mathrm{Supp} \left(\frac{s_0}{\exp^{\circ n_1}E(x-n_2)}\right) < \Gamma_{X-n_2,0}^{small}$.
\end{lemma}

\begin{proof}
    By Lemma \ref{dominant leading monomial}, there is some $m \in \mathbb{N}$ for which $s$ can be identified with an element $cM_0(1 + \epsilon) \in K_{X-m}$ with $c>0$, $\mathrm{Supp} (\epsilon) < \Gamma_{X-m,0}^{small}$, and $M_0$ of the form
        $$e_{X-m}(\alpha)\prod_{j=1}^pE'(x_j-m)^{a_{j}}$$
    with $\alpha \in A_{X-m,k}$, $a_1,\dots,a_p \in \kk$, and $x_1 > \cdots > x_p \in X$.
    We will also call this element $s$.
    Then 
    \begin{align*}
        \log s &= \log \big(cM_0(1+\epsilon)\big) \\
            &= \alpha + \sum_{j=1}^p a_j \log E'(x_j-m) + \log c + \sum_{l=1}^{\infty} \frac{(-1)^{l+1}}{l}\epsilon^l \\
            &\in K_{X-m}.
    \end{align*}
    First suppose $\alpha = 0$.
    Then we must have $a_1>0$ since $s > F$.
    We can compute that 
    \begin{align*}
        \log\varphi_m(\log s)
            &= E(x_1-m-2) + \log a_1 + 
            \log\varphi_m\left(\frac{\log s}{a_1E(x_1-m-1)}\right) 
    \end{align*}
    So $\log(\log s) \in K_{X-m-2}$, and $\mathrm{Supp} \left(\log\varphi_m\left(\frac{\log s}{a_1E(x_1-m-2)}\right)\right) < \Gamma_{X-m-2,0}^{small}$.
    Let $\epsilon = \log\varphi_m\left(\frac{\log s}{a_1E(x_1-m-1)}\right)$.
    Then
    \begin{align*}
        \log\varphi_{m+1}(\log\varphi_m(\log s)) &= E(x_1-m-3) + \log\left(1+\varphi_m\left(\frac{\log a_1 + \epsilon}{E(x_1-m-2)}\right)\right).
    \end{align*}
    So the lemma holds with $n_0 = 3$, $n_1 = 0$, and $n_2 = m+3$ in the case $\alpha = 0$.
    
    Now suppose $0 \ne \alpha \in A_{X-m,0}$.
    Then we can write
        $$\alpha = s_1 \prod_{j=1}^qE(y_j-m)^{b_j}$$
    with $b_1 > 0$ and $\mathrm{Lv}(s_1) = \Gamma_{X-m,0}^{small}$.
    By Lemma \ref{log closure}, there is some $m'$ such that the image of $s_1$ under $\varphi_{m+m'-1} \circ \cdots \circ \varphi_m$ has a logarithm.
    Then we can compute that
    \begin{align*}
        \log(\varphi_{m+m'-1} \circ \cdots \circ \varphi_m)(\log(s)) = b_1(e_{X-m-m'})^{\circ (m'-1)}\big(E(y_1-m-m')\big) +\epsilon
    \end{align*}
    with $\epsilon = \log\left(\frac{(\varphi_{m+m'-1} \circ \cdots \circ \varphi_m)(\log(s))}{b_1(e_{X-m-m'})^{\circ (m'-1)}\big(E(y_1-m-m')\big)}\right)$ and $\mathrm{Supp} (\epsilon) < \Gamma_{X-m-m',0}^{small}$.
    Taking two more logarithms gives us the conclusion of the lemma with $n_0 = 4$, $n_1 = m'-1$, and $n_2 = m+m'+2$.
    
    If $0 \ne \alpha \in A_{X-m,n+1}$, then let $\alpha_{n+1} = \alpha$.
    There are $c_l,\alpha_l, \epsilon_l$ for $l = 0,\dots,n$ such that
    \begin{align*}
        \alpha_{l+1} = c_l e_{X-m}(\alpha_l) (1 + \epsilon_l)
    \end{align*}
    with $c_l \in K_{X-m,l}$, $\alpha_l \in A_{X-m,l}$, and $\epsilon_l \in \mathfrak{m}\big(B_{X-m,l+1}\big)$.
    Write 
        $$\alpha_0 = s_1 \prod_{j=1}^qE(y_j-m)^{b_j}$$
    with $b_1 > 0$ and $\mathrm{Lv}(s_1) = \Gamma_{X-m,0}^{small}$.
    Again by Lemma \ref{log closure}, there is some $m'$ such that 
    the image of $s_1, c_0,\dots,c_n$ under $\varphi_{m+m'-1} \circ \cdots \circ \varphi_m$ have as many logarithms as needed.
    Then
    \begin{align*}
        \log^{\circ (1 + (n+1))}(\varphi_{m+m'-1} \circ \cdots \circ \varphi_m)(\log s) = b_1(e_{X-m-m'})^{\circ (m'-1)}\big(E(y_1-m-m')\big) + \epsilon
    \end{align*}
    where $\epsilon = \log\left(\frac{(\varphi_{m+m'-1} \circ \cdots \circ \varphi_m)(\log s)}{b_1(e_{X-m-m'})^{\circ (m'-1)}\big(E(y_1-m-m')\big)}\right)$ with $\Supp(\epsilon) < \Gamma_{X-m-m',0}^{small}$.
    Taking two more logarithms completes the lemma in this case with $n_0 = n+5$, $n_1 = m'-1$, and $n_2 = m+m'+2$.
\end{proof}

\begin{lemma}\label{L closure}
    For all $s \in (H_i)_{>F}$, we have $L(s) \in H_{i+1}$.
\end{lemma} 

\begin{proof}
    First suppose $i = 0$, so that $s \in F((\tau^{-1}))^{le}$.
    Then there are some $n,k \in \mathbb{N}$ such that 
        $$\mathrm{Lm}\left(\log^{\circ n}(s)\right) = \log^{\circ k}(\tau).$$
    So $\log^{\circ n}(s) = r\log^{\circ k}(\tau)(1 + \epsilon)$ with $r>0$.
    Then we can represent $L(s)$ as follows:
    \begin{align*}
        L(s) = L&\big(\log^{\circ(n+1)}(s)) + n+1 \\
            = L&\big(\log^{\circ(k+1)}(\tau) + \log r + \log(1 + \epsilon)\big) + n + 1 \\
            = L&\big(\log^{\circ(k+1)}(\tau)\big) + \sum_{l=1}^{\infty} \frac{L^{(l)}\left(\log^{\circ(k+1)}(\tau)\right)}{l!}\big(\log r + \log(1 + \epsilon)\big)^l  + n + 1.
    \end{align*}
    Now $L\left(\log^{\circ(k+1)}(\tau)\right) = L(\tau) -k-1$, and differentiating, we get $L'\left(\log^{\circ(k+1)}(\tau)\right) \cdot \left(\log^{\circ(k+1)}(\tau)\right)' = L'(\tau)$.
    So 
    \begin{align*}
        L'\left(\log^{\circ(k+1)}(\tau)\right) = \frac{L'(\tau)}{\left(\log^{\circ(k+1)}(\tau)\right)'} = L'(\tau)\log^{\circ k}(\tau) \cdots \log \tau \cdot \tau.
    \end{align*}
    We can continue differentiating to obtain expressions for $L^{(l)}\left(\log^{\circ(k+1)}(\tau)\right)$, $l \in \NN$, in terms of $L'(\tau),\dots,L^{(l)}(\tau)$ and $\log^{\circ k}(\tau), \dots, \log \tau,\tau$, all of which are elements of $C_{\varnothing,i}$.
    This allows us to finish expressing $L(s)$ as
    \begin{align*}
        L(s) =& n-k + L(\tau) + \big(L'(\tau)\log^{\circ k}(\tau) \cdots \log \tau \cdot \tau\big)\big(\log r + \log(1 + \epsilon)\big) \\
            &+\Big(L''(\tau)\log^{\circ k}(\tau) \cdots \log \tau \cdot \tau + L'(\tau)\big(\log^{\circ k}(\tau) \cdots \log \tau \cdot \tau\big)'\Big)\big(\log r + \log(1 + \epsilon)\big)^2 \\
            &+ \cdots.
    \end{align*}
    Since the $l$th derivative of $\log^{\circ k}(\tau) \cdots \log \tau \cdot \tau$ is infinitesimal with $\tau^{l-1}$ in the denominator for $l > 1$, this sum is an element of $C_{\varnothing,0}$.
    So $L(s) \in H_{1}$.
    
    Now suppose $i > 0$, and let $s \in (H_i)_{>F}$.
    Then we can identify $s$ with an element, which we also call $s$, of $C_{\bar{\alpha},i-1}$ for some $\bar{\alpha} \subset H_{i-1}/\mathrm{Fin} (H_{i-1})$.
    By Lemma \ref{dominant leading monomial is E(x-n)}, there are $n_0, n_1,n_2 \in \mathbb{N}$ such that we can express
        $$\log^{\circ n_0}(s) = \exp^{\circ n_1}E(x-n_2) + s_0$$
    with $\mathrm{Supp} \left(\frac{s_0}{\exp^{\circ n_1}E(x-n_2)}\right) < \Gamma_{X-n_2,0}^{small}$ and $x \in X$, where $X$ is either $X_i'$, $X_i$, or $X_{\alpha}$ for some $\alpha \in H_{i-1}/\mathrm{Fin} (H_{i-1})$.
    So identifying $\exp^{\circ n_1}E(x-n_2)$ with $E(x+n_1-n_2)$, we can represent $L(s)$ by
    \begin{align*}
        L(s) &= L\big(\log^{\circ n_0}(s)\big) + n_0 \\
            &= L\big(E(x+n_1-n_2) + s_0\big) + n_0 \\
            &= x+ n_0+n_1-n_2  + \sum_{l=1}^{\infty} \frac{L^{(l)}\big(E(x+n_1-n_2)\big)}{l!}(s_0)^l \\
            &= x+n_0+n_1-n_2 + \frac{s_0}{E'(x+n_1-n_2)} + \frac{(s_0)^2E''(x+n_1-n_2)}{E'(x+n_1-n_2)^3} + \cdots \\
            &\in K_{X-n_2}.
    \end{align*}
    So we have identified $L(s)$ with an element of $H_{i+1}$.
\end{proof}

We will use the next two lemmas to show that for all positive infinite (relative to $F$) elements $s \in H_i$, we can represent $\big(E^{(d)}\big)^{-1}(s)$ by an element of $H_{i+3}$.

\begin{lemma}\label{epsilon calculation}
    For any infinite $x \in H_i$ and $1 <a \in H_{i+1}$ such that $\frac{\log a}{E(x-1)}$ is infinitesimal, we can compute $\epsilon_{a,d}(x) \in H_{i+1}$ such that 
        $$E(x+\epsilon_{a,d}(x))=aE^{(d)}(x)$$
    in $H_{i+2}$.
\end{lemma}

\begin{proof}
    Note that $E(x+\epsilon_{a,d}(x))=aE^{(d)}(x)$ if and only if (taking log twice)
        $$E(x-2+\epsilon_{a,d}(x))=E(x-2)+\log\left(1+\frac{\log B_d(x-1)+\log a}{E(x-1)}\right).$$
    Let $C(x) = \log\left(1+\frac{\log B_d(x-1)+\log a}{E(x-1)}\right)$, an infinitesimal.
    Now
    \begin{align*}
       E\big(x-2+\epsilon_{a,d}(x)\big)&=E(x-2)+C(x) 
            &&\text{ if and only if } \\ E\big(L(x)+\epsilon_{a,d}(L(x)+2)\big)&=x+ C\big(L(x)+2\big) 
            &&\text{ if and only if } \\ L(x)+\epsilon_{a,d}\big(L(x)+2\big)&=L\big(x+ C(L(x)+2)\big).
    \end{align*}
    Rearranging terms, we get
    $$\epsilon_{a,d}\big(L(x)+2\big)=L\big(x+ C(L(x)+2)\big)-L(x) = \sum_{n=1}^{\infty} \frac{L^{(n)}(x)}{n!}C\big(L(x)+2\big)^n.$$
    And substituting $E(x-2)$ for $x$, we get 
    \begin{align*}
        \epsilon_{a,d}(x) 
            &= \frac{C(x)}{E'(x-2)}-\frac{ E''(x-2)C(x)^2}{2!E'(x-2)^3}+\cdots 
    \end{align*}
    
    We must show that $\epsilon_{a,d}(x) \in H_{i+1}$.
    Let $\bar{\beta} =\{\beta_0 < \cdots < \beta_j\}$ be the shortest sequence such that $E(x-1) \in C_{\bar{\beta},i}$, and let $\bar{\alpha} = \{\alpha_0 < \cdots < \alpha_k\}$ be the shortest sequence such that $a \in C_{\bar{\alpha},i}$.
    If $\beta_j > \alpha_k$, then we view $\epsilon_{a,d}(x)$ as a sum of monomials built from $X_{\beta_j}$ with coefficients in $C_{\{\beta_0 < \cdots < \beta_{j-1}\} \cup \bar{\alpha},i}$.
    If $\beta_j < \alpha_k$, then we view $\epsilon_{a,d}(x)$ as a sum of monomials built from $X_{\alpha_k}$ with coefficients in $C_{\bar{\beta} \cup \{\alpha_0 < \cdots < \alpha_{k-1}\},i}$.
    If $\beta_j = \alpha_k$, then we view $\epsilon_{a,d}(x)$ as a sum of monomials built from $X_{\alpha_k}$ with coefficients in $C_{\{\beta_0 < \cdots < \beta_{j-1}\} \cup \{\alpha_0 < \cdots < \alpha_{k-1}\},i}$.
    In any case, $\epsilon_{a,d}(x) \in H_{i+1}$.
\end{proof}

\begin{lemma}
    For all $s \in (H_i)_{>F}$ and $d \in \mathbb{N}$, $\big(E^{(d)}\big)^{-1}(s) \in H_{i+3}$.
\end{lemma}

\begin{proof}
    We will find $f \in H_{i+2}$ such that 
        $$0 \le \frac{s}{E^{(d)}(f)} -1 \le \frac{1}{E(f-1)^{1/3}}.$$
    Let $\mu = \frac{s}{E^{(d)}(f)} -1$.
    This suffices to prove the Lemma because then 
    \begin{align*}
        \left(E^{(d)}\right)^{-1}(s)&=\big(E^{(d)}\big)^{-1}\left(E^{(d)}(f)+ \mu E^{(d)}(f)\right) 
            \\
            &= \sum_{n=0}^{\infty}\frac{\left(\big(E^{(d)}\big)^{-1}\right)^{(n)}\big(E^{(d)}(f)\big)}{n!} \left(\mu E^{(d)}(f)\right)^n 
            \\
            &=f +\frac{\mu E^{(d)}(f)}{E^{(d+1)}(f)}-\frac{\mu^2E^{(d)}(f)^2E^{(d+2)}(f)}{2!E^{(d+1)}(f)^3}+\dots
            \\
            &= f  +\frac{\mu B_d(f-1)}{B_{d+1}(f-1)}-\frac{\mu^2B_d(f-1)^2B_{d+2}(f-1)}{2!B_{d+1}(f-1)^3}+\dots
            \\
            &= f +\mu\frac{B_d(f-1)}{E'(f-1)^{d+1}}\sum_{n=0}^{\infty}\left(1-\frac{B_{d+1}(f-1)}{E'(f-1)^{d+1}}\right)^n 
            \\
            &- \mu^2\frac{B_d(f-1)^2B_{d+2}(f-1)}{2!E'(f-1)^{3d+3}}\sum_{n=0}^{\infty}\left(1-\frac{B_{d+1}(f-1)^3}{E'(f-1)^{3d+3}}\right)^n.
    \end{align*}
    We will show the final sum is an element of $H_{i+3}$ if $f \in H_{i+2}$.
    We must show it is a valid sum in some $C_{\bar{\alpha},i+2}$ containing both $s$ and $E^{(d)}(f)$.
    Let $\bar{\alpha} = \{\alpha_0 < \cdots < \alpha_k\}$ and $m \in \mathbb{N}$ be such that $s$ and $E^{(d)}(f)$ are in the stage-$K_{X_{\alpha_k}-m}$ construction of $C_{\bar{\alpha},i+2}$.
    
    The above expansion of each $\Big(\big(E^{(d)}\big)^{-1}\Big)^{(n)}\big(E^{(d)}(f)\big)\cdot E^{(d)}(f)^n$ for $n >0$ is a valid infinite sum involving only integer powers of $E^{(k)}(f-1)$ for $k \in \mathbb{N}$.
    The sum of exponents in each term is a negative integer, and the largest term of each sum is $\frac{1}{E'(f-1)}$.
    Also, for each $l$, there are only finitely many terms with sum of exponents equal to $l$ that can appear in the expansion of $\Big(\big(E^{(d)}\big)^{-1}\Big)^{(n)}\big(E^{(d)}(f)\big)\cdot E^{(d)}(f)^n$ for some $n \in \mathbb{N}$.
    So if $0 \le \mu \le \frac{1}{E(f-1)^{1/3}}$, then the final sum representing $\big(E^{(d)}\big)^{-1}(s)$ is summable:
    \begin{enumerate}
        \item For each monomial $M$ appearing in the expression, there are only finitely many $n \in \mathbb{N}$ such that 
            $$N \in \mathrm{Supp} \left(\mu^n \Big(\big(E^{(d)}\big)^{-1}\Big)^{(n)}\big(E^{(d)}(f)\big)\cdot E^{(d)}(f)^n\right)$$
        for some $N \in M \Gamma_{X_{\alpha_k}-m,0}^{small}$.
        
        \item $\displaystyle \bigcup_{n > 1} \mathrm{Supp} \left(\mu^n \Big(\big(E^{(d)}\big)^{-1}\Big)^{(n)}\big(E^{(d)}(f)\big)\cdot E^{(d)}(f)^n\right)$ is reverse well-ordered.
    \end{enumerate}
    So we need only find an element $f$ satisfying the required inequalities.
    
    Define $\displaystyle f := L(s)-\frac{d\log E'(L(s)-1)+\frac{1}{E(L(s)-1)^{1/2}}}{E'\big(L(s)-2\big)E\big(L(s)-1\big)}$.
    Note that the exponent of $\frac{1}{E(L(s)-1)^{1/2}}$ in the numerator of $f$ is $1/2$, while the exponent of the bound $\frac{1}{E(f-1)^{1/3}}$ on $\mu$ is $1/3$.
    Instead of 1/2 and 1/3, any two numbers between 0 and 1 that are not infinitesimally close to 0, 1, or each other would work just as well.
    
    We will show that 
        $$E^{(d)}(f) \le s \le \left(1+\frac{1}{E(L(s)-1)^{1/3}}\right)E^{(d)}(f).$$
    This suffices since $\frac{1}{E(L(s)-1)^{1/3}}\le \frac{1}{E(f-1)^{1/3}}$.
    We can compute $\epsilon_{a,d}(f)$ such that $aE^{(d)}(f) = E\big(f + \epsilon_{a,d}(f)\big)$ for $a = 1$ and $a = 1 + \frac{1}{E(L(s)-1)^{1/3}}$ by Lemma \ref{epsilon calculation}.
    So it suffices to show that
    \begin{align*}
        f + \epsilon_{1,d}(f) &\le L(s) \le f + \epsilon_{1+\frac{1}{E(L(s)-1)^{1/3}},d}(f) \\
        \epsilon_{1,d}(f) &\le L(s) - f \le \epsilon_{1+\frac{1}{E(L(s)-1)^{1/3}},d}(f) \\
        \epsilon_{1,d}(f) &\le \frac{d\log E'(L(s)-1)+\frac{1}{E(L(s)-1)^{1/2}}}{E'\big(L(s)-2\big)E\big(L(s)-1\big)} \le \epsilon_{1+\frac{1}{E(L(s)-1)^{1/3}},d}(f).
    \end{align*}
    Let
    \begin{align*}
        \delta :&= L(s) - f 
            \\
            &= \frac{d\log E'(L(s)-1)+\frac{1}{E(L(s)-1)^{1/2}}}{E'\big(L(s)-2\big)E\big(L(s)-1\big)}
            \\
            &= \frac{dE(L(s)-2) + d\log E'(L(s)-2) + \frac{1}{E(L(s)-1)^{1/2}}}{E'\big(L(s)-2\big)E\big(L(s)-1\big)}.
    \end{align*}
    So the dominant monomial of $\delta^n$ is $\frac{d^n}{E'(L(s)-3)^nE(L(s)-1)^n}$.
    
    We can compute $\epsilon_{a,d}(f)$ using Lemma \ref{epsilon calculation}.
    First we compute $E^{(d)}(f-m)$ for $m \ge 1$:
    \begin{align*}
        E^{(d)}(f-m) 
            &= E^{(d)}\left(L(s)-\delta -m\right) 
            \\
            &= \sum_{n=0}^{\infty}\frac{E^{(d+n)}(L(s)-m)}{n!} (-\delta)^n .
    \end{align*}
    Now compute $C(f)$:
    \begin{align*}
        C(f) &= \log \left(1 + \frac{\log B_d(f-1) + \log a}{E(f-1)}\right) \\
            &= \sum_{k=1}^{\infty} \frac{(-1)^{k+1}}{k}\left(\frac{\log B_d(f-1) + \log a}{E(f-1)}\right)^k \\
            &= \sum_{k=1}^{\infty}\frac{(-1)^{k+1}}{k}\left(\frac{dE(x-2) + d\log E'(f-2) + \log a + \log\left(\frac{B_d(f-1)}{E'(f-1)^d}-1\right)}{E(f-1)}\right)^k 
            \\
            &= \frac{dE(f-2) + d\log E'(f-2) + \log a}{E(f-1)} + \cdots 
            \\
            &= \frac{d\sum_{n=0}^{\infty}\frac{E^{(n)}(L(s)-2)}{n!} \delta^n  + d\log \sum_{n=0}^{\infty}\frac{E^{(1+n)}(L(s)-2)}{n!} \delta^n  + \log a}{\sum_{n=0}^{\infty}\frac{E^{(n)}(L(s)-1)}{n!} \delta^n } + \cdots 
            \\
            &= \frac{d E(L(s)-2) + d\log E'(L(s)-2) + \log a}{E(L(s)-1)} +\cdots
    \end{align*}
   where the ``$\cdots$" has terms with $E(L(s)-1)^n$ in the denominator for $n \ge 2$.
    
    Finally, we compute $\epsilon_{a,d}(f)$:
    \begin{align*}
        \epsilon_{a,d}(f) 
            =& \frac{C(f)}{E'(f-2)} - \frac{C(f)^2E''(f-2)}{2!E'(f-2)^3} + \cdots 
            \\
            =& \frac{\frac{d E(L(s)-2) + d\log E'(L(s)-2) + \log a}{E(L(s)-1)} +\cdots}{\sum_{n=0}^{\infty}\frac{E^{(1+n)}(L(s)-2)}{n!} \delta^n} \\
            &- \frac{\left(\frac{d E(L(s)-2) + d\log E'(L(s)-2) + \log a}{E(L(s)-1)} +\cdots\right)^2\sum_{n=0}^{\infty}\frac{E^{(2+n)}(L(s)-2)}{n!} \delta^n}{2!\left(\sum_{n=0}^{\infty}\frac{E^{(1+n)}(L(s)-2)}{n!} \delta^n\right)^3} + \cdots 
            \\
            =& \frac{d E(L(s)-2) + d\log E'(L(s)-2) + \log a}{E(L(s)-1)E'(L(s)-2)} + \cdots
    \end{align*}
    where again the ``$\cdots$" has terms with $E(L(s)-1)^n$ in the denominator for $n \ge 2$, which are much smaller than the preceding terms.
    Replacing $a$ by $1$ and $1 + \frac{1}{E(L(s)-1)^{1/3}}$ we have 
    \begin{align*}
        \epsilon_{1,d}(f) < \delta < \epsilon_{1+\frac{1}{E(L(s)-1)^{1/3}},d}(f)
    \end{align*}
    as claimed.
    
    Since $L(s) \in H_{i+1}$ by Lemma \ref{L closure}, we have $f \in H_{i+2}$ by Lemma \ref{E closure}, and thus $\big(E^{(d)}\big)^{-1}(s) \in H_{i+3}$ again by Lemma \ref{E closure}.
\end{proof}

So we have shown that if $s \in (H_i)_{>F}$, then $E^{(d)}(s), L(s) \in H_{i+1}$ and $\left(E^{(d)}\right)^{-1}(s) \in H_{i+3}$.
If $s < 0$, then we define $E^{(d)}(s) = 0$, and if $s \le 1$ then we define $L(s) = \left(E^{(d)}\right)^{-1}(s) = 0$.
If $0 \le s \in H_i$ is bounded in $F$, then there are $r \in F_+$ and $\epsilon \in \mu(H_i)$, the infinitesimals of $H_i$ relative to $F$, such that $s = r + \epsilon$.
Since $F \vDash T_{\mathrm{transexp}}$, $E^{(d)}(r) \in F$ for all $d \in \mathbb{N}$, so we identify $E^{(d)}(s)$ with the series
    $$\sum_{n=0}^{\infty} \frac{E^{(d+n)}(r)}{n!} \epsilon^n \in H_i$$
which is an element of $H_{i}$ by Lemma \ref{R_an closure}.

If $s > 1$, then we can define
\begin{align*}
    L(s) &= \sum_{n=0}^{\infty} \frac{L^{(n)}(r)}{n!} \epsilon^n \in H_i \\
    \left(E^{(d)}\right)^{-1}(s) &=\sum_{n=0}^{\infty} \frac{\left(\left(E^{(d)}\right)^{-1}\right)^{(n)}(r)}{n!} \epsilon^n \in H_i.
\end{align*}
So $M_F$ can be made into a $\mathcal{L}_{\mathrm{transexp}}$-structure and a model of $T_{\mathrm{transexp}}$.

\subsection{A derivation on \texorpdfstring{$M_F$}{MF}}\label{subsection derivation}

There is a derivation $\partial_0$ on $H_0 = F((\tau^{-1}))^{le}$, which can be thought of as differentiation with respect to $\tau$.
Its field of constants is $F$.
We will show that given a derivation $\partial_i$ on $H_{i}$ that we think of as differentiation with respect to $\tau$, we can extend it to a derivation $\partial_{i+1}$ on $H_{i+1}$. 

Let $i>0$.
We first define how $\partial_{i+1}$ acts on generators built from $X_{i}'$.
We will extend it to $D_{X_i'}(F)$, then to $C_{\varnothing,i}$, and then to $C_{\bar{\alpha},i}$ inductively for each finite increasing sequence $\bar{\alpha} \subset \big(H_{i}/\mathrm{Fin} (H_{i})\big)_{>F}$.

\subsubsection{The derivation on \texorpdfstring{$D_{X_i'}(F)$}{DXi'(F)}}
First, let $x_i' = L^{\circ (i+2)}(\tau)$ and let
    $$y_{i}' = \frac{1}{E'\big(L^{\circ (i+1)}(\tau)\big) E'\big(L^{\circ i}(\tau)\big) \cdots E'\big(L(\tau)\big)}$$
which is intended to be the derivative of $x_i'$.
For all $d \in \mathbb{N}$ and $a \in F$, define 
\begin{align*}
    \partial_{i+1} \big(E^{(d)}(x_i'-m)^a\big)
        &:= a E^{(d)}(x_i'-m)^{a-1} \cdot  E^{(d+1)}(x_i'-m) \cdot y_i' \\
    \partial_{i+1} \big(\log E'(x_i'-m)^a\big)
        &:= a \log E'(x_i'-m)^{a-1} \cdot \frac{E''(x_i'-m)}{E'(x_i'-m)}\cdot y_i'
\end{align*}
which are elements of $C_{\{\alpha_i < \cdots < \alpha_{1}\},i}$ where $\alpha_j$ is such that $E'\big(L^{\circ j}(\tau)\big) \in X_{\alpha_j} \subset \alpha_j$ for $j=1,\dots,i$.
(Recall that $E'\big(L^{\circ (i+1)}(\tau)-m\big) \in C_{\varnothing,i}$.)
The indices in the sequence of $\alpha_i <\alpha_{i-1}<\dots,<\alpha_1$ are decreasing because $E'(L^{\circ i}(\tau)) < \cdots < E'(L(\tau))$.
Extend $\partial_{i+1}$ to products so that it satisfies the Leibniz rule.
Extend $\partial_{i+1}$ to sums in $K_{X_i'-m,0}$ by
\begin{align*}
    \partial_{i+1} \left(\sum_{M \in \Gamma_{X_i'-m,0}} c_M M\right) := \sum_{M \in \Gamma_{X_i'-m,0}} c_M\partial_{i+1}(M).
\end{align*}
Extend the derivation to monomials with exp by defining $\partial_{i+1}(e(a)) = e(a)\partial_{i+1}(a)$.
We must show that $\partial_{i+1}$ maps $D_{X_i'}(F)$ to $H_{i+1}$ and that it is well defined.

\begin{lemma}\label{derivative maps to target}
    For each $s \in K_{X_i'-m}$, we have $\partial_{i+1}(s) \in C_{\{\alpha_i < \cdots < \alpha_{1}\},i}$ where $\alpha_j$ is such that $E'\big(L^{\circ j}(\tau)\big) \in X_{\alpha_j} \subset \alpha_j$ for $j=1,\dots,i$.
\end{lemma}

\begin{proof}
    If $s \in K_{X_i'-m}$, then every monomial of $\partial_{i+1}(s)$ is a product of a monomial of $K_{X_i'-m}$ with $y_i'$.
    We will show that if $s \in K_{X_i'-m}$, then $\frac{\partial_{i+1}(s)}{y_i'}$ is a valid sum in $K_{X_i'-m}$ and thus an element of $C_{\varnothing,i}$.
    
    Note that the sum of exponents of generators of $\frac{\partial_{i+1}\Big(E^{(d)}\big(x_i'-m\big)^a\Big)}{y_i'}$ is still $a$.
    Using this, we make two observations:
    \begin{enumerate}
        \item If $g \in \Gamma_{X_i'-m,0}^{small}$, then $\frac{\partial_{i+1}(g)}{y_i'} \in F[\Gamma_{X_i'-m,0}^{small}]$.
        
        \item If $t \in T_{X_i'-m}$, then $\frac{\partial_{i+1}(t)}{y_i'} \in T_{X_i'-m} \subset F[\Gamma_{X_i'-m,0}^{small}]$.
    \end{enumerate}
    So if $M = e_{X_i'-m}(t) g$ is in a coset $w \in \Gamma_{X_i'-m,0}/\Gamma_{X_i'-m,0}^{small}$, then
    \begin{align*}
        \frac{\partial_{i+1}(M)}{y_i'} = e_{X_i'-m}(t)\left(\frac{\partial_{i+1}(t)}{y_i'} g + \frac{\partial_{i+1}(g)}{y_i'}\right)
    \end{align*}
    and $\mathrm{Supp} \left(\frac{\partial_{i+1}(M)}{y_i'}\right)$ is a finite subset of $w$.
    So if $s \in K_{X_i'-m,0}$, then $\frac{\partial_{i+1}(s)}{y_i'}$ is a valid sum in $K_{X_i'-m,0}$.
    
    Now assume that for all $l=0,\dots,n$, if $s \in K_{X_i'-m,l}$, then $\frac{\partial_{i+1}(s)}{y_i'} \in K_{X_i'-m,l}$.
    Let $s = \sum c_ae_{X_i'-m}(a) \in K_{X_i'-m,n+1}$ where $a \in A_{X_i'-m,n}$ and $c_a \in K_{X_i'-m,n}$.
    Then
    \begin{align*}
        \frac{\partial_{i+1}(s)}{y_i'} 
            = \sum_{a \in A_{X_i'-m,n}} \left(\frac{\partial_{i+1}(c_a)}{y_i'} + c_a \frac{\partial_{i+1}(a)}{y_i'}\right)e_{X_i'-m}(a).
    \end{align*}
    By assumption, $\frac{\partial_{i+1}(c_a)}{y_i'} + c_a \frac{\partial_{i+1}(a)}{y_i'} \in K_{X_i'-m,n}$, so $\frac{\partial_{i+1}(s)}{y_i'} \in K_{X_i'-m,n+1}$.
\end{proof}

\begin{lemma}\label{derivative well defined}
    Let $s \in K_{X_i'-m}$. Then 
        $$\varphi_m\left(\frac{\partial_{i+1}(s)}{y_i'}\right) = \frac{\partial_{i+1}(\varphi_m(s))}{y_i'}.$$
\end{lemma}

\begin{proof}
    It suffices to prove the maps commute on generators built from $X_i'$.
    We can compute
    \begin{align*}
        \varphi_m\left(\frac{\partial_{i+1}\big(E^{(d)}(x_i')^a\big)}{y_i'}\right) 
            =& \varphi_m\left(aE^{(d)}(x_i')^{a-1}E^{(d+1)}(x_i')\right) \\
            =& ae^{aE(x_i'-1)}E'(x_i'-1)^{d(a-1)} \sum_{n=0}^{\infty} \binom{a-1}{n} \left(\frac{B_d(x_i'-1)}{E'(x_i'-1)^d}-1\right)^n\\
                & \cdot B_{d+1}(x_i'-1) \\
        \frac{\partial_{i+1}(\varphi_m(E^{(d)}(x_i')^a))}{y_i'} 
            =& \frac{\partial_{i+1}\left(e^{aE(x_i'-1)}E'(x_i'-1)^{da} \sum_{n=0}^{\infty} \binom{a}{n} \left(\frac{B_d(x_i'-1)}{E'(x_i'-1)^d}-1\right)^n\right)}{y_i'} \\
            =& e^{aE(x_i'-1)}\Bigg(aE'(x_i'-1)\cdot E'(x_i'-1)^{da}\sum_{n=0}^{\infty} \binom{a}{n} \left(\frac{B_d(x_i'-1)}{E'(x_i'-1)^d}-1\right)^n \\
                &+ daE'(x_i'-1)^{da-1} \sum_{n=0}^{\infty} \binom{a}{n} \left(\frac{B_d(x_i'-1)}{E'(x_i'-1)^d}-1\right)^n \\
                &+ \bigg(E'(x_i'-1)^{da}\sum_{n=1}^{\infty} n \binom{a}{n} \left(\frac{B_d(x_i'-1)}{E'(x_i'-1)^d}-1\right)^{n-1} \\
                &\cdot \frac{B_d(x_i'-1)'E'(x_i'-1)^d-dB_d(x_i'-1)E'(x_i'-1)^{d-1}}{E'(x_i'-1)^{2d}}\bigg)\Bigg).
    \end{align*}
    To show these two expressions are equal, observe first that $e^{aE(x_i'-1)}$ appears in all monomials of both, so we may divide it out.
    Next, we will factor 
        $$a\varphi_m\big(E^{(d)}(x_i')^{a-1}\big) = aE'(x_i'-1)^{d(a-1)}\sum_{n=0}^{\infty} \binom{a-1}{n} \left(\frac{B_d(x_i'-1)}{E'(x_i'-1)^d}-1\right)^n$$
    out of the second expression, and show that what remains is equal to $B_{d+1}(x_i'-1)$.
    Note
    \begin{align*}
        \sum_{n=0}^{\infty} \binom{a}{n} \left(\frac{B_d(x_i'-1)}{E'(x_i'-1)^d}-1\right)^n 
            &= \left(\sum_{n=0}^{\infty} \binom{a-1}{n} \left(\frac{B_d(x_i'-1)}{E'(x_i'-1)^d}-1\right)^n\right)\cdot \frac{B_d(x_i'-1)}{E'(x_i'-1)^d} \\
        \sum_{n=1}^{\infty} n\binom{a}{n} \left(\frac{B_d(x_i'-1)}{E'(x_i'-1)^d}-1\right)^{n-1}
            &= a\sum_{n=1}^{\infty} \left(\binom{a}{n} -\binom{a-1}{n}\right) \left(\frac{B_d(x_i'-1)}{E'(x_i'-1)^d}-1\right)^{n-1} \\
            &= a\sum_{n=0}^{\infty} \left(\binom{a}{n+1} -\binom{a-1}{n+1}\right) \left(\frac{B_d(x_i'-1)}{E'(x_i'-1)^d}-1\right)^{n} \\
            &= a\sum_{n=0}^{\infty} \binom{a-1}{n} \left(\frac{B_d(x_i'-1)}{E'(x_i'-1)^d}-1\right)^{n}.
    \end{align*}
    So when we factor out $a\varphi_m\big(E^{(d)}(x_i')^{a-1}\big)$ of the second expression we are left with
    \begin{align*}
        \frac{\frac{\partial_{i+1}(\varphi_m(s))}{y_i'}}{a\varphi_m\big(E^{(d)}(x_i')^{a-1}\big)} 
            =& E'(x_i'-1)B_d(x_i'-1) + d\frac{B_d(x_i'-1)}{E'(x_i'-1)} \\
                &+ E'(x_i'-1)^d\left(\frac{B_d(x_i'-1)'}{E'(x_i'-1)^d}-d\frac{B_d(x_i'-1)}{E'(x_i'-1)^{d+1}}\right) \\
            =& E'(x_i'-1)B_d(x_i'-1) + B_d(x_i'-1)' \\
            =& B_{d+1}(x_i'-1).
    \end{align*}
    The computation showing the maps commute on generators of the form $\log E'(x_i')^a$ is very similar.
\end{proof}

Lemma \ref{derivative well defined} shows that $\partial_{i+1}$ is well defined on $D_{X_i'}(F)$.

\subsubsection{The derivation on $C_{\varnothing,i}$}

Let $x_i = L^{\circ (i+1)}(\tau)$ and
\begin{align*}
    y_i &= \frac{1}{E'\big(L^{\circ i}(\tau)\big) E'\big(L^{\circ (i-1)}(\tau)\big) \cdots E'\big(L(\tau)\big)}
\end{align*}
which is intended to be the derivative of $x_i$.
For all $d \in \mathbb{N}$ and $a \in D_{X_i'}(F)$, define 
\begin{align*}
    \partial_{i+1} \big(E^{(d)}(x_i-m)^a\big) =
        & E^{(d)}(x_i-m)^{a}\bigg(\partial_{i+1}(a)\log\big(E^{(d)}(x_i-m)\big)+ a\frac{E^{(d+1)}(x_i-m)}{E^{(d)}(x_i-m)}y_{i}\bigg) 
        \\
    \partial_{i+1} \big(\log E'(x_i-m)^a\big) =
        & \log E'(x_i-m)^{a}
        \\
        &\bigg(\partial_{i+1}(a)\log\big(\log E'(x_i-m)\big)+ a\frac{E''(x_i-m)}{\log E'(x_i-m)E'(x_i-m)}y_{i}\bigg)
\end{align*}
which are elements of $C_{\{\alpha_i < \cdots < \alpha_{1}\},i}$ where $\alpha_j$ is such that $E'\big(L^{\circ j}(\tau)\big) \in X_{\alpha_j} \subset \alpha_j$.
Extend $\partial_{i+1}$ to products so that it satisfies the Leibniz rule.
Extend $\partial_{i+1}$ to sums in $K_{X_i-m,0}$ by
\begin{align*}
    \partial_{i+1} \left(\sum_{M \in \Gamma_{X_i-m,0}} c_M M\right) 
        = \sum_{M \in \Gamma_{X_i-m,0}} \partial_{i+1}(c_M)M + c_M\partial_{i+1}(M).
\end{align*}
Extend the derivation to monomials with exp by defining $\partial_{i+1}(e(a)) = e(a)\partial_{i+1}(a)$.
We must show that $\partial_{i+1}$ maps $C_{\varnothing,i}$ to $H_{i+1}$ and that it is well defined.

\begin{rmk}\label{split derivative of s}
    Compute that
    \begin{align*}
        \log\big(E^{(d)}(x_i-m)\big) = 
            &E(x_i-m-1) + d\log E'(x_i-m-1) 
            \\
            &+ \sum_{k=1}^{\infty} \frac{(-1)^{k+1}}{k}\left(\frac{B_d(x_i-m-1)}{E'(x_i-m-1)^d}-1\right)^k\\
        \log\big(\log E'(x_i-m)\big) =
            &E(x_i-m-2) + \sum_{k=1}^{\infty} \frac{(-1)^{k+1}}{k}\left(\frac{\log E'(x_i-m-1)}{E(x_i-m-1)}\right)^k.
    \end{align*}
    So if $s = \sum_{M \in \Gamma_{X_i-m,0}} c_M M$, then every monomial of $\sum_{M \in \Gamma_{X_i-m,0}} c_M\partial_{i+1}(M)$ is a product of a monomial of $K_{X_i-m-2}$ with either $y_i'$ or $y_i$, i.e., we can split this sum into
        $$\sum_{M \in \Gamma_{X_i-m,0}} c_M\partial_{i+1}(M) = s_0y_i' + s_1y_i$$
    with $\mathrm{Supp}(s_0),\mathrm{Supp}(s_1) \subset \Gamma_{X_i-m-2,2}$.
        
    Inductively, if $s = \sum_{a \in A_{X_i-m,n}} c_a e_{X_i-m}(a) \in K_{X_i-m,n+1}$, then we have
        $$\sum_{a \in A_{X_i-m,n}} c_a\partial_{i+1}\big(e_{X_i-m}(a)\big) = \sum_{a \in A_{X_i-m,n}} c_ae_{X_i-m}(a)\partial_{i+1}(a).$$
    Since we can split each $\partial_{i+1}(a)$, we can also split 
        $$\sum_{a \in A_{X_i-m,n}} c_a\partial_{i+1}\big(e_{X_i-m}(a)\big) = s_0y_i' + s_1y_i$$
    with $\mathrm{Supp}(s_0),\mathrm{Supp}(s_1) \subset \Gamma_{X_i-m-2,n+3}$.
\end{rmk}

\begin{lemma}\label{derivative maps to target 2}
    For each $s \in K_{X_i-m}$, we have $\partial_{i+1}(s) \in C_{\{\alpha_i < \cdots < \alpha_{1}\},i}$, where $\alpha_j$ is such that $E'\big(L^{\circ j}(\tau)\big) \in X_{\alpha_j} \subset \alpha_j$.
\end{lemma}

\begin{proof}
    If $s = \sum_{M \in \Gamma_{X_i-m,0}} c_M M \in K_{X_i-m,0}$, then so is $\sum_{M \in \Gamma_{X_i-m,0}} \frac{\partial_{i+1}(c_M)}{y_i'}M$ since Lemmas \ref{derivative maps to target} and \ref{derivative well defined} show each $\frac{\partial_{i+1}(c_a)}{y_i'} \in D_{X_i'}(F)$, the field of coefficients.
    So we need only worry about 
        $$\sum_{M \in \Gamma_{X_i-m,0}} c_M\partial_{i+1}(M).$$
    Write $\sum_{M \in \Gamma_{X_i-m,0}} c_M\partial_{i+1}(M) = s_0y_i' + s_1y_i$ as in Remark \ref{split derivative of s}.
    To prove the Lemma, we must show $s_0,s_1$ are valid sums in $K_{X_i-m-2}$.
    
    By the exact same argument as Lemma \ref{derivative maps to target}, the preimage of $s_1$ is a valid sum in $K_{X_i-m}$, and thus $s_1$ is a valid sum in $K_{X_i-m-2}$. 
    
    Now we show $s_0$ is a valid sum in $K_{X_i-m-2}$.
    Let $w \in \Gamma_{X_i-m,0}/\Gamma_{X_i-m,0}^{small}$ be a coset with representatives appearing in $s$.
    Since $s|_w$ is finite, there is some largest derivative $d$ appearing in $s|_w$.
    Thus the part of $\partial_{i+1}\big(s|_w\big)$ that contributes to $s_0$ can be split into 
        $$t_{-1}\log \log E'(x_i-m) + t_0 \log E(x_i-m) + \cdots + t_d \log E^{(d)}(x_i-m)$$
    This expression gives a valid sum in $K_{X_i-m-2}$ because each $t_j$ is a finite sum.
    Since $\varphi_m(s|_v) > \varphi_m(s|_w)$ if $v > w$, the whole of $s_0$ is a valid sum in $K_{X_i-m-2}$.
    
    Now assume the lemma holds for all $s \in K_{X_i-m,l}$ for $l=0,\dots,n$.
    Suppose 
        $$s = \sum_{a \in A_{X_i-m,n}} c_ae_{X_i-m}(a) \in K_{X_i-m,n+1}$$
    where $a,c_a \in A_{X_i-m,n}$.
    Just as above, $\sum_{a \in A_{X_i-m,n}} \frac{\partial_{i+1}(c_a)}{y_i'}e_{X_i-m}(a) \in K_{X_i-m,n+1}$ by Lemmas \ref{derivative maps to target} and \ref{derivative well defined}.
    By Remark \ref{split derivative of s}, write 
        $$\sum_{a \in A_{X_i-m,n}} c_a\partial_{i+1}\big(e_{X_i-m}(a)\big) = s_0y_i' + s_1y_i.$$
    Again, $s_1 \in K_{X_i-m-2}$ by the same argument as Lemma \ref{derivative maps to target}.
    And
        $$s_0 = \sum_{a \in A_{X_i-m,n}} (\varphi_{m+1} \circ \varphi_{m})\big(e_{X_i-m}(a)\big)a_0$$
    where $a_0$ is such that $\partial_{i+1}(a) = a_0y_i' + a_1y_i$.
    Since $a_0 \in K_{X_i-m-2,n+2}$, $s_0$ is a valid sum in $K_{X_i-m-2,n+3}$.
\end{proof}

\begin{lemma}\label{derivative well defined logs}
    Let $s \in K_{X_i-m}$. 
    Write 
    \begin{align*}
        \partial_{i+1}(s) &= s_0y_i' + s_1y_i \\
        \partial_{i+1}\big(\varphi_m(s)\big) &= t_0y_i' + t_1y_i
    \end{align*}
    as in Remark \ref{split derivative of s}, where $s_0,s_1,t_0,t_1 \in C_{\varnothing,i}$ by Lemma \ref{derivative maps to target 2}.
    Then $\varphi_m(s_0) = t_0$.
\end{lemma}

\begin{proof}
    It suffices to show that the maps commute on generators built from $X_i$.
    First let $s = E^{(d)}(x_i)^a$.
    Let $Y = \left(\frac{B_d(x_i-1)}{E'(x_i-1)^d}-1\right)$.
    Then
    \begin{align*}
        \varphi_m(s_0)y_i' 
            =& \varphi_m\big(E^{(d)}(x_i)^a\log E^{(d)}(x_i)\big) \partial_{i+1}(a)\\
            =& e_{X_i}(aE(x_i-1)) E'(x_i-1)^{da} \partial_{i+1}(a)
            \\&\sum_{n=0}^{\infty}\binom{a}{n} Y^n\left(E(x_i-1) + dE'(x_i-1) + \sum_{k=1}^{\infty}\frac{(-1)^{k+1}}{k}Y^k \right) \\
        t_0y_i' =& E(x_i)^a \log E(x_i) \partial_{i+1}(a) \cdot E'(x_i-1)^{da}\sum_{n=0}^{\infty}\binom{a}{n} Y^n \\
            &+ E(x_i)^a E'(x_i-1)^{da}\log E'(x_i-1)\partial_{i+1}(da)\sum_{n=0}^{\infty}\binom{a}{n} Y^n \\
            &+ E(x_i)^aE'(x_i-1) \sum_{n=0}^{\infty}\partial_{i+1}\left(\binom{a}{n}\right) Y^n \\
            =& E(x_i)^aE'(x_i-1)^{da}
            \\
            &\left(\sum_{n=0}^{\infty}\binom{a}{n} Y^n\big(E(x_i-1) + dE'(x_i-1)\big)\partial_{i+1}(a) + \sum_{n=0}^{\infty}\partial_{i+1}\left(\binom{a}{n}\right) Y^n\right)
    \end{align*}
    Matching like terms, all that remains to show is that
    \begin{align*}
        \left(\sum_{n=0}^{\infty}\binom{a}{n} Y^n\right)\left(\sum_{k=1}^{\infty}\frac{(-1)^{k+1}}{k}Y^k\right)\partial_{i+1}(a) = \sum_{n=0}^{\infty}\partial_{i+1}\left(\binom{a}{n}\right) Y^n
    \end{align*}
    which follows from the identity
    \begin{align*}
        \frac{d}{dX}\binom{X}{n} = \sum_{j=0}^{n-1}\frac{(-1)^{n-j-1}}{n-j}\binom{X}{j}
    \end{align*}
    for binomial coefficients.
    The argument for $s = \log E'(x_i)^a$ is very similar.
\end{proof}

To see that $\partial_{i+1}$ is well defined on $C_{\varnothing,i}$, write
\begin{align*}
    \partial_{i+1}(s) &= s_0y_i' + s_1y_i \\
    \partial_{i+1}\big(\varphi_m(s)\big) &= t_0y_i' + t_1y_i
\end{align*}
as in Remark \ref{split derivative of s}, where $s_0,s_1,t_0,t_1 \in C_{\varnothing,i}$ by Lemma \ref{derivative maps to target 2}.
By Lemma \ref{derivative well defined logs}, $\varphi(s_0) = t_0$.
By the same argument as in Lemma \ref{derivative well defined}, $\varphi_m(s_1) = t_1$.
So $\partial_{i+1}$ commutes with the maps $\varphi_m$.

\subsubsection{The derivation on $C_{\bar{\alpha},i}$}

Before extending the derivation, we will associate to every $s \in H_{i+1}$ a finite sequence $\chi(s) \subset H_{i+1}$.
$\chi(s)$ will list all the elements that some $E^{(d)}(\cdot)$ appearing in $s$ is composed with, including instances in the coefficients and exponents of $s$.
\begin{enumerate}
    \item For $s \in D_{X_i'}(F)$, let $\chi(s) := (y_i')$.
    
    \item For $s \in C_{\varnothing,i}\setminus D_{X_i'}(F)$, let $\chi(s) := (y_i',y_i)$.
    
    \item Assume we've defined $\chi(a)$ for all $a \in C_{\bar{\alpha},i}$.
    Assume also that if $s_0,s_1 \in \boldsymbol{k} \in \kappa_{\bar{\alpha},i}$ and neither $s_0$ nor $s_1$ is in any $\boldsymbol{k}_0 \subsetneq \boldsymbol{k}$ for any $\boldsymbol{k}_0 \in \kappa_{\bar{\alpha},i}$, then $\chi(s_0) = \chi(s_1)$.
    
    Let $\beta > \bar{\alpha}$ and $\bar{\alpha}$ minimal such that $s \in C_{\bar{\alpha} \cup \{\beta\},i} \setminus C_{\bar{\alpha},i}$.
    Then $s$ must be in $D_{X}(\boldsymbol{k})$ for a unique smallest finite $X = \{x_1,\dots,x_p\} \subset X_{\beta}$ and unique minimal $\boldsymbol{k} \in \kappa_{\bar{\alpha},i}$.
    Define
        $$\chi(s) := \chi(s_0) \frown (\iota_i(x_1),\dots,\iota_i(x_p))$$
    where $s_0 \in \boldsymbol{k} \in \kappa_{\bar{\alpha},i}$ with $s_0 \not \in \boldsymbol{k}_0 \subsetneq \boldsymbol{k}$ for any $\boldsymbol{k}_0 \in \kappa_{\bar{\alpha},i}$.
\end{enumerate}

Now we extend the derivation.
Let $\bar{\alpha}$ be an increasing sequence in $\big(H_i/\mathrm{Fin} (H_i)\big)_{>F}$, and assume we have defined $\partial_{i+1}$ on $C_{\bar{\alpha},i}$.
Let $\beta > \bar{\alpha}$.
For all $d \in \mathbb{N}$, $x \in X_{\beta}$, and $a \in C_{\bar{\alpha},i}$, define
\begin{align*}
    \partial_{i+1} \big(E^{(d)}(x-m)^a\big)
        &= E^{(d)}(x-m)^{a}\bigg(\partial_{i+1}(a)\log\big(E^{(d)}(x-m)\big)+ a\frac{E^{(d+1)}(x-m)}{E^{(d)}(x-m)}\iota_i\big(\partial_{i}(x)\big)\bigg) 
\end{align*}
and
\begin{multline*}
    \partial_{i+1} \big(\log E'(x-m)^a\big)
        = \log E'(x-m)^{a} \bigg(\partial_{i+1}(a)\log\big(\log E'(x-m)\big) \\+ a\frac{E''(x-m)}{\log E'(x-m)E'(x-m)}\iota_i\big(\partial_{i}(x)\big)\bigg)
\end{multline*}
which are elements of $C_{\bar{\alpha} \cup \{\beta\} \cup \bar{\gamma},i}$ where $\bar{\gamma}$ is such that $\partial_{i+1}(y) \in C_{\bar{\gamma},i}$ for all 
    $$y \in \chi\big(E^{(d)}(x-m)^a\big) = \chi \big(\log E'(x-m)^a\big).$$
Extend $\partial_{i+1}$ to products so that it satisfies the Leibniz rule.
Extend $\partial_{i+1}$ to sums in $K_{X_i-m,0}$ by
\begin{align*}
    \partial_{i+1} \left(\sum_{M \in \Gamma_{X_i-m,0}} c_M M\right) 
        = \sum_{M \in \Gamma_{X_i-m,0}} \partial_{i+1}(c_M)M + c_M\partial_{i+1}(M).
\end{align*}
Extend the derivation to monomials with exp by defining $\partial_{i+1}(e(a)) = e(a)\partial_{i+1}(a)$.
Just like earlier, we must show that $\partial_{i+1}$ maps to $H_{i+1}$ and that it is well defined.

\begin{rmk}\label{split derivative of s inductive}
    Let $s = \sum_{M \in \Gamma_{X-m,0}} c_M M$ where $X = \{x_1,\dots,x_p\} \subset X_{\beta}$.
    Write $\chi(s) = (z_1,\dots,z_q,x_1,\dots,x_p)$.
    For any $x \in X$, we can compute $\log E^{(d)}(x-m)$ and $\log \big(\log E'(x-m)\big)$ just as in Remark \ref{split derivative of s}.
    So every monomial of $\sum_{M \in \Gamma_{X-m,0}} c_M \partial_{i+1}(M)$ is a product of a monomial of $K_{X_{\beta}-m-2}$ with one of $\partial_{i+1}(z_1),\dots,\partial_{i+1}(z_1),\partial_{i+1}(x_1),\dots,\partial_{i+1}(x_p)$.
    So we can split
        $$\sum_{M \in \Gamma_{X-m,0}} c_M \partial_{i+1}(M) = s_1\partial_{i+1}(z_1) + \cdots + s_p\partial_{i+1}(z_q) + s_{p+1}\partial_{i+1}(x_1) + \cdots + s_{p+q}\partial_{i+1}(x_p)$$
    with $\mathrm{Supp}(s_j) \subset \Gamma_{X-m-2}$ for $j=1,\dots,p+q$.
\end{rmk}

\begin{lemma}
    For each $s \in K_{X_{\beta}-m}$, we have 
        $$\partial_{i+1}(s) \in C_{\bar{\alpha} \cup \{\beta\} \cup \bar{\gamma},i}$$
    for some $\bar{\gamma} \subset \big(H_i/\mathrm{Fin} (H_i)\big)_{>F}$.
\end{lemma}

\begin{proof}
    Let $s \in K_{X_{\beta}-m}$, and write $\chi(s) = (z_1,\dots,z_q;x_1,\dots,x_p)$.
    Let $X = \{x_1,\dots,x_p\}$.
    Let $\bar{\gamma}$ be such that $\partial_{i+1}(z_j),\partial_{i+1}(x_l) \in C_{\bar{\alpha} \cup \{\beta\} \cup \bar{\gamma},i}$ for $j=1,\dots,q$, $l = 1,\dots,p$.
    
    First suppose $s = \sum_{M \in \Gamma_{X-m,0}} c_M M \in K_{X_{\beta}-m,0}$.
    Then for each $M$, we can write 
    \begin{align*}
        \partial_{i+1}(c_M) &= c_{M,1}\partial_{i+1}(z_1) + \cdots + c_{M,q}\partial_{i+1}(z_q) 
    \end{align*}
    with $c_{M,1},\dots,c_{M,q} \in C_{\bar{\alpha},i}$.
    Since $s$ is a valid sum, so is $\sum_{M \in \Gamma_{X-m,0}} c_{M,j} M$ for each $j=1,\dots,q$.
    Thus
        $$\sum_{M \in \Gamma_{X-m,0}} \partial_{i+1}(c_M) M = \left(\sum_{M \in \Gamma_{X-m,0}} c_{M,1} M\right) \partial_{i+1}(z_1) + \cdots + \left(\sum_{M \in \Gamma_{X-m,0}} c_{M,q} M\right) \partial_{i+1}(z_q)$$
    is a valid sum in $C_{\bar{\alpha} \cup \{\beta\} \cup \bar{\gamma},i}$.
    
    Following Remark \ref{split derivative of s inductive}, write
        $$\sum_{M \in \Gamma_{X-m,0}} c_M \partial_{i+1}(M) = s_1\partial_{i+1}(z_1) + \cdots + s_p\partial_{i+1}(z_q) + s_{p+1}\partial_{i+1}(x_1) + \cdots + s_{p+q}\partial_{i+1}(x_p).$$
    To show $\sum_{M \in \Gamma_{X-m,0}} c_M \partial_{i+1}(M)$ is a valid sum, we must show that each $s_j \in K_{X_{\beta}-m-2}$.
    For $s_{q+1},\dots,s_{q+p}$, this follows from the argument of Lemma \ref{derivative maps to target}.
    For $s_{p+1},\dots,s_{p+q}$, this follows from the argument of Lemma \ref{derivative maps to target 2}.
    In fact, each $s_j \in K_{X_{\beta}-m-2,2}$ since $s \in K_{X_{\beta}-m,0}$.
    
    Now assume the result holds for all $s \in K_{X_{\beta}-m,l}$ for $l = 1,\dots,n$.
    Suppose
        $$s = \sum_{a \in A_{X-m,n}} c_a e_{X_{\beta}-m}(a) \in K_{X_{\beta}-m,n+1}.$$
    The same argument as above shows 
        $$\sum_{a \in A_{X-m,n}} \partial_{i+1}(c_a) e_{X_{\beta}-m}(a) \in C_{\bar{\alpha} \cup \{\beta\} \cup \bar{\gamma},i}.$$ 
    To show $\sum_{a \in A_{X-m,n}} c_a e_{X_{\beta}-m}(a)\partial_{i+1}(a)$ is a valid sum, note that by induction, we can write 
        $$\partial_{i+1}(a) = a_1\partial_{i+1}(z_1) + \cdots + a_p\partial_{i+1}(z_q) + a_{p+1}\partial_{i+1}(x_1) + \cdots + a_{p+q}\partial_{i+1}(x_p)$$
    with each $a_j \in K_{X_{\beta}-m-2,n+2}$.
    And 
        $$\sum_{a \in A_{X-m,n}} c_a e_{X_{\beta}-m}\big((\varphi_{m+1} \circ \varphi_m)(a)\big)a_j$$
    is a valid sum in $K_{X_{\beta}-m-2,n+3}$ for each $j=1,\dots,p+q$.
\end{proof}

\begin{rmk}
    To show that $\partial_{i+1}$ is well defined on $C_{\bar{\alpha} \cup \{\beta\},i}$, it suffices to show that $\partial_{i+1}$ and $\varphi_m$ commute on generators, i.e., for each $x \in X_{\beta}$ and $a \in C_{\bar{\alpha},i}$, if we write 
        $$\partial_{i+1}\big(\varphi_m(E^{(d)}(x-m))\big) = t_0\partial_{i+1}(a) + t_1 \partial_{i+1}(x)$$
    then
    \begin{align*}
        t_0 &= \varphi_m\left(E^{(d)}(x-m)^a \log E^{(d)}(x-m)\right) \\
        t_1 &= \varphi_m\left(aE^{(d)}(x-m)^a\frac{E^{(d+1)}(x-m)}{E^{(d)}(x-m)}\right)
    \end{align*}
    and similarly for $\log E'(x-m)$.
    The first equality follows from the computations of Lemma \ref{derivative well defined logs}, and the second equality follows from the computations of Lemma \ref{derivative well defined}.
    So $\partial_{i+1}$ is well defined on $C_{\bar{\alpha} \cup \{\beta\},i}$.
\end{rmk}

Having defined $\partial_{i+1}$ on all $C_{\bar{\alpha},i}$, we may extend it to the direct limit $H_{i+1}$ of the $C_{\bar{\alpha},i}$'s.
It is well defined on the direct limit $M_F$ of the $H_i$'s by Lemma \ref{C's form directed system}.
So $M_F$ is a differential field.

\subsection{Ordering of germs of terms at \texorpdfstring{$+\infty$}{+infinity}}

Let $\mathcal{G}$ be the ring of germs at $+\infty$ of functions $f : \mathbb{R} \to \mathbb{R}$, and let $\mathcal{T}(x)$ be the algebra of $\mathcal{L}_{\mathrm{transexp}}$-terms in a single variable $x$ over $\mathbb{R}$. Define $\theta : \mathcal{T}(x) \to \mathcal{G}$ by sending each term $t(x)$ to the germ of the function $x \mapsto t(x)$ at $+\infty$.
Let $M_0$ be the $\mathcal{L}_{\mathrm{transexp}}$-substructure of $M_{\mathbb{R}}$ generated by $\tau$.

\begin{lemma}\label{theta well defined}
    For all $t_1(x),t_2(x) \in \mathcal{T}(x)$, we have $T_{\transexp} \vdash \exists y \big((x>y) \rightarrow (t_1(x)=t_2(x))\big)$ if and only if $M_{\RR} \vDash t_1(\tau) = t_2(\tau)$.
\end{lemma}

\begin{proof}
    If $t_1 \ne t_2$ but $M_{\RR} \vDash t_1(\tau) = t_2(\tau)$, then $t_1(\tau)$ and $t_2(\tau)$ must be identified in some direct limit during the construction of $M_{\RR}$.
    The identification made in each step of the definition of each $\varphi_{m,0}$ and $\iota_i$ is a consequence of $T_{\transexp}$.
    So the identification between $t_1(x)$ and $t_2(x)$ in $M_{\RR}$ directly corresponds to a proof via consequences of $T_{\transexp}$ that $t_1(x) = t_2(x)$ for all large enough $x$.
\end{proof}

\begin{theorem}
    $\theta\big(\mathcal{T}(x)\big)$ is totally ordered.
\end{theorem}

\begin{proof}
    Let $M_0$ be the $\mathcal{L}_{\mathrm{transexp}}$-substructure of $M_{\mathbb{R}}$ generated by $\tau$.
    Define $\psi : M_0 \to \mathcal{G}$ as follows:
    \begin{enumerate}
        \item $\psi(\tau)$ is the germ of the identity function.
        
        \item If $\psi(s_i) = g_i \in \mathcal{G}$ for $i=1,\dots,n$ and $\tilde{f} \in \mathcal{L}_{\mathrm{transexp}}$ is an $n$-ary function symbol corresponding to the function $f$, then $\psi(\tilde{f}(s_1,\dots,s_n)) = f(g_1,\dots,g_n)$.
    \end{enumerate}
    $\psi$ is a well defined map by Lemma \ref{theta well defined}, and its image in $\mathcal{G}$ coincides with $\theta(\mathcal{T}(x))$.
    Since $M_0$ is a field, $\psi$ is injective.
    Furthermore, $\psi(M_0)$ is formally real, and since every $s > 0$ in $M_0$ is a square, the order on $\psi(M_0)$ arising from $M_0$ is unique.
    So $\theta\big(\mathcal{T}(x)\big) = \psi(M_0)$ is ordered.
\end{proof}

\bibliographystyle{alpha}
\bibliography{STSeries}

\end{document}